# SCALAR CURVATURE, METRIC DEGENERATIONS, AND THE STATIC VACUUM EINSTEIN EQUATIONS ON 3-MANIFOLDS, II.

MICHAEL T. ANDERSON

## CONTENTS

## 0. INTRODUCTION.

Consider the problem of finding a constant curvature metric on a given closed 3-manifold $M$. It is known that essential spheres and tori, (except in the case of flat manifolds), prevent the existence of such metrics, but it is unknown if there are other topological obstructions. We approach this issue by recasting it as a natural variational problem on the space of metrics $\mathbb{M}$ on $M$. The problem then becomes two-fold. First, find a natural curvature functional $I$ on $\mathbb{M}$ whose minima are necessarily constant curvature metrics. Second, try to detect topological information on $M$ from the geometric behavior of minimizing sequences for $I$.

First, which curvature functional $I$ should one choose that is best suited for these tasks? In dimension 3, the full curvature tensor $R$ is determined algebraically by the Ricci curvature $r$, and so functionals involving $R$ or $r$ have essentially the same properties. However, one does not expect to be able to detect the existence of essential spheres in $M$ from the behavior of minimizing sequences for $\int |r|^2 dV$ for example, or from the geometric structure of limits of minimizing sequences. Choosing other $L^p$ norms of $r$ either does not help or leads to severe analytic problems. These difficulties are explained in detail in [An3, §7]; we also point out that there exist critical points of $\int |r|^2 dV$ on the space of unit volume metrics $\mathbb{M}_1$ which are not constant curvature, c.f. [La].

This leads one to consider functionals involving the scalar curvature $s$. (It would also be interesting to consider the Chern-Simons functional in this respect, but we will not do so here). The Einstein-Hilbert action $\mathcal{S}$ was studied in some detail in the predecessor paper [AnI]. Here, we focus on functionals related to the $L^2$ norm of $s$, i.e.

$$\mathcal{S}^2 = (v^{1/3} \int s^2 dV)^{1/2}, \tag{0.1}$$

Partially supported by NSF Grants DMS 9505744, 9802722, 0072591.

where $v$ is the volume and $dV$ is the volume form. The volume factor is needed to make $\mathcal{S}^2$ scale-invariant. The $L^2$ norm is chosen mainly for the simplicity of the Euler-Lagrange equations, but it is crucial that $2 > \frac{3}{2}$, the latter being the critical curvature exponent in dimension 3.

For several, partly more technical, reasons it turns out that there are numerous advantages to consider instead the related functional

$$\mathcal{S}^2_- = (v^{1/3} \int (s^-)^2 dV)^{1/2}, \tag{0.2}$$

where $s^- = \min(s, 0)$. In §1 we discuss in detail the reasons why $\mathcal{S}^2_-$ can be considered as the optimal functional for the problems above and in particular why it behaves better than the Einstein-Hilbert functional $\mathcal{S}$ or $\mathcal{S}^2$ in (0.1). The functional $\mathcal{S}^2_-$ is of course not smooth across regions where $s$ changes sign. This introduces some technical complications, which however are comparatively minor and are dealt with in detail in §3.

It is obvious that $\mathcal{S}^2(g) \geq 0$ and $\mathcal{S}^2_-(g) \geq 0$. Define the Sigma constant $\sigma(M)$ of $M$ to be the supremum of the scalar curvatures of unit volume Yamabe metrics on $M$; this is a topological invariant of $M$. In case

$$\sigma(M) \leq 0, \tag{0.3}$$

we have, c.f. Proposition 3.1,

$$inf_{\mathbb{M}} \mathcal{S}^2 = inf_{\mathbb{M}} \mathcal{S}^2_- = |\sigma(M)|, \tag{0.4}$$

so that a minimizing sequence $\{g_i\}$ for $\mathcal{S}^2$ or $\mathcal{S}^2_-$ should have the same essential characteristics as a maximizing sequence of Yamabe metrics for the Einstein-Hilbert action $\mathcal{S}$. We note that (0.4) is not true when $\sigma(M) > 0$, since any such manifold has scalar-flat metrics so that inf $\mathcal{S}^2 = $ inf $\mathcal{S}^2_- = 0$.

Thus, throughout this paper, we assume (0.3). By a result of Gromov-Lawson [GL2, Thm.8.1], any closed 3-manifold $M$ which has a $K(\pi, 1)$ in its prime decomposition satisfies $\sigma(M) \leq 0$, so that (0.3) is satisfied for a large variety of 3-manifolds. One (elementary) reason for preferring $\mathcal{S}^2_-$ to $\mathcal{S}^2$ is that the *only* critical points of $\mathcal{S}^2_-$ on $\mathbb{M}$ are constant curvature metrics when $\sigma(M) \leq 0$, c.f. Lemma 1.3; this is not known to be the case for $\mathcal{S}^2$.

Having chosen the functional, consider now the behavior of minimizing sequences for $\mathcal{S}^2_-$. Here, the main difficulty is of course that there are many possible choices of such which could have vastly different behaviors. However, there is a basic dichotomy allowing one to make more specific choices. Thus, suppose first that there exists a minimizing sequence $\{g_i\}$ for $\mathcal{S}^2_-$ which does not (curvature) degenerate in the sense that there exists a global bound of the form

$$\int_M |z_{g_i}|^2 dV_{g_i} \leq \Lambda, \tag{0.5}$$

for some constant $\Lambda < \infty$. Here $z$ is the trace-free Ricci curvature, $z = r - \frac{s}{3}g$. Then the limiting behavior of $\{g_i\}$ can be completely determined from an $L^2$ version of the Cheeger-Gromov theory [CG1,2] of convergence and collapse of Riemannian manifolds. In fact, the limiting structure of $\{g_i\}$ essentially induces the geometrization of the 3-manifold $M$ in the sense of Thurston [Th]; this is proved in [An6, §2], c.f. also [An2] for an outline. In particular, one can detect the presence of essential tori in $M$ from the near-limiting behavior of $\{g_i\}$ for $i$ sufficiently large.

Note that the existence of a minimizing sequence $\{g_i\}$ satisfying (0.5) is a topological condition on $M$. Under (0.5) and under the standing assumption (0.3), it is shown in [An6] that $\sigma(M) = 0$ if and only if $M$ is a graph manifold with infinite $\pi_1$. In fact, $\sigma(M)$ is closely related with the hyperbolic volume of the hyperbolic pieces in the geometric decomposition of $M$.

However, the topological meaning of the bound (0.5) is not at all clear. It is not known if it holds for *any* closed 3-manifold $M$, (except graph manifolds). Similarly, there is not a single example of a 3-manifold known to have $\sigma(M) < 0$, c.f. [An2, §4], [Gr2, Rmk. 3.B], so that this invariant could possibly be trivial among manifolds satisfying $\sigma(M) \leq 0$.

Thus, the main purpose of this paper is to begin the analysis of the situation where (0.5) does not hold, i.e. study the behavior of minimizing sequences $\{g_i\}$ for $\mathcal{S}^2_-$ which necessarily degenerate in the sense that

$$\int_M |z_{g_i}|^2 dV_{g_i} \to \infty, \quad \text{as} \quad i \to \infty. \tag{0.6}$$

Again, one does not expect to be able to understand the behavior of an arbitrary minimizing sequence $\{g_i\}$ satisfying (0.6); all kinds of complicated behavior could occur, c.f. also §1(III). Thus, from the outset, we work with a preferred or optimal choice of minimizing sequence, which in a precise sense has the least possible amount of (non-scalar) curvature in $L^2$.

To do this, consider the following scale-invariant perturbation of $\mathcal{S}^2_-$ :

$$I_{\varepsilon}{}^- = \varepsilon v^{1/3} \int_M |z|^2 dV + \left(v^{1/3} \int_M (s^-)^2 dV\right)^{1/2}. \tag{0.7}$$

We assume that $\varepsilon$ is a free positive parameter, so that $I_{\varepsilon}{}^-$ is a perturbation of the $L^2$ norm of the scalar curvature $s^-$ in the direction of the $L^2$ norm of the curvature $z$. If $g_{\varepsilon}$ is a unit volume minimizer for $I_{\varepsilon}{}^-$, then as $\varepsilon \to 0$, the metrics $\{g_{\varepsilon}\}$ form a mimimizing family for $\mathcal{S}^2_-$, and $g_{\varepsilon}$ has smaller value for the $L^2$ norm of $z$ than any other unit volume metric on $M$ with the same value of $\mathcal{S}^2_-$, c.f. §3.1. The choice of the particular perturbing functional, i.e. the $L^2$ norm of $z$, is also optimal in many respects, c.f. again §1 for further discussion.

By way of background, it is perhaps useful to recall that the solution to the Yamabe problem also requires the use of a suitable perturbation of the Einstein-Hilbert functional. In this case, the perturbation is chosen to make use of the compactness of the Sobolev embedding $L^{1,2} \subset L^p, p < 2n/(n-2)$. In this manner, one obtains a preferred minimizing sequence of metrics $\{g_p\}$ for the Yamabe problem and the resolution of the problem is obtained by understanding the behavior of these metrics as $p$ converges to the critical exponent $2n/(n-2)$.

The perturbation in the direction of the $L^2$ norm of curvature serves a similar purpose. The space of unit volume metrics on $M$ which have a uniform bound on the $L^2$ norm of curvature is essentially compact; more precisely, while it is non-compact, the allowable degenerations are well-understood in terms of the collapse theory of Cheeger-Gromov [CG1,2]. In particular, we prove in Theorem 3.9 that minimizers $g_{\varepsilon}$ for $I_{\varepsilon}{}^-$ exist, (when $M$ is not a graph manifold), and that they have good regularity and completeness properties.

Given this, we then need to understand the behavior of $\{g_{\varepsilon}\}$ as $\varepsilon \to 0$, or some sequence $\varepsilon = \varepsilon_i \to 0$. The behavior (0.6) implies that $\{g_{\varepsilon}\}$ degenerates on any sequence $\varepsilon = \varepsilon_i \to 0$ so that there is no smooth convergence to a limit metric $g_o$. The second problem above then becomes one of understanding if such degenerations are related in a simple way to the topology of $M$, in particular whether they are caused by the presence of essential 2-spheres or not.

To study this degeneration in general, we blow-up or rescale the metrics $g_{\varepsilon}$ at points $x_{\varepsilon}$ by a natural measure of the size of the curvature near $x_{\varepsilon}$. Thus, consider the rescalings

$$g_{\varepsilon}' = \rho(x_{\varepsilon})^{-2} \cdot g_{\varepsilon}, \tag{0.8}$$

where $\rho(x_{\varepsilon})$ is the $L^2$ curvature radius at $x_{\varepsilon}$, c.f. (2.1) for the exact definition. Very roughly, $\rho(x)$ is the largest radius such that the function $\rho^2(x)|R|(y)$ has bounded $L^2$ average over $y \in B_x(\rho(x))$ and so is a measure of the degree of curvature concentration in $L^2$ near $x$. We assume that $\rho(x_{\varepsilon}) \to 0$ as $\varepsilon = \varepsilon_i \to 0$, meaning that the curvature is blowing up, (i.e. diverging to $+\infty$), in $L^2$, in smaller and smaller neighborhoods of $x_\varepsilon$. The radius $\rho$ scales as a distance, so that w.r.t. $g'_\varepsilon$,

$$\rho'(x_\varepsilon) = 1. \tag{0.9}$$

Since (0.9) holds, we may pass to a limit metric $g'$ at least locally, i.e. within $B'_{x_\varepsilon}(1)$, passing to subsequences if necessary; this requires an unwrapping of any collapse in general, c.f. §2. These local limits, sometimes called geometric limits, model the small-scale degeneration of the metrics $g_\varepsilon$ near $x_\varepsilon$ in case one has (at least) strong $L^{2,2}$ convergence of $g'_\varepsilon$ to $g'$. The limits may or may not extend to complete Riemannian manifolds; this depends strongly on the choice of the base points $x_\varepsilon$.

One of the main purposes of this paper is to establish the existence of, and strong convergence to, such blow-up limits, classify all of their possible forms, and prove the existence of complete limits $(N, g')$. This is accomplished in §4, but we will not detail all of these results here, c.f. Theorem 4.2, Proposition 4.4 and Theorem 4.9 for example. By way of contrast, we note that the non-trivial blow-up limits constructed in [AnI] are *never* complete.

The form of the local blow-up limits depends crucially on the rate at which the curvature blows up to $\infty$ near $x_\varepsilon$, as $\varepsilon \to 0$. This is measured by the ratio $\varepsilon/\rho^2(x_\varepsilon)$. Both terms here go to 0 as $\varepsilon \to 0$ and so apriori, after passing to subsequences, one may have $\varepsilon/\rho^2(x_\varepsilon)$ either diverges to $+\infty$, converges to a positive number $\alpha$, or converges to 0. Each case here gives rise to blow-up limits with different geometric properties. In Theorem 5.2, we prove (essentially) that the first case cannot occur, i.e. $\varepsilon/\rho^2 \leq K$, for some $K$ independent of $x_\varepsilon$. Thus, there is a maximal or absolute scale of curvature blow-up, depending only on $M$. In the opposite direction, in Theorem 5.5 it is proved that there are base points $y_\varepsilon \in M$ satisfying the opposite inequality, i.e. $\varepsilon/\rho^2(y_\varepsilon) \to \alpha > 0$. Not all base points satisfy this property, but those that do are the most significant.

We summarize some of the discussion above in the following result, one of the main conclusions of the paper; it should be compared with Theorem A of [AnI].

**Theorem B.** *Suppose $\sigma(M) \leq 0$, and let $\{g_\varepsilon\}$ be a sequence of minimizers for $I_\varepsilon^-$, for $\varepsilon = \varepsilon_i \to 0$. Suppose the sequence $\{g_\varepsilon\}$ degenerates, in the sense that*

$$\int |z_{g_\varepsilon}|^2 dV_{g_\varepsilon} \to \infty, \quad as \ \ \varepsilon \to 0. \tag{0.10}$$

*Then there exists a sequence of points $\{y_\varepsilon\} \in (M, g_\varepsilon)$, and scale-factors $\rho(y_\varepsilon) \to 0$, such that the blow-up metrics*

$$g'_\varepsilon = \rho(y_\varepsilon)^{-2} \cdot g_\varepsilon, \tag{0.11}$$

*based at $y_\varepsilon$, have a subsequence converging in the strong $L^{2,2}$ topology to a limit $(N, g', y)$.*

*The limit $(N, g')$ is a complete non-flat $C^{2,\alpha} \cap L^{3,p}$ Riemannian manifold, for any $\alpha < 1$, $p < \infty$, with uniformly bounded curvature, and non-negative scalar curvature.*

*Further, $(N, g')$ minimizes the $L^2$ norm of the curvature $z$ over all metrics $\bar{g}$ of non-negative scalar curvature on $N$, satisfying $vol_{\bar{g}}K \leq vol_g K$ and $\bar{g}|_{N \setminus K} = g'|_{N \setminus K}$, for some (arbitrary) compact set $K \subset N$.*

*The metric $g'$ satisfies the following Euler-Lagrange equations, called* the $\mathcal{Z}_c^2$ equations:

$$\alpha \nabla \mathcal{Z}^2 + L^*(\tau) = 0, \tag{0.12}$$

$$\Delta(\tau + \frac{1}{12}\alpha s) = -\frac{1}{4}\alpha |z|^2,$$

*for some constant $\alpha > 0$. Here $L^*$ is the adjoint of the linearization of the scalar curvature, c.f. (1.6), and $\nabla \mathcal{Z}^2$ is the gradient of the $L^2$ norm of $z$, c.f. (3.8).*

*The potential function $\tau$ is non-positive and locally Lipschitz, the scalar curvature $s$ is non-negative and locally Lipschitz and supp $s$ is disjoint from supp $\tau$, in the sense that*

$$s \cdot \tau \equiv 0. \tag{0.13}$$

*The metric $g'$ is smooth in any region where $\tau < 0$ or $s > 0$ and the junction set $\Sigma = \{\tau = 0\} = \{s = 0\}$ is a $C^{1,\beta}$ smooth surface in $N$.*

*Any smooth compact domain in the manifold $N$ is (naturally) smoothly embedded as a compact submanifold of $M$.*

The potential function $\tau$ plays a central role throughout the paper. Roughly it is the *scale-invariant* limit of the $L^2$ average of the scalar curvature $s_\varepsilon^-$ of the 'base' metrics $g_\varepsilon$, c.f. §4.1. The fact that the metric $g'$ is not $C^\infty$ is directly related to the use of the cutoff $s^-$ in $I_\varepsilon^-$, giving rise to the junction $\Sigma$. (Across the surface $\Sigma$ where $g'$ is not $C^\infty$ smooth, (0.12) is to be understood as holding weakly).

We point out that it is part of the conclusion of Theorem B that the blow-up sequence $(M, g'_\varepsilon, y_\varepsilon)$ does not collapse in the sense of Cheeger-Gromov.

Theorem B describes the degeneration of the metrics $(M, g_\varepsilon)$ near certain preferred base points $y_\varepsilon$ in terms of the global limit $(N, g')$. Such limits are quite rigid, in that they must satisfy a number of rather strong geometric conditions. However, in order to obtain useful topological information about $M$, one needs to understand the global topological and geometric structure of these limits $(N, g')$. In this respect, we prove the following result in §7.

**Theorem C.** *Let $(N, g)$ be a complete $\mathcal{Z}_c^2$ solution, i.e. $(N, g)$ satisfies the data in (0.12) through (0.13). Suppose that there exists a compact set $K \subset N$ and a constant $\omega_o < 0$ such that the potential function*

$$\omega = \tau + \frac{\alpha}{12} s : N \to \mathbb{R}, \tag{0.14}$$

*satisfies*

$$\omega \leq \omega_o < 0, \tag{0.15}$$

*on $N \setminus K$, and that the level sets of $\omega$ in $N$ are compact.*

*Then $N$ is an open 3-manifold, topologically of the form*

$$N = P \#_1^q (\mathbb{R}^3), \tag{0.16}$$

*where $P$ is a closed 3-manifold, (possibly empty), admitting a metric of positive scalar curvature, and $1 \leq q < \infty$.*

*Each end $E = E_k$ of $N$ is asymptotically flat, in the sense it is diffeomorphic to $\mathbb{R}^3 \setminus B$, for some ball $B \subset \mathbb{R}^3$, and the metric $g$ in a suitable chart for $E$ has the form*

$$g_{ij} = (1 + \frac{2m}{t})\delta_{ij} + O(t^{-2}), \tag{0.17}$$

*with $|\partial^k g_{ij}| = O(t^{-2-k}), k > 0.$. Here $t(x) = ||x||$ and $m > 0$ is a constant depending only on $(E, g)$. The potential function $\omega$ is bounded on $E$, and has the asymptotic expansion*

$$\omega = \omega_E + \frac{m|\omega_E|}{t} + O(t^{-2}), \tag{0.18}$$

*where $\omega_E < 0$ is a constant, depending on $(E, g)$.*

We point out that all known examples or constructions of complete $\mathcal{Z}_c^2$ solutions, at least on manifolds $N$ of finite topological type, satisfy the conclusions (0.16)-(0.17). (It is not known if there are any examples of infinite topological type).

The main significance of Theorem C is in the precise determination of the asymptotic geometry and topology of the limit $(N, g')$. In particular this result gives the existence of 'canonical' 2-spheres in each asymptotically flat end of $N$, which hence give distinguished 2-spheres in $M$. In [An6, §3], we prove that if $(N, g')$ is any blow-up limit of $\{g_\varepsilon\}$, $\varepsilon = \varepsilon_i \to 0$, containing an asymptotically flat end $E$, then both $N$ and $M$ are reducible. In fact the canonical 2-sphere from the asymptotically flat end $E$ gives an essential $S^2$ in $M$. These results thus give some partial confirmation at this stage that degeneration of the metrics $g_\varepsilon$ is related to the existence of essential 2-spheres.

In Proposition 7.17, we show that the hypothesis that $\omega$ has compact levels in Theorem C can be replaced by the assumption that the curvature $z$ of $g'$ is sufficiently small outside $K \subset N$. In work to follow, we will discuss the validity of Theorem C in general, i.e. without such assumptions on $\omega$ or $z$.

Finally, it is of interest to construct concrete models or examples of degenerating sequences for $\mathcal{S}^2_-$ that clearly exhibit the various forms of the possible blow-up limits, the differences between convergence and collapse behaviors and the results of Theorems B and C. Such models concretely illustrate that the technical considerations leading to these theoretical conclusions actually arise in very simple situations. Such models are constructed in §6. We recommend the reader to examine these constructions, especially as an aid in understanding some of the motivation for the theoretical arguments.

In Appendix A, we discuss briefly some of the similarities and differences between the approach taken in this paper and its predecessor [AnI]. While this paper is a natural continuation of the issues addressed in [AnI], it is logically independent of [AnI] in the sense that, with one exception, none of the results obtained there are required for the work here.

This completes the discussion of the contents of the paper. I'd like to thank Misha Gromov and referees of the paper for their suggestions in improving the manuscript. A special thanks to Vitali Milman for his patience in seeing the paper through the review process.

## 1. Remarks on the Choice of the Curvature Functional.

In this section, we motivate the choice of the curvature functional $\mathcal{S}^2_-$ and its perturbation $I_\varepsilon{}^-$. This discussion is not logically necessary for the work to follow.

Consider in general scalar curvature functionals of the form

$$\mathcal{S}_f = \int f(s)dV, \tag{1.1}$$

where $f : \mathbb{R} \to \mathbb{R}$ is at least Lipschitz continuous and $\mathcal{S}_f$ is restricted to $\mathbb{M}_1$ to normalize the volume. We assume that $f$ is bounded below and that

$$inf\mathcal{S}_f = |\sigma(M)|,$$

as in (0.4). Typical examples are $f(s) = |s|^p$, for $p \geq 3/2$, c.f. Proposition 3.1 below, but there is a broad range of possibilities.

For three basic reasons, it is necessary to make a specific, and so in effect optimal, choice of a minimizing sequence $\{g_i\}$ for $\mathcal{S}_f$. It is worth discussing these reasons briefly.

**(I). (Strong Convergence).** In order to have limits of blow-ups by the curvature, as in (0.8), accurately describe the local degenerations of the sequence, it is necessary that the blow-up sequence converge strongly to the limit. Without such strong convergence, the limits may always be flat for instance and thus bear no relation with the degeneration itself. This requires that the minimizing sequence satisfy certain regularity conditions. In practice, this can be obtained if the minimizing sequence satisfies a suitable family of P.D.E.'s, or possibly by a smoothing process, but cannot be expected to hold for an arbitrary minimizing sequence.

**(II). (Collapse).** If a minimizing sequence $\{g_i\}$ for $\mathcal{S}_f$, or a sequence of its blow-ups, collapses in some region, in the sense of Cheeger-Gromov [CG1,2], then limits of such a sequence will not even exist, (as Riemannian metrics on some 3-manifold). In this situation, it is possible to pass locally to the universal cover to unwrap the collapse, and thus pass to a limit, c.f. §2. However, since collapse is associated to loss of volume and the functional $\mathcal{S}_f$ is an integral, properties of an arbitrary minimizing sequence are not likely to be preserved when passing to local covers. Thus for example, the methods of [AnI], (which deal with general sequences), are not able to handle collapse situations, c.f. [AnI, §4.1].

**(III).(Topological Degenerations).** It does not seem possible to attach any meaning to the limit of an arbitrary minimizing sequence $\{g_i\}$ for $\mathcal{S}_f$. The functional $\mathcal{S}_f$ for any $f$ is too weak to bound or control the geometry and topology of $(M, g_i)$. For example, there are simple examples of sequences which degenerate in a topologically invisible way, (i.e. the degeneration occurs in 3-balls), on arbitrarily dense sets in $M$, c.f. [AnI, Rmk 6.1(i),(ii)].

To rule out such 'superfluous' degenerations, one would like to find preferred minimizing sequences whose degeneration is caused solely by the underlying topology of $M$, as discussed in §0. A common and useful method to prove that such an optimal sequence does not degenerate for topologically trivial reasons is to construct comparison sequences which would have better properties than the optimal sequence in such situations.

As we will see from Remark 3.2, metrics $g_\varepsilon = g_{\varepsilon,f}$ which minimize the perturbed functional

$$I_{\varepsilon,f} = \varepsilon v^{1/3} \int |z|^2 + \mathcal{S}_f, \tag{1.2}$$

form a minimizing sequence or family of $\mathcal{S}_f$, as $\varepsilon \to 0$, at least for many reasonable choices of functions $f$ as following (1.1). However, they differ markedly in their ability to resolve the issues (I),(III), (although much less so for (II)).

We begin with some general remarks on the functionals $\mathcal{S}_f$.

**Lemma 1.1.** *Suppose $f : \mathbb{R} \to \mathbb{R}$ is a $C^1$ function. Then $\mathcal{S}_f : \bar{\mathbb{M}}_1 \to \mathbb{R}$ is a differentiable functional on the $L^{2,2}$ completion $\bar{\mathbb{M}}_1$ of $\mathbb{M}_1$.*

*Proof.* If $f$ is a $C^1$ function on $\mathbb{R}$, then it is easy to see that $\mathcal{S}_f$ is (Frechet) differentiable on $\mathbb{M}$. Thus, if $h$ is a smooth symmetric form on $M$, and $g_t = g + th$, then

$$d(\mathcal{S}_f)_g(h) = \frac{d}{dt}\mathcal{S}_f(g + th)|_{t=0} = \int f'(s)L(h)dV + \int \frac{1}{2}f(s) < g, h > dV,$$

where $L(h) = \frac{d}{dt}(s_{g+th})|_{t=0}$ is the linearization of the scalar curvature, given by

$$L(h) = -\Delta tr h + \delta\delta h - < r, h > . \tag{1.3}$$

Here, we require that $\int < g, h > dV = 0$, to preserve the condition that $g_t \in \mathbb{M}_1$, to first order. If $g$ is an $L^{2,2}$ positive definite symmetric form, i.e. an $L^{2,2}$ metric, then $s$ and $r$ are well-defined $L^2$ tensors on $M$. Similarly, if $h \in L^{2,2}$ then $L(h)$ is well-defined in $L^2$. Thus, the right side of the equation above is well-defined for $g$ and $h$ in $L^{2,2}$. It is then standard that $\mathcal{S}_f$ extends to a differentiable function on $\bar{\mathbb{M}}_1$. ∎

**Remark 1.2.** It is worth pointing out that this result is false if $f$ is only Lipschitz. As a concrete example of interest, for a given number $s_o < 0$, let $f(s) = \min(s, s_o)^2$. Then $f$ is Lipschitz, in fact piecewise $C^1$, but is not a $C^1$ function of $s$ *unless* $s_o = 0$.

Let $g$ be a metric on $M$ with $s_g \equiv s_o$ on some open set $\Omega \subset M$ and let $g_t = g + th$ be a volume preserving variation of $g$ supported in $\Omega$. If $t > 0$, i.e. if one computes the right derviative, then

one easily calculates

$$(\frac{d}{dt})^+ \mathcal{S}_f(g_t) = 2s_o \int_{\Omega^-} L(h)dV,$$

where $\Omega^- = \{L(h) < 0\} \subset \Omega$, while if $t < 0$, i.e. computing the left derivative gives

$$(\frac{d}{dt})^- \mathcal{S}_f(g_t) = -2s_o \int_{\Omega^+} L(h)dV,$$

where $\Omega^+ = \{L(h) > 0\} \subset \Omega$. Here we have used the fact that tr $h$ has mean value 0 on $M$ and supp $h \subset \Omega$. However, since $L(h)$ need not have mean value 0, $\mathcal{S}_f$ is not differentiable at such $g$ unless $s_o = 0$, (so $\mathcal{S}_f$ is essentially $\mathcal{S}_-^2$).

Of course if $f$ is a $C^1$ function of $s$ and $g$ is a metric realizing inf $\mathcal{S}_f$, then

$$\frac{d}{dt}\mathcal{S}_f(g + th) = 0, \quad \forall h \in T_g\mathbb{M}_1. \tag{1.4}$$

Next, we consider the $L^2$ gradient $\nabla \mathcal{S}_f$ of $\mathcal{S}_f$, i.e. the gradient w.r.t. the $L^2$ metric on $\mathbb{M}$. If $f$ is a $C^2$ function of $s$, for example $f(s) = s^p$ for $p \geq 2$, then $\nabla \mathcal{S}_f$ is essentially determined by $f$ and the formal adjoint of $L$. More precisely, since $L$ is a second order divergence operator, one may apply the divergence theorem twice to obtain

$$\nabla \mathcal{S}_f = L^*(f'(s)) + \tfrac{1}{2}f(s) \cdot g, \tag{1.5}$$

where $L^*$ is the $L^2$ adjoint of $L$, given by

$$L^*u = D^2u - \Delta u \cdot g - u \cdot r. \tag{1.6}$$

However, to perform the integration by parts twice, it is necessary that $f$ have at least two weak derivatives w.r.t. $s$. This is not the case for the function $f(s) = (s^-)^2$ giving the functional $\mathcal{S}_-^2$ in (0.2), (without the square root). We will return to this issue in more detail at the end of §1.

For this section, let

$$\mathsf{S}_-^2 = \int_M (s^-)^2 dV, \tag{1.7}$$

so that $\mathsf{S}_-^2$ differs from $\mathcal{S}_-^2$ only by absence of volume and square root factors. Thus, $\mathsf{S}_-^2 = \mathcal{S}_f$, for $f(s) = (s^-)^2$. For any functional $\mathcal{S}_f$, the function

$$\phi = f'(s) \tag{1.8}$$

from (1.5) plays a central role in the analysis. (It will become the function $\tau$ from §3 on). When $f = s^2, \phi = 2s$ while for $f = (s^-)^2, \phi = 2s^-$.

One elementary reason for preferring $\mathsf{S}_-^2$ or $\mathcal{S}_-^2$ is the following result.

**Lemma 1.3.** *Let $g$ be a $C^2$ metric on $M$ which is a critical point for $\mathsf{S}_-^2$ on $\mathbb{M}_1$, i.e. $(\mathsf{S}_-^2, g)$ satisfies (1.4). If $\sigma(M) \leq 0$, then $g$ is of constant curvature.*

*Proof.* Suppose first that $s = s_g \geq 0$ everywhere in $M$. Since $M$ admits no metric of positive scalar curvature, it is standard that $g$ must be flat, c.f. [B, (4.49)]. Thus, we may suppose that $U^- = \{x \in M : s_g < 0\} \neq \emptyset$. The functional $\mathsf{S}_-^2$ is smooth in regions of negative scalar curvature and hence in $U^-$ we have from (1.5) and the usual Lagrange multiplier argument

$$L^*\phi + \tfrac{1}{2}(s^-)^2 \cdot g = c \cdot g. \tag{1.9}$$

Taking the trace of (1.9) gives $-2\Delta\phi - \phi s + \tfrac{3}{2}(s^-)^2 = 3c$ and so since $\phi = 2s^-$,

$$-2\Delta\phi - \tfrac{1}{2}(s^-)^2 = 3c. \tag{1.10}$$

It follows of course that $c = -\tfrac{1}{6}\mathsf{S}_-^2(g)$.

Let $q$ be a point realizing min $s^-$ on $M$. Then $q \in U^-$ and evaluating (1.10) at $q$ gives $(s^-)^2(q) \leq \mathsf{S}^2_-(g)$. However, since $s^- \leq 0$ everywhere on $M$ and $vol_g M = 1$, we must have $(\min s^-)^2 \geq \mathsf{S}^2_-(g)$. It follows then from the maximum principle that $s^- = \text{const} < 0$ and hence by (1.9) and (1.6) that $z = 0$, i.e. $(M, g)$ is of constant (negative) curvature.

∎

Of course Lemma 1.3 also holds for $\mathcal{S}^2_-$. It is not known however if this result holds for $\mathcal{S}^2$ for instance, c.f. [An3, Prop. 1.1]. We also point out that the $2^{\text{nd}}$ variational formulas of these functionals are very complicated; it is not clear from such formulas that a local minimum of $\mathcal{S}_f$, for $f \neq (s^-)^2$, can be shown to be of constant curvature.

We now discuss briefly three basic reasons for preferring the functional $\mathcal{S}^2_-$ and its perturbation $I_\varepsilon^-$ in (0.7), over any other choice of $f$ in (1.1) or (1.2).

**(i). (Upper bound on $\phi$ ).** From results in §3, c.f. Theorem 3.11 and the discussion following it, for many natural choices of $f$ in (1.1), minimizers $g_{\varepsilon,f}$ for $I_{\varepsilon,f}$ have a uniform lower bound on the potential function $\phi$ in (1.8) as $\varepsilon \to 0$. (Actually this is done in §3 only for $I_\varepsilon^-$, but the reader may easily verify that the argument generalizes to a large class of functionals $I_{\varepsilon,f}$).

However in general, it seems difficult to obtain an upper bound. The choice $f(s) = (s^-)^2$ of course immediately gives an upper bound on $\phi$. This has two advantages. First, it simplifies the form of the possible blow-up limits. Second, and more importantly, it allows one to prove strong convergence *everywhere*, i.e. for all base points, c.f. Theorem 4.2 and compare with (I) above. Without such an upper bound, it is not clear if this is the case, c.f. Remark 4.3(ii). In turn, if strong convergence does not hold everywhere, then it is not clear that there exist *any* blow-up limits which are complete Riemannian manifolds, c.f. (ii) below.

**(ii). (Horizon problem).** Recall that in [AnI], the blow-up limits of sequences of Yamabe metrics are solutions of the static vacuum Einstein equations $L^*u = 0$. However, the construction of non-trivial limits is quite delicate and involved. In particular, the limits are not complete, and one expects that there is singularity formation at the horizon where $u$ approaches 0, c.f. [AnI, §5].

It turns out that these difficulties, which essentially arise from the need to avoid possible super-trivial solutions, ($u \equiv 0$), of the vacuum Einstein equations, persist for most all natural choices of $f$, e.g. for $f(s) = s^2$. Here the zero set of $\phi$ in (1.8) plays the role of the horizon.

However, for $f(s) = (s^-)^2$, and the associated minimizers $g_\varepsilon$ of $I_\varepsilon^-$, this problem does not occur. This will become clear in §4, where the horizon becomes the junction surface $\Sigma$, c.f. the statement of Theorem B.

Thus, one can avoid the involved construction of [AnI, Thm. 3.10] and quite directly and easily obtain non-trivial blow-up limits, c.f. §4. By construction, such limits will automatically be complete and without singularities. To obtain limits with such properties however, one must enlarge the class of limit solutions, i.e. allow limits which are not necessarily solutions of the static vacuum Einstein equations. But this will be the case for minimizers $g_{\varepsilon,f}$ of $I_{\varepsilon,f}$, for most any $f$ anyway.

**(iii). (Comparison arguments).** This issue is related to the issue in (III) above. Thus consider minimizing sequences $g_\varepsilon = g_{\varepsilon,f}$ for $\mathcal{S}_f$, with $\varepsilon \to 0$, with $g_\varepsilon$ a minimizer of $I_{\varepsilon,f}$ in (1.2). In order to prove that degenerations $\{g_\varepsilon\}$ arise solely from the topological structure of $M$, as noted in (III) above, one typically needs to carry out replacement or comparison arguments, i.e. if there is a topologically trivial degeneration, then find a comparison sequence $\{\tilde{g}_\varepsilon\}$ which satisfies

$$I_{\varepsilon,f}(\tilde{g}_\varepsilon) < I_{\varepsilon,f}(g_\varepsilon), \tag{1.11}$$

at least for $\varepsilon$ sufficiently small. Such a contradiction will then imply the degeneration is topologically non-trivial.

It turns out that it is much simpler to do comparison or replacement arguments, i.e. make use of (1.11), with $I_\varepsilon{}^-$ instead of $I_{\varepsilon,f}$, $f = s^2$ for instance, because the regions of positive scalar curvature in $(M, g)$ do not contribute at all to the value of $\mathcal{S}_-^2$. Namely, it turns out that for $f(s) = s^2$ and similar choices, all blow-up limits $(N, g')$ of $\{g_\varepsilon\}$ are scalar-flat, (c.f. Remark 3.2 and §4.1). Hence it is not clear if they admit any compact perturbations or comparisons; even if such do exist, they are hard to construct and are comparatively rare. This implies that it is not possible to use the analogue of (1.11) on the blow-up limits, (at least in an elementary way), since the comparison metrics must be scalar-flat. For similar reasons, it is difficult to construct useful comparisons to the blow-up metrics $g_\varepsilon'$ in (0.8) for $\varepsilon$ small.

On the other hand, the use of the weaker and hence more flexible functional $\mathcal{S}_-^2$ allows comparisons of the blow-up limits having non-negative scalar curvature. In many circumstances, there will exist a large and natural variety of such comparisons.

Similarly, it turns out that the choice of the perturbation by the $L^2$ norm of $z$ in (1.2) is optimal. For analytic reasons, c.f. §3, the $L^2$ norm is preferred to any other $L^p$ norm while the particular choice of $z$ allows for certain replacement arguments which are not possible even when $r$ is used in place of $z$.

Although such comparison arguments will not be used in this paper, they play an important role in [An6].

We close this section with a discussion of the difference between the derivative of $\mathsf{S}_-^2$ and the formal gradient given by (1.5) with $\phi = f'(s) = 2s^-$.

First, a straightforward computation gives

$$\frac{d}{dt}\mathsf{S}_-^2(g + th) = \int_M \phi \cdot L(h) + \frac{1}{2}(s^-)^2 dV = \int_M < L^*(\phi) + \frac{1}{2}(s^-)^2, h > dV + \int_\Sigma Q^+ + Q^-. \tag{1.12}$$

Here we have used the divergence theorem twice, which introduces the boundary term over $\Sigma = \{s = 0\}$. The term $Q^+$ may be computed to be

$$Q^+ = -\phi(< \nabla tr h + \delta h, \nu >) + < (< \nabla\phi, \nu > \cdot g - \nabla\phi \otimes \nu), h >, \tag{1.13}$$

where $\nu$ is the unit outward normal of $\partial M^+ = \Sigma$, for $M^+ = \{s \geq 0\}$. Similarly $Q^-$ is given by the same expression, with $\nu$ the unit outward normal of $\partial M^- = \Sigma$, for $M^- = \{s < 0\}$. To carry out these computations, we have assumed that $g$ is (sufficiently) smooth and that $\Sigma$ is a smooth hypersurface in $M$, with $\nabla\phi$ well-defined on each side of $\Sigma$.

Now $\phi = 0$ at $\Sigma$, so that the first bracketed term in (1.13) vanishes, both for $Q^+$ and $Q^-$. The field $\nabla\phi$ may be discontinuous at $\Sigma$, in that $\nabla\phi = 0$ from the right, i.e. from the $M^+$, side, but $\nabla\phi = \nabla s$ in general is not 0 from the left, i.e. from the $M^-$, side. Thus, only this boundary term remains.

Let $\nabla^-\phi$ be the gradient of $\phi$ in $M^-$, $\nu^-$ the unit normal to $\Sigma$ into $M^-$ and $dA$ Lebesque measure on the junction $\Sigma$. Also let

$$\nabla_{for}\mathsf{S}_-^2 = L^*(2s^-) + \frac{1}{2}(s^-)^2 \cdot g \tag{1.14}$$

denote the *formal* $L^2$ gradient of $\mathsf{S}_-^2$. With this understood, we have proved the following:

**Lemma 1.4.** *For $g$ as above, we have*

$$\nabla\mathsf{S}_-^2 = \nabla_{for}\mathsf{S}_-^2 + (< \nabla^-\phi, \nu^- > \cdot g - \nabla^-\phi \otimes \nu^-)dA. \tag{1.15}$$

■

In particular, at a critical point of $\mathsf{S}_-^2$, assumed smooth, the metric is not *necessarily* a solution of the Euler-Lagrange equations $\nabla_{for}\mathsf{S}_-^2 = L^*(2s^-) + \frac{1}{2}(s^-)^2 \cdot g = c \cdot g$, (and similarly for $\mathcal{S}_-^2$). However, we will see later in §3 that this deficiency can be overcome when the perturbation $I_\varepsilon{}^-$ is used, and in fact this leads to useful information.

## 2. Brief Background.

We summarize here some of the main background results that are needed, mainly related to the $L^2$ Cheeger-Gromov theory, c.f. [An3, §3]. $N$ will always denote either an oriented, open, connected 3-manifold, or a compact connected and oriented 3-manifold, possibly with boundary. $M$ is always a fixed closed, oriented and connected 3-manifold, without boundary.

The $L^2$ *curvature radius* $\rho(x)$ of $(N, g, x)$ is the radius of the largest geodesic ball centered at $x$, with $\rho(x) \leq diam N$, contained in the interior of $N$, such that for all $y \in B_x(\rho(x))$, and $B_y(s) \subset B_x(\rho(x)), s \leq dist(y, \partial B_x(\rho(x)))$, one has

$$\frac{s^4}{vol B_y(s)} \int_{B_y(s)} |r|^2 \leq c_o; \tag{2.1}$$

here $c_o > 0$ is a small parameter, c.f. Theorem 2.3 below, which is fixed once and for all. Note that from the definition, for $y \in B_x(\rho(x))$, we have

$$\rho(y) \geq dist(y, \partial B_x(\rho(x))). \tag{2.2}$$

We will occasionally also use the $L^{1,2}$ curvature radius $\rho^{1,2}$; this is defined in the same way as $\rho$, where in (2.1), $s^4$ is replaced by $s^6$ and $|r|^2$ by $|\nabla r|^2$. (Similar definitions hold w.r.t. $L^{k,p}$).

A sequence $x_i \in (N_i, g_i)$ is $(\rho, c)$ *buffered* if

$$\frac{\rho(x_i)^4}{vol B_y(\rho(x_i))} \int_{B_{x_i}((1-c)\rho(x_i))} |r|^2 \geq c \cdot c_o. \tag{2.3}$$

This condition prevents all of the curvature from concentrating in $L^2$ near the boundary of $B_{x_i}(\rho(x_i))$. The same sequence is *strongly* $(\rho, c)$ buffered if for all $y_i \in \partial B_{x_i}(\rho_i(x_i)), \rho_i(y_i) \geq c \cdot \rho_i(x_i)$. Similarly, this stronger condition prevents any of the curvature from concentrating in $L^2$ near the boundary, c.f. [AnI, Def.3.7] for further discussion.

The *volume radius* $\nu(x)$ at $x$ is given by

$$\nu(x) = sup\{r : \frac{vol\{B_y(s)\}}{s^3} \geq \mu, \forall B_y(s) \subset B_x(r), s \leq dist(y, B_x(\rho(x)))\}, \tag{2.4}$$

where again $\mu$ is a free (but fixed) small parameter, which measures the degree of volume collapse of the metric near $x$. To be concrete, we assume that $\mu = 10^{-1}$ throughout the paper.

We remark that $\rho$ and $\nu$ were defined slightly differently in [AnI] and [An3]. The correct definitions here in (2.1) and (2.4) should have been used throughout [AnI] and [An3].

The $L^{2,2}$ *harmonic radius* $r_h(x)$ of $(N, g, x)$ is the radius of the largest geodesic ball at $x$ on which one has a harmonic coordinate chart in which the metric coefficients $g_{\alpha\beta}$ of $g$ are bounded in $L^{2,2}$ norm by a fixed constant $C < \infty$. We observe that in dimension 3, Sobolev embedding gives $L^{2,2} \subset C^{1/2}$. Note also that all radii above scale as distances.

A sequence of metrics $g_i$ on $N$ (or $N_i$) is *non-collapsing* at $x_i$ if there is a constant $\nu_o > 0$ such that

$$\nu_i(x_i) \geq \nu_o \cdot \rho_i(x_i). \tag{2.5}$$

If (2.5) does not hold, for some $\nu_o > 0$, we say the sequence is *collapsing* at $x_i$. Note that the condition (2.5) is scale-invariant.

While the discussion in [AnI] basically ruled out, by assumption, any collapse behavior, we will not exclude the possibility of collapse in this paper. We summarize here some of the essential facts of collapse in dimension 3, and refer to [CG1,2], [Fu], [O], and especially [R], for further details.

- Any F-structure on a 3-manifold admits a polarized substructure, i.e. an F-structure for which the local torus actions are locally free. Thus, we will essentially always work with polarized F-structures.

- A polarized F-structure on $N$ determines a graph manifold structure on $N$, and vice versa.
- An F-structure $\mathcal{F}$ on $N$ is called *injective*, if for any orbit $\mathcal{O}_x$ in $\mathcal{F}$, the inclusion $\mathcal{O}_x \subset N$ induces an injection $\pi_1(\mathcal{O}_x, x) \to \pi_1(N, x)$. A polarized F-structure is injective if and only if there is no solid torus component in the graph manifold decomposition.
- Given arbitrary $\Lambda, D < \infty$, there exists $\delta = \delta(\Lambda, D) > 0$ such that if $(N, g)$ is a Riemannian 3-manifold, with metric satisfying

$$|R| \leq \Lambda, \quad diam N \leq D, \quad \text{and} \quad vol N \leq \delta,$$

  then there is a manifold $U \subset N$, with $dist(\partial U, \partial N) \leq \mu = \mu(\delta)$, such that $U$ is a Seifert fibered space, or a torus bundle over a 1-manifold, (i.e. a Sol-manifold); here $\mu(\delta) \to 0$ as $\delta \to 0$. Further $U$ admits a pure polarized F-structure $\mathcal{F} = \mathcal{F}_\delta$, possibly depending on $\delta$, with orbits given by the fibers, for which the diameter of the orbits in the metric $g$ is bounded by $\delta' = \delta'(\Lambda, D)$, where $\delta' \to 0$ as $\delta \to 0$. In either of these cases, (Seifert fibered or torus bundle), the fibers of $\mathcal{F}_\delta$ are injective in $\pi_1(U)$, provided $U$ is not covered by $S^3$. In particular, this holds if $U$ is an open 3-manifold.

**Remark 2.1.** Let $(N, g_i)$ be a collapsing sequence of metrics on $N$, and assume $\rho_i(x) \geq \rho_o, \forall x \in N$, for some $\rho_o > 0$ and diam $N \leq D < \infty$. If $\{g_i\}$ collapses along a sequence of injective F-structures, then one may unwrap the collapsing sequence $\{g_i\}$ by passing to a sufficiently large finite cover $\bar{N}$, or to the universal cover $\widetilde{N}$. The metrics $g_i$, lifted to $\widetilde{N}$ say, do not collapse and thus, (see below), have a subsequence converging in the weak $L^{2,2}$ topology to a limit $L^{2,2}$ metric $g$ on $\widetilde{N}$. In other words, one is able to resolve or remove the degeneration of $\{g_i\}$ in this case, c.f. [An3, Cor.3.14ff] for further discussion.

The limit metric $g$ on $\bar{N}$ or $\widetilde{N}$ has some further important structure in this case. Namely, in passing to the universal cover, the generators of the action of $\pi_1(F)$, where $F$ is a typical fiber of the F-structure, have displacement functions converging to 0 on $\widetilde{N}$, corresponding to the collapse of the fibers under $\{g_i\}$. Thus, these actions converge to smooth, free and proper actions of $\mathbb{R}$ or $\mathbb{R}^2$ on $\widetilde{N}$, according to whether the sequence of F-structures is of rank 1 or 2, (or a sequence of rank 1 structures converging to a rank 2 structure). This feature that one may unwrap collapsing sequences under the bounds above is special to dimension 3.

It is straightforward to verify from the definition that the curvature radius behaves naturally under coverings; in fact, if $\pi : \bar{N} \to N$ is a locally isometric covering and $\bar{x}$ a lift of $x$, then there exists $c_1 = c_1(c_o) > 0$ such that

$$c_1^{-1} \rho(x) \geq \bar{\rho}(\bar{x}) \geq c_1 \rho(x).$$

**Remark 2.2.** Let $\{g_i\}$ be a sequence of metrics on $M$ with $\mathrm{vol}_{g_i} M \leq 1$ and

$$\int_M |r_{g_i}|^2 dV_{g_i} \leq \Lambda < \infty. \tag{2.6}$$

Then a subsequence of $\{g_i\}$ either converges to a limit $L^{2,2}$ metric $g_\infty$ on $M$ in the weak $L^{2,2}$ topology, (modulo diffeomorphisms), or it collapses everywhere, or it forms cusps, a mixture of these two cases. More precisely, in the last case, there is a maximal open domain $\Omega$, weakly embedded in $M$ in the sense that any compact subset of $\Omega$ with smooth boundary embeds smoothly in $M$, such that the metrics $g_i$ converge weakly in $L^{2,2}$ to a limit $L^{2,2}$ metric $g_\infty$ on $\Omega$; the convergence is uniform on compact subsets. For sufficiently large $K \subset \Omega$ as above, the complement $M \setminus K$ is a graph manifold, which is partially collapsed by $\{g_i\}$ along a sequence of F-structures.

A similar version of this result holds for open 3-manifolds $N$, provided one stays a fixed distance away from the boundary.

We refer to [An3, Thm. 3.19], [An4, Rmk. 5.5] for further details and discussion of this result.

We will often use the following result from [An3, Cor.3.14-Rmk.3.18], which characterizes the possible behavior of metrics within the $L^2$ curvature radius. In particular, one has no cusp formation within the $L^2$ curvature radius.

**Theorem 2.3.** *Let $(B_{x_i}(1), g_i)$ be a sequence of smooth metrics defined on a geodesic ball about $x_i$, of radius 1, complete up to the boundary. Suppose the $L^2$ curvature radius satisfies $\rho(x_i) \geq 1$, so that for all $B_y(s) \subset B_x(1), s \leq 1,$*

$$\frac{s^4}{vol B_y(s)} \int_{B_y(s)} |r|^2 \leq c_o,$$

*where $c_o$ is sufficiently small. Then there are diffeomorphisms $\psi_i$ of $B_{x_i}(1)$ such that a subsequence of $\psi_i^* g_i$ either converges in the weak $L^{2,2}$ topology, uniformly on compact subsets, to a limit $L^{2,2}$ metric $g$ on $B_x(1)$, or the subsequence collapses along a sequence of injective F-structures defined on $B_{x_i}(1 - \delta_i)$, where $\delta_i \to 0$ as $i \to \infty$.*

We note that this result is false if $c_o$ is too large. Thus, assuming $c_o$ is sufficiently small, $c_o = 10^{-3}$ suffices, this result implies in particular that the unit ball $\widetilde{B}_{x_i}(1 - \delta)$ in the universal cover is diffeomorphic to a 3-ball, and the metric $g_i$ is close, (depending on $c_o$), in the $L^{2,2}$ topology, to the flat metric on $B(1 - \delta) \subset \mathbb{R}^3$, for any fixed $\delta > 0$. This result will be used frequently in particular to obtain regularity estimates.

We note that one also has the following estimate on the volume radius $\nu$ in (2.4):

$$sup\nu \leq c \cdot inf\nu,$$

where the sup and inf are taken over $B_x(\rho(x))$; here $c$ depends only on $c_o$, c.f [An3,(3.62)]. Hence $\nu$ decreases at most geometrically w.r.t. the radius $\rho$. Thus, if one has a sequence of geodesic balls $(B_{x_i}(R), g_i)$ with a uniform lower bound

$$\rho_i(y) \geq \rho_o > 0 \text{ for } y \in B_{x_i}(R), \tag{2.7}$$

then either there is a constant $\nu_o = \nu_o(R)$, such that

$$\nu_i(y) \geq \nu_o > 0, \tag{2.8}$$

everywhere in $B_{x_i}(R)$, or

$$\nu_i(y) \to 0, \tag{2.9}$$

everywhere in $B_{x_i}(R)$. This implies in particular that cusps cannot form in $B_{x_i}(R)$ under the bound (2.7).

In the case of (2.8) it follows from repeated applications of Theorem 2.3 that a subsequence of $(B_{x_i}(R), g_i)$ converges in the weak $L^{2,2}$ topology to a limit $L^{2,2}$ metric $g$ on $B_x(R)$.

In the case of (2.9), the metrics $g_i$ collapse along a sequence of F-structures on $B_{x_i}(R)$, and the F-structures are injective on $B_{x_i}(R)$ provided $B_{x_i}(R)$ is not covered by $S^3$.

We will often work with the trace-free Ricci curvature $z = r - \frac{s}{3} \cdot g$ in place of $r$. In this regard, let $\rho_z$ denote the $L^2$ curvature radius w.r.t. $z$, i.e. (2.1) with $z$ in place of $r$. Then from [An4, (5.26)], if

$$s \geq \lambda > -\infty, \tag{2.10}$$

throughout $B_x(\rho_z(x))$, one has the estimate

$$\rho(x) \leq \rho_z(x) \leq \rho_o \cdot \rho(x), \tag{2.11}$$

where $\rho_o$ depends only on $c_o$ in (2.1) and $\lambda$. Thus these curvature radii are essentially equivalent under the bound (2.10).

Further, the result discussed in Remark 2.2 holds with bounds of the form

$$\int_M |z_{g_i}|^2 dV_{g_i} \leq \Lambda, \quad \int_M (s_{g_i}^-)^2 dV_{g_i} \leq \Lambda, \tag{2.12}$$

in place of (2.6). This is proved in [An4, Prop.5.3] and the discussion following it concerning [An3, Thm.3.19]. (In effect, this means that, under the bounds (2.12), the scalar curvature can become large, positively, only in regions which have very small volume radius).

Finally, in working with sequences $\{a_i\}$ or families, for instance $\{g_\varepsilon\}$, we will often pass to subsequences of $\{a_i\}$ or subsequences within the family $\{g_\varepsilon\}$, as for example in Theorem 2.3, without explicitly indicating the subsequence. The choice of subsequences will only be made explicit where it is necessary to be careful in this respect. A sequence is said to sub-converge if a subsequence converges.

In Theorem 4.2, we will need the following elliptic estimates, proved in [An4, Prop. 3.7, Lem. 3.8] and the remark following those results. These estimates are stated here for convenience.

**Proposition 2.4.** *Let $(B, g)$, $B = B_x(1)$ be a geodesic ball of radius 1, with $r_h(x) \geq 1$. For a given smooth function $f$ on $B$, let $u$ be a smooth solution of*

$$\Delta u = f,$$

*on $B$. Let $B_\delta = B_x(\delta)$, where $\delta > 0$ is (arbitrarily) small, depending only on the choice of $c_o$ in (2.1). Then there is a constant $c > 0$, depending only on $\delta^{-1}$, such that*

$$\int_B u \leq c(\delta^{-1}) \int_B |f| + 2\pi u_{av}(\delta) + \int_{B_\delta} u, \tag{2.13}$$

*where $u_{av}$ is the average of $u$ on $S_\delta = \partial B_\delta$. Further for any ball $\bar{B} \subset B$, $d = dist(\partial \bar{B}, \partial B)$, and any $\mu > 0$, there is a constant $c = c(\mu, d, \delta^{-1})$ such that*

$$||u||_{L^{3-\mu}(\bar{B})} \leq c(\mu, d, \delta^{-1}) \cdot \left( \int_B |f| + 2\pi u_{av}(\delta) + \int_B |u^-| \right). \tag{2.14}$$

### 3. Basic Properties of the Curvature Functional.

In this section, we discuss the functional $\mathcal{S}_-^2$ and its perturbation $I_\varepsilon^-$,

$$I_\varepsilon^- = \varepsilon v^{1/3} \mathcal{Z}^2 + \mathcal{S}_-^2 = \varepsilon v^{1/3} \int |z|^2 + \left( v^{1/3} \int (s^-)^2 dV \right)^{1/2}. \tag{3.0}$$

and the initial geometric properties of the minimizers $g_\varepsilon$ of $I_\varepsilon^-$. The power $1/2$ in (3.0) is needed in order to handle the degeneracy of $\int (s^-)^2$ at scalar-flat metrics, i.e. when $\sigma(M) = 0$. In the case $\sigma(M) < 0$, there is no degeneracy, and one may dispense with this power. However, in order to deal with both cases together, we work with (3.0). The volume powers are chosen of course so that $I_\varepsilon^-$ is scale invariant.

**§3.1.** First, we need to prove the statement (0.4).

**Proposition 3.1.** *Suppose $\sigma(M) \leq 0$. Then*

$$inf_{\mathbb{M}} \mathcal{S}_-^2 = |\sigma(M)|. \tag{3.1}$$

*Proof.* Recall that by definition, $|\sigma(M)| = \inf |s_{\bar{g}}|$, where $\bar{g}$ is a unit volume Yamabe metric, when $\sigma(M) \leq 0$. Thus it is obvious that $\inf \mathcal{S}_-^2 \leq |\sigma(M)|$. To obtain the opposite inequality, for any metric $g \in \mathbb{M}_1$, let $\bar{g}$ be the associated (unique) unit volume Yamabe metric conformal to $g$. Thus write $g = u^4 \bar{g}$, so that $u > 0$ and $dV = u^6 d\bar{V}$. The equation relating the scalar curvatures $\bar{s}$ and $s$ of $\bar{g}$ and $g$ is

$$u^5 s = -8\bar{\Delta} u + \bar{s} u.$$

Hence

$$\mathcal{S}_-^2(\bar{g}) = |\bar{s}| = -\int \bar{s} d\bar{V} = -\int s u^4 d\bar{V} - 8\int u^{-1} \bar{\Delta} u d\bar{V} = -\int s u^4 d\bar{V} - 8\int |\bar{\nabla} log u|^2 d\bar{V} \leq$$

$$\leq |\int (s^-) u^4 d\bar{V}| \leq (\int (s^-)^{3/2} u^6 d\bar{V})^{2/3} = (\int (s^-)^{3/2} dV)^{2/3} \leq (\int (s^-)^2 dV)^{1/2} = \mathcal{S}_-^2(g).$$

This gives the opposite inequality $|\sigma(M)| \leq inf \mathcal{S}_-^2$, and hence (3.1). ∎

Clearly, a natural choice of minimizing sequence for $\mathcal{S}_-^2$ is to choose a sequence with the least curvature in $L^2$. Thus, given $\varepsilon > 0$, let $g_\varepsilon$ be a unit volume metric realizing the infimum of $I_\varepsilon^-$ in (3.0) in the sense that $g_\varepsilon$ is a $L^{2,2}$ metric satisfying $I_\varepsilon^-(g_\varepsilon) = \inf_{g \in \mathbb{M}} I_\varepsilon^-(g)$; the existence of such metrics will be proved in Theorem 3.8 below. Let $\delta = \delta(g_\varepsilon) = \mathcal{S}_-^2(g_\varepsilon) - |\sigma(M)|$. Then it is clear that

$$\mathcal{Z}^2(g_\varepsilon) \leq \mathcal{Z}^2(g), \tag{3.2}$$

for all metrics $g$ on $M$ satisfying $vol_g M = 1$ and $\mathcal{S}_-^2(g) - |\sigma(M)| \leq \delta$. (The problem of minimizing $\mathcal{Z}^2$ on a given sublevel set of $\mathcal{S}_-^2$ is related to that of minimizing $I_\varepsilon^-$ by the usual Lagrange multiplier method).

To see that the family $\{g_\varepsilon\}$ forms a minimizing family for $\mathcal{S}_-^2$ as $\varepsilon \to 0$, note that as a function of $\varepsilon$, $I_\varepsilon^-(g)$ is monotone increasing in $\varepsilon$. Hence for any $\varepsilon > 0$ and $g \in \mathbb{M}$,

$$lim_{\varepsilon \to 0} I_\varepsilon^-(g_\varepsilon) \leq I_\varepsilon^-(g_\varepsilon) \leq I_\varepsilon^-(g). \tag{3.3}$$

Clearly $I_\varepsilon^-(g) \geq \mathcal{S}_-^2(g)$, so that $I_\varepsilon^-(g_\varepsilon) \geq |\sigma(M)|$ for all $\varepsilon$. On the other hand, we may choose $g$ and $\varepsilon$ so that $I_\varepsilon^-(g) \leq |\sigma(M)| + \delta$, for any given $\delta > 0$. Thus, we see $I_\varepsilon^-(g_\varepsilon) \to |\sigma(M)|$, which implies

$$\left( \int_M (s_{g_\varepsilon}^-)^2 dV_{g_\varepsilon} \right)^{1/2} \to |\sigma(M)|, \quad \text{and} \tag{3.4}$$

$$\varepsilon \int |z_{g_\varepsilon}|^2 \to 0. \tag{3.5}$$

Instead of letting $\varepsilon \to 0$, consider for a moment the situation when $\varepsilon \to \infty$. Then the $\mathcal{S}_-^2$ term in $I_\varepsilon^-$ becomes negligible. In the limit $\varepsilon \to \infty$, the variational problem of minimizing $I_\varepsilon^-$ becomes that of minimizing $\mathcal{Z}^2$ among *all* unit volume metrics on $M$. Thus the family of functionals $I_\varepsilon^-$ as $\varepsilon$ varies between 0 and $\infty$ interpolates in this sense between $\mathcal{S}_-^2$ and $\mathcal{Z}^2$.

**Remark 3.2.** All of the discussion above holds equally well for the functionals,

$$I_\varepsilon = \varepsilon v^{1/3} \int |z|^2 + \left( v^{1/3} \int s^2 \right)^{1/2}, \quad I_\varepsilon' = \varepsilon v^{1/3} \int |r|^2 + \left( v^{1/3} \int s^2 \right)^{1/2}, \tag{3.6}$$

and related functionals $I_{\varepsilon,f}$ as discussed in §1, in place of $I_\varepsilon^-$. It should be kept in mind however that as $\varepsilon \to 0$, there is not necessarily a uniform upper bound on $\mathcal{S}^2(g_\varepsilon)$, for $g_\varepsilon$ a minimizer of $I_\varepsilon^-$, as is obviously the case for minimizers of $I_\varepsilon$ or $I_\varepsilon'$. Thus, it is possible that the $L^2$ norm of the scalar curvature of $g_\varepsilon$ goes to infinity as $\varepsilon \to 0$, c.f. also Remark 3.5 below.

We also point out that some of the arguments to follow are somewhat simpler when $I_\varepsilon^-$ has $z$ replaced by the full Ricci curvature $r$. For instance, in this situation, a bound on $I_\varepsilon^-$ trivially gives a bound on the $L^2$ norm of $s$, (depending on $\varepsilon$). However, for reasons which will become only fully clear in [An6], (related to issue (iii) in §1), we need to work with $I_\varepsilon^-$ as given in (3.0).

The functional $I_\varepsilon'$ in (3.6) and the existence and structure of minimizers for $I_\varepsilon'$ has been treated in full detail in [An3,§8], (with the insignificant difference that the full curvature $R$ was used in place of the Ricci curvature $r$, and the $\mathcal{S}^2$ term had no square root). Much of this work carries over without change to the functional $I_\varepsilon{}^-$. This is because these results are primarily deduced from the dominant curvature terms $\mathcal{R}^2$ in $I_\varepsilon'$ and $\mathcal{Z}^2$ in $I_\varepsilon{}^-$.

**(I). Existence of $L^{2,2}$ Minimizers.**

Choose a fixed $\varepsilon > 0$. If $\{g_k\}$ is a unit volume minimizing sequence for $I_\varepsilon{}^-$, then clearly

$$||z_{g_k}||_{L^2(g_k)} \leq \Lambda, \quad \text{and} \quad ||(s^-)_{g_k}||_{L^2(g_k)} \leq \Lambda,$$

for some $\Lambda = \Lambda(\varepsilon) < \infty$.

By Remark 2.2 and (2.12), it follows that a subsequence of $\{g_k\}$ converges, modulo diffeomorphisms of $M$, on a maximal domain $\Omega_\varepsilon$ to a limit $L^{2,2}$ metric $g_\varepsilon$ defined on $\Omega_\varepsilon$. Here the convergence is in the weak $L^{2,2}$ topology on $(\Omega_\varepsilon, g_\varepsilon)$, uniform on compact subsets. The manifold $\Omega_\varepsilon$ weakly embeds in $M$, denoted as

$$\Omega_\varepsilon \subset\subset M,$$

in the sense that any domain $K$ with smooth and compact closure in $\Omega_\varepsilon$ embeds as such in $M$.

As in Remark 2.2, for sufficiently large $K \subset \Omega_\varepsilon$, the complement $M \setminus K$ in $M$ is a graph manifold while part of $M \setminus K$ is collapsed along a sequence of F-structures or graph manifold structures to a lower dimensional space under $\{g_k\}$. The domain of convergence $\Omega_\varepsilon$ is empty only if $M$ itself is a graph manifold.

If $M$ is a graph manifold, then inf $I_\varepsilon{}^- = 0$, for all $\varepsilon \geq 0$, and a minimizing sequence collapses all of $M$ along a sequence of F-structures, except possibly in case $M$ is a flat manifold, c.f. also §6.3. Conversely, if inf $I_\varepsilon{}^- = 0$, for some $\varepsilon > 0$, then $M$ is a graph manifold. The topological structure and geometrization of graph manifolds is completely understood, c.f. [St], [Wa] or also [An6, §2].

Thus, we assume from now on that $M$ is not a graph manifold. It follows that either $\Omega_\varepsilon = M$, or $\Omega_\varepsilon \neq \emptyset$ and any large compact set in $\Omega_\varepsilon$ with smooth boundary embeds in $M$.

In this case, one has

$$0 < \varepsilon \int_{\Omega_\varepsilon} |z_{g_\varepsilon}|^2 + (\int_{\Omega_\varepsilon} (s_{g_\varepsilon}^-)^2)^{1/2} \leq inf_{\mathbb{M}_1} I_\varepsilon{}^-, vol_{g_\varepsilon}\Omega_\varepsilon \leq 1. \tag{3.7}$$

Such a limit $(\Omega_\varepsilon, g_\varepsilon)$ of a minimizing sequence will be called a *minimizing pair* for $I_\varepsilon{}^-$ and $g_\varepsilon$ a *minimizer* of $I_\varepsilon{}^-$; (in Theorem 3.9 below it is proved that the inequalities in (3.7) are in fact equalities).

**(II). The Euler-Lagrange equations.**

The form of $\nabla I_\varepsilon{}^-$ is important for the work in later sections. As noted in §1, the functional $\mathcal{S}_-^2$, although differentiable, does not have a "smooth" gradient, i.e. by (1.15), the $L^2$ gradient is not the same as the formal gradient. Thus, the initial considerations below are formal. They will be made rigorous following this in Theorem 3.3.

Let $g$ be any smooth metric on $M$ and let $g_t = g + th$ be a smooth variation of $g$, with $g' = h = \frac{dg_t}{dt}|_{t=0}$. On the curve $g_t$, we have

$$(v^{1/3}\mathcal{Z}^2)' = (v^{1/3})' \int |z|^2 dV + v^{1/3}(\int (|z|^2 dV)'.$$

and $dV' = \frac{1}{2} < g, h > dV$, so that

$$(v^{1/3})' = \frac{1}{6}v^{-2/3} \int < g, h > dV.$$

Next $\int(|z|^2 dV)' = \int < \nabla \mathcal{Z}^2, h >$, and from [B,4.66],

$$\nabla \mathcal{Z}^2 = D^* Dz + \frac{1}{3} D^2 s - 2 \overset{\circ}{R} \circ z + \frac{1}{2}(|z|^2 - \frac{1}{3}\Delta s) \cdot g. \tag{3.8}$$

Also $(\mathcal{S}_-^2)' = \frac{1}{2}(\mathcal{S}_-^2)^{-1}[(v^{1/3})'\mathsf{S}_-^2 + v^{1/3}(\mathsf{S}_-^2)]'$, for $\mathsf{S}_-^2$ as in (1.7). Assembling these computations, one easily computes that the formal gradient of $I_\varepsilon{}^-$ at $g \in \mathbb{M}$ is,

$$\nabla I_\varepsilon{}^- = v^{1/3}\left(\varepsilon \nabla \mathcal{Z}^2 + L^*(\tau) + \chi \cdot g\right). \tag{3.9}$$

Here $\chi = (\frac{s\tau}{4} + c)$ and

$$\tau = \frac{s^-}{\sigma}, \quad \sigma = (v^{1/3}\int(s^-)^2 dV)^{1/2}, \tag{3.10}$$

$$c = \frac{1}{12\sigma v}\int(s^-)^2 dV + \frac{\varepsilon}{6v}\int|z|^2 dV, \tag{3.11}$$

(compare with (1.6)). We have used the fact that $s\tau = s^-\tau$. and all quantities are w.r.t. the metric $g$. Taking the trace gives the associated trace equation

$$tr\nabla I_\varepsilon{}^- = -v^{1/3}\left(2\Delta(\tau + \frac{\varepsilon s}{12}) + \frac{1}{4}s\tau + \frac{\varepsilon}{2}|z|^2 - 3c\right). \tag{3.12}$$

The function $\tau$, the $L^2$ average of $s^-$ on $M$, plays a crucial role throughout this paper.

As noted in §1, a minimizer $g_\varepsilon$ of $I_\varepsilon{}^-$ may be, but is not necessarily, even a weak solution of the equation $\nabla I_\varepsilon{}^- = 0$. We address this issue next.

**(III). Regularity.**

By way of background, it is proved in [An3, Thms.4.1, 8.1] that a minimizer of the smooth functional $I'_\varepsilon$ in (3.6) is a $C^\infty$, in fact real-analytic, solution of the Euler-Lagrange equations $\nabla I'_\varepsilon = 0$, for any $\varepsilon > 0$ fixed. The expression for $\nabla I'_\varepsilon$ is almost identical to the expressions (3.9)-(3.12); the only difference is that $\mathcal{Z}^2$ in (3.9)-(3.12) should be replaced by $\mathcal{R}^2 = \int|r|^2$, and more importantly, $\tau$ in (3.10) is replaced by $\tau' = s/\sigma$, with $\sigma$ the $L^2$ norm of $s$, i.e. there is no cutoff at 0.

As mentioned above, such regularity need not hold for $I_\varepsilon{}^-$, since this functional is not $C^\infty$. Nevertheless, using a regularization of $I_\varepsilon{}^-$, we prove that there are minimizers $g_\varepsilon$ of $I_\varepsilon{}^-$, (in the sense following (3.7)), which have sufficient partial regularity to ensure that they are weak solutions of the Euler-Lagrange equations.

**Theorem 3.3.** *For any fixed $\varepsilon > 0$, there exist minimizing pairs $(\Omega_\varepsilon, g_\varepsilon)$ of $I_\varepsilon{}^-$ with $g \in L^{3,p}$ locally, for any $p < \infty$. In particular, $g_\varepsilon \in C^{2,\alpha}$, for any $\alpha < 1$, and the Ricci curvature $r_\varepsilon$ of $g_\varepsilon$ is locally in $L^{1,p} \cap C^\alpha$. The potential function $\omega$,*

$$\omega = \tau + \frac{\varepsilon}{12}s, \tag{3.13}$$

*as in (3.12) is locally in $L^{3,p}$, while each summand $\tau$ and $\frac{\varepsilon}{12}s$ is locally Lipschitz smooth.*

*The metric $g_\varepsilon$ is a $L^{3,p}$ weak solution of the Euler-Lagrange equations*

$$\nabla I_\varepsilon{}^- = 0, \tag{3.14}$$

*defined in (3.9)-(3.12), for some constant $c$ in (3.11). Further, away from the locus*

$$\{s_\varepsilon = 0\},$$

*i.e. in a neighborhood of any point $x \in \Omega_\varepsilon$ for which $s_\varepsilon(x) \neq 0$, the metric $g_\varepsilon$ is real-analytic, and is a smooth solution there of the Euler-Lagrange equations.*

By an $L^{3,p}$ weak solution, it is meant that (3.14) holds when paired by integration with any $L^{1,q}$ 2-tensor of compact support in $\Omega_\varepsilon$, $\frac{1}{p} + \frac{1}{q} = 1$, and integration by parts is used to transfer $4^{\text{th}}$ derivatives of $g_\varepsilon$ to $3^{\text{rd}}$ derivatives. We will see later in Theorem 3.9 that the constant $c$ is in fact given by (3.11).

The proof of Theorem 3.3 will proceed in several steps. First, let $\phi_\delta : \mathbb{R} \to \mathbb{R}$ be a $C^\infty$ smoothing of the function $s^- = \min(s, 0)$. We choose $\phi_\delta$ so that $\phi_\delta(s) = s, s \leq 0$, $\phi_\delta(s) = \delta s, s \geq \delta$, and on $[0, \delta]$, $\phi_\delta$ is a smooth concave interpolation so that $\phi_\delta^{-1}(0) = 0$. (Here we are considering $s$ as a parameter on $\mathbb{R}$). We also require that $((\phi_\delta)^2)'' > 0$ everywhere. It is easy to see that such $C^\infty$ functions exist. One may also choose $\phi_\delta$, if desired, to be a real-analytic approximation to this $C^\infty$ smoothing of $s^-$. It follows that the $C^\infty$ functional

$$I_{\varepsilon,\delta} = \varepsilon \int |z|^2 dV + \left( \int \phi_\delta(s)^2 dV \right)^{1/2} : \mathbb{M}_1 \to \mathbb{R}, \tag{3.15}$$

converges to $I_\varepsilon^-|_{\mathbb{M}_1}$ as $\delta \to 0$. Note that

$$\frac{d(s^-)^2}{ds} = 2s^-, \tag{3.16}$$

and that $s^-$ is monotone non-decreasing as a function of $s$. By construction, the derivative $(\phi_\delta^2)' = d(\phi_\delta^2)/ds$ is strictly monotone increasing in $s$, for each $\delta > 0$. Hence the function

$$\tau_\delta = \frac{\phi_\delta \phi_\delta'}{\sigma_\delta}, \tag{3.17}$$

where $\sigma_\delta$ is the constant $\left( \int \phi_\delta(s)^2 dV \right)^{1/2}$, is invertible as a function of $s$. Observe here that $\tau_\delta = s^-/\sigma_\delta$ when $s \leq 0$.

For the same reasons as in (I) above, there exists an $L^{2,2}$ minimizer $g_\delta$ of $I_{\varepsilon,\delta}$ defined on a maximal domain $\Omega_{\varepsilon,\delta}$, for any $\delta > 0$. Since $I_{\varepsilon,\delta}$ is a smooth functional on $\mathbb{M}_1, g_\delta$ is a $L^{2,2}$ weak solution of the Euler-Lagrange equations

$$\nabla I_{\varepsilon,\delta} = 0.$$

Following the derivation in (II) above, it is easy to see that these equations are of the form

$$\varepsilon \nabla \mathcal{Z}^2 + L^*(\tau_\delta) + \chi_\delta \cdot g = 0, \tag{3.18}$$

$$2\Delta(\tau_\delta + \frac{\varepsilon s}{12}) + (s\tau_\delta - \frac{3}{4}\frac{\phi_\delta^2}{\sigma_\delta}) = -\frac{\varepsilon}{2}|z|^2 + 3c, \tag{3.19}$$

where $\chi_\delta = \frac{1}{4}\frac{\phi_\delta^2}{\sigma_\delta} + c$, and $c$ is some constant. One easily checks that these expressions become (3.9) and (3.12) in the limit $\delta \to 0$, except possibly for the value of the (global) constant $c$. These equations hold when paired with $L^{2,2}$ test functions or forms of compact support, in the usual sense of weak solutions.

We claim that, for any fixed $\delta > 0$, any $L^{2,2}$ weak solution of (3.18)-(3.19) is $C^\infty$ smooth, and is real-analytic if the function $\phi_\delta$ is a real-analytic function of $s$. The proof of this is exactly the same as the proof of the smooth regularity of minimizers of $I_\varepsilon'$ in [An3, §4,§8], and follows the usual bootstrap method for proving regularity of elliptic equations. As noted above, the only difference in the equations (3.18)-(3.19) from the Euler-Lagrange equations for $I_\varepsilon'$, (c.f. [An3,(8.1)-(8.2)]) is the difference of $\mathcal{R}^2$ and $\mathcal{Z}^2$, and, more importantly, the form of the potential function

$$\omega_\delta = \tau_\delta + \frac{\varepsilon s}{12}. \tag{3.20}$$

Since regularity is a purely local issue, in the following, we work on small balls $B \subset \Omega_{\varepsilon,\delta}$ within the $L^{2,2}$ harmonic radius.

First, observe that $g_\delta$ has $L^{2,2}$ components in a local harmonic coordinate chart, and that $z = z_{g_\delta}$ is (locally) in $L^2$. By Sobolev embedding, $L^{2,2} \subset C^{1/2}$. The bootstrap or iteration process starts with a value of $k < \frac{1}{2}$ and $p = 2$ and uses the Sobolev embedding $L^1 \subset L^{k-2,2} = (L_o^{2-k,2})^*$, where $(L_o^{2-k,2})^*$ is the dual space of the Sobolev space $L_o^{2-k,2}$ of functions of compact support. We refer to [Ad,Thm.7.48] or [Tr, Rmk.25.3] for the definition of fractional-order Sobolev spaces, and also to these texts for Sobolev spaces with negative exponents.

Thus, $|z|^2$ belongs to some Sobolev space $L^{k-2,2}$ and consider the Laplacian as an operator $\Delta : L^{k,2} \to L^{k-2,2}$. From (3.19), regularity theory for $\Delta$ implies that $\omega_\delta \in L^{k,2}$; (this uses the fact that $g_\delta \in L^{2,2}$). Now from the assumptions on $\tau_\delta$ following (3.17), it follows from the form of $\omega$ in (3.20) that $\omega_\delta$ is a monotone increasing function of $s$, and hence invertible. Since $\phi_\delta$ is assumed to be a smooth, or real-analytic, function of $s$, it follows that $s \in L^{k,2}$, i.e. the second summand in (3.20) has an improved regularity. It then follows of course that also $\tau_\delta \in L^{k,2}$. This in turn implies that $L^*(\tau_\delta) \in L^{k-2,2}$, so that from (3.18), $\varepsilon \nabla \mathcal{Z}^2 \in L^{k-2,2}$. From the form of (3.8), since we have $D^2 s \in L^{k-2,2}$ and the curvature terms are in $L^1$, it follows that $D^* D z \in L^{k-2,2}$. Using again the $L^{2,2}$ regularity of $g_\delta$, this implies that $z, r \in L^{k,2}$. Since the Ricci curvature is to leading order the Laplacian of the metric in harmonic coordinates, it follows that $g \in L^{k+2,2}$. Further, we also now have, by Sobolev embedding, $|z|^2 \in L^{(3/2)-\mu}$, $\mu > 0$ and so have boosted the initial regularity $|z|^2 \in L^1$ to $|z|^2 \in L^{(3/2)-\mu}$. Continuing in the same way with $L^p$ elliptic estimates gives the $C^\infty$ regularity of $g_\delta$ and $\omega_\delta$. We refer to [An3, §4, §8] for further details if desired.

If $\phi_\delta$ is a real-analytic function of $s$, then $g_\delta$ is real-analytic; this follows again from elliptic regularity, since the coefficients in (3.18)-(3.19) are real-analytic functions of the unknown $g_\delta$, c.f. [Mo, Ch.6.7]. This proves the claim above.

We now derive uniform estimates on the behavior of $\{g_\delta\}$ as $\delta \to 0$. The key to this is the following simple Lemma.

**Lemma 3.4.** *Suppose the potential function $\omega_\delta$ in (3.20) satisfies (locally)*

$$||\omega_\delta||_{L^{k,p}} \leq C_{k,p}, \tag{3.21}$$

*where $C_{k,p}$ is independent of $\delta$. If $0 \leq k \leq 1$ and $1 \leq p \leq \infty$, then each summand $\tau_\delta$ and $\frac{\varepsilon}{12}s$ of $\omega_\delta$ in (3.20) also satisfies the bound (3.21). In particular, if $\omega_\delta$ is Lipschitz, with Lipschitz constant $C$, then so are $\tau_\delta$ and $\frac{\varepsilon s}{12}$, with the same Lipschitz constant.*

*Proof.* From the assumption above on the smoothing $\phi_\delta$, the functions $\tau_\delta$ in (3.17) and $s$ have the same sign everywhere. Thus the Lemma is obvious for $k = 0$. Next suppose $k = 1$. We have

$$\int |\nabla \omega_\delta|^p = \int (\tau_\delta' + \frac{\varepsilon}{12})^p |\nabla s|^p \leq C_{k,p}.$$

Again by assumption on $\tau_\delta, \tau_\delta' > 0$, for all $\delta > 0$, so that one obtains bounds on both summands $\int |\nabla \tau_\delta|^p$ and $\int |\nabla \frac{\varepsilon}{12} s|^p$. Note that this argument also remains valid if $p = \infty$.

The same argument also holds in $L^{k,p}$, for fractional values of $k < 1$; c.f. the references to [Ad] or [Tr] above. ∎

**Remark 3.5.** Lemma 3.4 is not true when $k > 1$, so that the conclusion that the summands are Lipschitz is optimal. In fact, as a simple model, one may take, (in 1 variable $t$),

$$s = t, \text{ for } t \leq 0 \text{ and } s = \frac{12}{\varepsilon}(\frac{1}{\sigma} + \frac{1}{12}\varepsilon)t, \text{ for } t > 0, \tag{3.22}$$

where $\sigma > 0$ is any fixed constant. For $\omega$ as in (3.13), (i.e. setting $\delta = 0$ in (3.20)), this gives $\omega = (\frac{1}{\sigma} + \frac{1}{12}\varepsilon)t \in C^\infty$, but $s$ itself is only Lipschitz. As $\varepsilon \to 0$ here, note that the derivative of $s$ blows up, as does $s$ itself, at any given $t > 0$.

This Lemma allows one to prove the stated regularity of the limit $g_\varepsilon = \lim_{\delta \to 0} g_\delta$ in Theorem 3.3 using the same method as following the equations (3.18)-(3.19).

**Proof of Theorem 3.3.**

The trace-free curvature $z$ and negative part of scalar curvature $s^-$ of $g_\delta$ are uniformly bounded in $L^2$ on $\Omega_{\varepsilon,\delta}$ as $\delta \to 0$; (recall throughout this discussion, $\varepsilon > 0$ is fixed). Hence, by the discussion in (I) above, a subsequence converges in the weak $L^{2,2}$ topology to a $L^{2,2}$ limit metric $g_\varepsilon$, defined on a maximal domain $\Omega_\varepsilon$. This means in particular that the metrics $g_\delta$ are uniformly bounded in $L^{2,2}$ in local harmonic coordinates in the following sense. Given any ball $B \subset (\Omega_\varepsilon, g_\varepsilon)$ of radius equal to the harmonic radius, the metrics $g_\delta$ may be viewed as metrics on $B$ and are uniformly bounded in $L^{2,2}$ in local $g_\delta$ harmonic coordinates on $B$. Further, the $g_\delta$ harmonic coordinates converge in $L^{3,2}$ to those of $g_\varepsilon$.

To prove uniform regularity, $\varepsilon|z|^2$ is uniformly bounded in $L^1$ as $\delta \to 0$. As before, by Sobolev embedding $L^1 \subset L^{t-2,2}$, for any given $t < \frac{1}{2}$; (recall again that we are working only locally). Regularity theory for the Laplacian $\Delta : L^{t,2} \to L^{t-2,2}$ in (3.19) implies that $\omega_\delta$ is uniformly bounded in $L^{t,2}$ as $\delta \to 0$. Now Lemma 3.4 implies that each summand $\tau_\delta$ and $\frac{\varepsilon}{12}s$ is uniformly bounded in $L^{t,2}$, and consequently $L^*\tau_\delta$ is bounded in $L^{t-2,2}$ as $\delta \to 0$. As before, consider the elliptic operator $\nabla \mathcal{Z}^2$ mapping $L^{t,2} \to L^{t-2,2}$. It follows from (3.18) and (3.8) that $\varepsilon z$, and hence $\varepsilon r$, is uniformly bounded in $L^{t,2} \subset L^{3-\mu}$, where the last inclusion follows from Sobolev embedding, $\mu = \mu(t) > 0$. This bound on $r$ implies a uniform $L^{2,3-\mu}$ bound on the metric $g_\delta$ in local $g_\delta$-harmonic coordinates as $\delta \to 0$. Further, $\varepsilon|z|^2$ is uniformly bounded in $L^{(3-\mu)/2}$, (since $\varepsilon$ is fixed). One may now repeat the analysis above, having improved the initial regularity on $\varepsilon|z|^2$, using $L^p$ elliptic estimates. This process may be continued as long as Lemma 3.4 is applicable.

In this way, one concludes that $\varepsilon|z|^2$ is uniformly bounded in $L^{1,p}$, and hence the function $\omega_\delta$ is uniformly bounded in $L^{3,p}$ as $\delta \to 0$, for any given $p < \infty$. However, the summands $\tau_\delta$ and $s_\delta$ are uniformly bounded only in $L^{1,\infty}$, i.e. in the Lipschitz topology. Similarly, the metrics $g_\delta$ remain uniformly bounded in the $L^{3,p}$ topology.

It follows that at the limit $\delta = 0$, $g_\varepsilon \in L^{3,p}, \omega = \lim_{\delta \to 0} \omega_\delta \in L^{3,p}$ while $\tau = \lim_{\delta \to 0} \tau_\delta$ and $s$ are (locally) Lipschitz. Since the metrics $g_\delta$ are smooth solutions of the Euler-Lagrange equations (3.18)-(3.19), $g_\varepsilon$ is a weak $L^{3,p}$ solution of its Euler-Lagrange equation (3.14). It is clear from the regularity arguments above that $g_\varepsilon$ is $C^\infty$, and in fact real-analytic, away from the set where $s$ vanishes. This completes the proof of Theorem 3.3.

■

Define the *junction set* of $(\Omega_\varepsilon, g_\varepsilon)$ to be the set

$$\Sigma = \partial\{s = 0\}, \tag{3.23}$$

where $s = s_\varepsilon$. The junction $\Sigma$ decomposes $\Omega_\varepsilon$ into a disjoint union

$$\Omega_\varepsilon = U^- \cup \Sigma \cup U^+, \tag{3.24}$$

where $U^- = \{x \in \Omega_\varepsilon : s(x) < 0\}$ and $U^+$ is the interior of the set $\{s \geq 0\}$, so that $\partial U^- = \partial U^+ = \Sigma$.

**Remark 3.6.** We conjecture that $\Sigma$ is a Lipschitz surface in $\Omega_\varepsilon$, of course possibly disconnected, and that $U^+ = \{s > 0\}$. In particular, this would imply that $\{s = 0\}$ has empty interior and equals $\Sigma$. On the other hand, apriori, the closed set $\Sigma$ could be quite complicated in $\Omega_\varepsilon$.

Observe that in the interior of $\{s = 0\}$, the (vanishing of the) trace equation (3.12) implies that $|z|^2 = const = 6c/\varepsilon$, for $c$ as in (3.11). It seems likely that this, together with the full minimizing property of $g_\varepsilon$ should imply that $g_\varepsilon$ is flat in this interior, (giving then a contradiction), but this remains an open question.

Note also that the proof of Theorem 3.3 shows that $g_\varepsilon$ is smooth, in fact real-analytic, on the domain $U^+ \cup U^-$, i.e. away from $\Sigma$.

Next consider the behavior of $\nabla\omega$ at $\Sigma$. By Theorem 3.3 and Sobolev embedding, $\omega \in C^{2,\alpha}$. Hence the gradient $\nabla\omega$ is a $C^{1,\alpha}$ vector field on $\Omega_\varepsilon$, and so is well-defined at every point in $\Omega_\varepsilon$. In $U^-$, we have

$$\nabla\omega = \nabla(\tau + \frac{\varepsilon}{12}s) = (\frac{1}{\sigma} + \frac{1}{12}\varepsilon)\nabla s,$$

while in $U^+$,

$$\nabla\omega = \frac{1}{12}\varepsilon\nabla s.$$

Since $\nabla\omega$ is continuous, if $\{x_i\}$ is any sequence in $U^-$ converging to a point $x \in \Sigma$, then $\nabla\omega(x_i) = (\frac{1}{\sigma} + \frac{1}{12}\varepsilon)\nabla s(x_i) \to \nabla\omega(x)$. Thus,

$$(\frac{1}{\sigma} + \frac{1}{12}\varepsilon)\nabla^- s(x) \equiv lim_{i\to\infty}(\frac{1}{\sigma} + \frac{1}{12}\varepsilon)\nabla s(x_i) = \nabla\omega(x),$$

is well-defined at any $x \in \Sigma$. Applying the same argument in $U^+$ shows that,

$$\frac{1}{12}\varepsilon\nabla^+ s(x) \equiv lim_{j\to\infty}\frac{1}{12}\varepsilon\nabla s(x_j) = \nabla\omega(x),$$

is well defined, for $x_j$ any sequence in $U^+$.

This proves the following result, compare with Remark 3.5 and the discussion in (1.12)-(1.15).

**Proposition 3.7.** *In the notation above, for any point $x \in \Sigma$, we have*

$$\nabla^-(\tau + \frac{\varepsilon}{12}s)(x) = \frac{\varepsilon}{12}\nabla^+ s(x). \tag{3.25}$$

∎

**(IV). Completeness.**

The proof that $(\Omega_\varepsilon, g_\varepsilon)$ is a complete Riemannian manifold is exactly the same as the proof of completeness of minimizers of $I'_\varepsilon$, c.f. (3.6), given in [An3,Thms. 5.4, 8.2], given the following two assertions.

(i). The scalar curvature $s_\varepsilon$ of $(\Omega_\varepsilon, g_\varepsilon)$ is uniformly bounded below, for any fixed $\varepsilon > 0$, i.e.

$$s_\varepsilon \geq \lambda > -\infty; \tag{3.26}$$

this is the exact analogue of [An3, (5.0)].

(ii). By the estimate (2.11), it follows from (i) that the $L^2$ curvature radii w.r.t. the curvatures $r$ and $z$ are equivalent.

The reason that the proof is the same is because, given (3.26), the proof of completeness of minimizers of $I'_\varepsilon$ relies completely on the behavior of the curvature radius $\rho$ near any metric boundary of $\Omega_\varepsilon$; on approach to $\partial\Omega_\varepsilon$ the metric must become 'cusp-like'. Now by (i) and (ii), together with the global statement (2.12), one obtains exactly the same control on $\rho$ for $I_\varepsilon^-$ minimizers as one does with $I'_\varepsilon$ minimizers; see also [An3, Remark 5.5] for a summary of the proof.

The following Lemma establishes (3.26), at least for some minimizing pairs $(\Omega_\varepsilon, g_\varepsilon)$.

**Lemma 3.8.** *For any $\varepsilon > 0$ fixed, there exist minimizing pairs $(\Omega_\varepsilon, g_\varepsilon)$ of $I_\varepsilon^-$ as in Theorem 3.3 such that (3.26) holds on $\Omega_\varepsilon$, for some $\lambda = \lambda(\varepsilon) < \infty$.*

**Proof.** Recall from the proof of Theorem 3.3 that $(\Omega_\varepsilon, g_\varepsilon)$ is obtained as a limit of minimizers $(\Omega_{\varepsilon,\delta}, g_{\varepsilon,\delta})$ of $I_{\varepsilon,\delta}$ in (3.15), with $\varepsilon > 0$ fixed and $\delta \to 0$. Fix $\delta > 0$ small, and let $\{g_j\}$ be a minimizing sequence of unit volume $C^\infty$ metrics on $M$ for $I_{\varepsilon,\delta}$, so that $g_j \to g_{\varepsilon,\delta}$ in $L^{2,2}$, as in (I) above. Let $\hat{g}_j$ be unit volume metrics in the conformal class $[g_j]$ realizing the infimum of $I_{\varepsilon,\delta}$ restricted to $[g_j]$. From the form of $\phi_\delta$ in (3.15), for any $\varepsilon > 0$ and $\delta > 0$ fixed, a bound on $I_{\varepsilon,\delta}$ gives a bound on the $L^2$ norm of the Ricci curvature, (depending on $\varepsilon/\delta$). Hence, Remark 2.2 implies that $\hat{g}_j$ exists. Further, and this is the main point, a result of Gursky [Gu] implies that $\hat{g}_j$ is a globally defined $L^{2,2}$ metric on $M$, (and not just on some $\Omega_{\varepsilon,\delta}(j)$). In addition, completely analogous to the regularity arguments following (3.20), the work in [An3, p.239-240] implies that $\hat{g}_j$ is $C^\infty$, since $g_j$ is.

The metric $\hat{g}_j$ minimizes $I_{\varepsilon,\delta}$ in its conformal class, and hence each $\hat{g}_j$ satisfies the trace Euler-Lagrange equation (3.19) on $M$, with

$$3c = \frac{\varepsilon}{2}\mathcal{Z}^2 + (\int_M s\tau_\delta) - \frac{3}{4}\sigma_\delta \tag{3.27}$$

where all the terms are evaluated at $\hat{g}_j$. We then evaluate (3.19) at a point $q \in M$ realizing the minimum of $\tau_\delta + \frac{\varepsilon}{12}s$. Notice that both $\tau_\delta$ and $s$ have the same sign, necessarily negative, at the minimum, and so $\tau_\delta(q) = min_M \tau_\delta$. This gives

$$(s\tau_\delta - \frac{3}{4}\frac{\phi_\delta^2}{\sigma_\delta})(q) \leq 3c. \tag{3.28}$$

Since $\phi_\delta(s) = s$ when $s \leq 0$, it follows from the definition (3.17) that

$$\frac{1}{2}s_{(\hat{g}_j)} \geq -(3c\sigma_\delta)^{1/2}. \tag{3.29}$$

Clearly, $c(\hat{g}_j)$ and $\sigma_\delta(\hat{g}_j)$ are uniformly bounded as $j \to \infty$.

To conclude the proof, we may carry out the process above for any fixed $\varepsilon > 0$ and any varying $\delta > 0$ tending to 0. Choosing a suitable diagonal subsequence then gives minimizers $\{\hat{g}_j\}$ of $I_{\varepsilon,\delta_j}$, with $\delta_j \to 0$. As in the proof of Theorem 3.3 above, it then follows that a subsequence converges to a limit minimizing pair $(\hat{\Omega}_\varepsilon, \hat{g}_\varepsilon)$ for $I_\varepsilon^-$, with the same regularity properties as in Theorem 3.3. Further, the constant $c(\hat{g}_j) \cdot \sigma_\delta(\hat{g}_j)$ is uniformly bounded as $\delta \to 0$, which gives the bound (3.26). ∎

For the remainder of the paper, we will always assume that $(\Omega_\varepsilon, g_\varepsilon)$ is obtained via Lemma 3.8, and so satisfies (3.26). We remove however the hat notation.

Summarizing, the following result expresses the geometric decomposition of $M$ w.r.t. a minimizer $g_\varepsilon$ of $I_\varepsilon^-$. c.f. again [An3, §5,§8] for the statements (3.31)-(3.32) and their proofs, which are identical for $I_\varepsilon'$ and $I_\varepsilon^-$.

**Theorem 3.9. (Geometric Decomposition for $I_\varepsilon^-$).** *Suppose $\sigma(M) \leq 0$. For any $\varepsilon > 0$, there is a complete, $L^{3,p} \cap C^{2,\alpha}$ Riemannian metric $g_\varepsilon$, defined on an open domain $\Omega_\varepsilon$, which satisfies (weakly) the Euler-Lagrange equations*

$$\nabla I_\varepsilon^- = 0, \tag{3.30}$$

*as in (3.9)-(3.12). The scalar curvature $s_\varepsilon$ and potential $\tau_\varepsilon = \lim_{\delta \to 0} \tau_{\varepsilon,\delta}$ are locally Lipschitz functions on $\Omega_\varepsilon$ and $g_\varepsilon$ is real-analytic on the complement of the junction $\Sigma = \partial\{s = 0\}$. The pair $(\Omega_\varepsilon, g_\varepsilon)$ realizes $inf_{\mathbb{M}_1} I_\varepsilon^-$ in the sense that*

$$I_\varepsilon^-(g_\varepsilon) = \varepsilon \int_{\Omega_\varepsilon} |z_{g_\varepsilon}|^2 + (\int_{\Omega_\varepsilon}(s_{g_\varepsilon}^-)^2)^{1/2} = inf_{\mathbb{M}_1} I_\varepsilon^-. \tag{3.31}$$

*Further, the curvature $r_{g_\varepsilon}$ of $g_\varepsilon$ is uniformly bounded on $\Omega_\varepsilon$ and $\Omega_\varepsilon$ consists of a finite number $n = n(\varepsilon, M)$ of components, with*

$$vol_{g_\varepsilon}\Omega_\varepsilon = \sum_{j=1}^N vol_{g_\varepsilon}\Omega_j = 1. \tag{3.32}$$

∎

As above in (I), any domain with smooth and compact closure in $\Omega_\varepsilon$ embeds as a smooth domain in $M$, and a neighborhood of infinity of $\Omega_\varepsilon$ is a graph manifold. Observe that (3.31) shows that $(\Omega_\varepsilon, g_\varepsilon)$ is a minimizer of $I_\varepsilon^-$ in the usual sense of the word, (even though $g_\varepsilon$ may not be a Riemannian metric on $M$). Further, from the construction in Lemma 3.8, the definition (3.17), and lower semi-continuity of the $L^2$ norm, we have

$$\int_{\Omega_\varepsilon} \tau_\varepsilon^2 dV \leq 1.$$

More importantly, (3.31)-(3.32) and lower semicontinuity of the functionals $\varepsilon \mathcal{Z}^2$ and $\mathcal{S}_-^2$ imply that $\varepsilon \mathcal{Z}^2(\hat{g}_j) \to \varepsilon \mathcal{Z}^2(g_\varepsilon)$ and $\sigma_{\delta_j}(\hat{g}_j) = (\int (\phi_{\delta_j}(s))^2)^{1/2}(\hat{g}_j) \to \sigma(g_\varepsilon)$, i.e. these values pass continuously to the limit. In particular, it follows that the constant $c$ for the Euler-Lagrange equations, obtained from the limit as $\delta \to 0$ of (3.27), is given as in (3.11) by

$$c = \frac{1}{6}(\varepsilon \mathcal{Z}^2(g_\varepsilon) + \frac{1}{2}\sigma(g_\varepsilon)). \tag{3.33}$$

It is worth pointing out that if a minimizer $g_\varepsilon$ of $I_\varepsilon^-$ is a (unit volume) constant curvature metric on $M$, then clearly $I_\varepsilon^-(g_\varepsilon) = I_0^-(g_\varepsilon)$. If there is then a sequence $g_{\varepsilon_i}$ of such metrics, with $\varepsilon_i \to 0$, and some $g_{\varepsilon_i}$ is not flat, it follows from the Mostow rigidity of constant curvature metrics that $g_{\varepsilon_i} = g_0$, for all $\varepsilon_i > 0$, so that the sequence $\{g_{\varepsilon_i}\}$, and in fact the full family $\{g_\varepsilon\}$ is constant, (modulo diffeomorphisms). Conversely, if $g_\varepsilon$ is flat for some $\varepsilon$, then all metrics $g_\varepsilon$ are flat and the family $g_\varepsilon$, (if non-empty, i.e. if not collapsed), varies in the moduli space of flat metrics. In particular, in either case thre is no degeneration of $\{g_\varepsilon\}$, in the sense of (0.6).

**Remark 3.10.** Next, we make several remarks on continuity and differentiability properties of the values $I_\varepsilon^-(g_\varepsilon)$ as $\varepsilon$ varies.

First, it is not known if the metrics $g_\varepsilon$ or their domains $\Omega_\varepsilon$ are unique. Some results on the moduli space of minimizers of $I_\varepsilon^-$ can be deduced from Theorem 3.9; in particular, the moduli space is compact in a natural sense, for any given $\varepsilon$.

It is obvious from Theorem 3.9 that $I_\varepsilon^-(g_\varepsilon)$ is a well-defined function of $\varepsilon$, independent of any choice of minimizer $g_\varepsilon$, (arising as in Lemma 3.8). Further, the constant $c$ in (3.33) is uniquely determined by $\varepsilon$ among the family of minimizers $g_\varepsilon$. Since $3c \neq I_\varepsilon^-(g_\varepsilon)$, it follows that both terms $\mathcal{Z}^2(g_\varepsilon)$ and $\sigma(g_\varepsilon)$ are well-defined functions of $\varepsilon$, denoted as $\mathcal{Z}^2(\varepsilon)$, $\sigma(\varepsilon)$.

Since limits of minimizers are minimizers, it is easily verified that $I_\varepsilon^-(g_\varepsilon)$ is continuous in $\varepsilon$. In addition, the minimizing property of $g_\varepsilon$ implies that $I_\varepsilon^-(g_\varepsilon)$ is strictly monotone increasing in $\varepsilon$, (unless $(M, g_\varepsilon)$ is flat).

Next we observe that $I_\varepsilon^-(g_\varepsilon)$ is differentiable in $\varepsilon$. To see this, for any fixed $\varepsilon > 0$ and any variable $\delta$, the minimizing property implies that $I_{\varepsilon+\delta}^-(g_{\varepsilon+\delta}) = (\varepsilon+\delta)\mathcal{Z}^2(\varepsilon+\delta) + \sigma(\varepsilon+\delta) \leq (\varepsilon+\delta)\mathcal{Z}^2(\varepsilon) + \sigma(\varepsilon)$. Hence

$$I_{\varepsilon+\delta}^-(g_{\varepsilon+\delta}) - I_\varepsilon^-(g_\varepsilon) \leq \delta \mathcal{Z}^2(\varepsilon), \tag{3.34}$$

and so $(\frac{d}{d\varepsilon})^+ I_\varepsilon^-(g_\varepsilon) \leq \mathcal{Z}^2(\varepsilon)$, where $+$ denotes the right-hand derivative. Replacing $\delta > 0$ by $-\delta$ in (3.34) gives the opposite inequality for the left-hand derivative, and so

$$\frac{d}{d\varepsilon} I_\varepsilon^-(g_\varepsilon) = \mathcal{Z}^2(\varepsilon). \tag{3.35}$$

Further, (3.34) implies that $I_\varepsilon^-(g_\varepsilon)$ is *concave* in $\varepsilon$, so that $\mathcal{Z}^2(\varepsilon)$ is a monotone non-increasing function of $\varepsilon$. On the other hand, one sees immediately from (3.2) that $\mathcal{Z}^2(\varepsilon)$ is a non-increasing function of $\sigma$. Hence $\sigma(\varepsilon) = \sigma(\mathcal{Z}^2(\varepsilon))$ is monotone non-decreasing in $\varepsilon$.

Finally, we observe that $\sigma(\varepsilon)$ is (strictly) monotone increasing in $\varepsilon$, unless $\sigma(\varepsilon)$ is identically constant in $\varepsilon$, (and hence $\sigma(\varepsilon) = |\sigma(M)|$, by (3.4)). For suppose there is an interval $[\varepsilon_1, \varepsilon_2]$ such that $\sigma(\varepsilon) = \sigma_o$ is constant on $[\varepsilon_1, \varepsilon_2]$, and $\sigma(\varepsilon) > \sigma_o$, for $\varepsilon > \varepsilon_2$. The minimizing property (3.2) implies that $\mathcal{Z}^2(\varepsilon) = z_o$ is constant on $[\varepsilon_1, \varepsilon_2]$ and $\mathcal{Z}^2(\varepsilon) \leq z_o$ for $\varepsilon \geq \varepsilon_2$. Hence $I_\varepsilon^-(g_\varepsilon)$ is affine on $[\varepsilon_1, \varepsilon_2]$. We then claim that for $\varepsilon > 0$ small,

$$I_{\varepsilon_1+\varepsilon}^-(g_{\varepsilon_2+\varepsilon}) < I_{\varepsilon_1+\varepsilon}^-(g_{\varepsilon_1+\varepsilon}).$$

To see this, both sides agree at $\varepsilon = 0$. The derivative of the right side w.r.t. $\varepsilon$ is $\mathcal{Z}^2(\varepsilon_1) = z_o$ for $\varepsilon \in [0, \varepsilon_2 - \varepsilon_1]$. This also equals the derivative of the left side when $\varepsilon = 0$. However since $I_\varepsilon^-(g_\varepsilon)$ is concave, the derivative of the left side is smaller than $z_o$, for $\varepsilon > 0$ small, which gives the claim. This claim however contradicts the minimizing property, and hence proves the observation.

**§3.2.** We now derive estimates on the behavior of the scalar curvature $s_\varepsilon$ and its normalization $\tau_\varepsilon$ for the metric $(\Omega_\varepsilon, g_\varepsilon)$, as $\varepsilon$ tends to 0. It turns out we need to separate the cases $\sigma(M) < 0$ and $\sigma(M) = 0$.

**Theorem 3.11.** *Suppose $\sigma(M) < 0$. Then for all $\varepsilon > 0$, the potential function $\tau_\varepsilon = \lim_{\delta \to 0} \tau_{\varepsilon,\delta}$ satisfies*

$$\tau_\varepsilon = \frac{s_\varepsilon^-}{\sigma(g_\varepsilon)}, \quad \text{and} \quad \int_{\Omega_\varepsilon} \tau_\varepsilon^2 = 1. \tag{3.36}$$

*As $\varepsilon \to 0$, we have*

$$inf_{\Omega_\varepsilon} \tau_\varepsilon \to -1, \tag{3.37}$$

*and, for any $p < \infty$,*

$$\int_{\Omega_\varepsilon} |\tau_\varepsilon + 1|^p dV_\varepsilon \to 0. \tag{3.38}$$

*as well as*

$$\int_{\Omega_\varepsilon} |\nabla \tau_\varepsilon|^2 dV_\varepsilon \to 0. \tag{3.39}$$

*Proof.* Recall from the construction of $(\Omega_\varepsilon, g_\varepsilon)$ in Theorem 3.9 that the constants $c(\hat{g}_j)$ in (3.27) converge to the constant $c$ in (3.33) and that $\sigma_{\delta_j}(\hat{g}_j) \to \sigma(g_\varepsilon)$ as preceding (3.33). Since $\sigma(g_\varepsilon) \geq |\sigma(M)| > 0$, and since $\phi_\delta(s) \to s^-$, it follows from the definition (3.17) that the limit $\tau_\varepsilon$ is given by (3.36). Further, the second equation in (3.36) follows, since it holds for each $\tau_{\varepsilon,\delta}$, and the $L^2$ norm passes continuously to the limit, (as preceding (3.33) again).

The estimate (3.37) now follows essentially from Lemma 3.8. Namely, from (3.28), we obtain on the limit $(\Omega_\varepsilon, g_\varepsilon)$,

$$\frac{1}{4} inf s_\varepsilon \cdot inf \tau_\varepsilon \leq \frac{\varepsilon}{2} \mathcal{Z}^2(g_\varepsilon) + \frac{1}{4} \sigma(g_\varepsilon).$$

From the first equation in (3.36), this gives

$$(inf \tau_\varepsilon)^2 \leq 1 + \frac{2\varepsilon \mathcal{Z}^2}{\sigma}(g_\varepsilon).$$

Since the right hand side converges to 1 as $\varepsilon \to 0$, by (3.3) and (3.4), this gives (3.37).

The estimate (3.38) is now essentially obvious. We have $\tau_\varepsilon \leq 0$, with $L^2$ average equal to 1 and by (3.37), inf $\tau_\varepsilon \geq -1 - \delta$, where $\delta \to 0$ as $\varepsilon \to 0$. It follows that $\tau_\varepsilon$ must converge to $-1$ almost everywhere as $\varepsilon \to 0$, and since $\tau_\varepsilon$ is bounded, (3.38) follows.

To obtain the estimate on the $L^2$ norm of $\nabla \tau_\varepsilon$, the trace equation for $(\Omega_\varepsilon, g_\varepsilon)$ is given by

$$2\Delta(\tau + \frac{\varepsilon}{12} s) + \frac{1}{4} s\tau = -\frac{\varepsilon}{2} |z|^2 + \frac{\varepsilon}{2} \int |z|^2 dV + \frac{1}{4\sigma} \int (s^-)^2 dV. \tag{3.40}$$

Multiply this by $\tau \leq 0$ and integrate by parts over $(\Omega_\varepsilon, g_\varepsilon)$. Since $\tau$ is uniformly bounded and $(\Omega_\varepsilon, g_\varepsilon)$ is of uniformly bounded curvature and finite volume, a standard cutoff argument near infinity shows the boundary term vanishes at infinity, and hence

$$2 \int |\nabla \tau|^2 dV \leq -(\frac{\varepsilon}{2} \int |z|^2 dV)(\int \tau dV) + \frac{\sigma}{4} \int |\tau| dV - \frac{1}{4} \int |s| \tau^2 dV,$$

The first term on the right here clearly goes to 0, as $\varepsilon \to 0$. A simple use of the Hölder inequality implies that the sum of the last two terms on the right is non-positive. This gives the estimate (3.39). ∎

We now turn to the case $\sigma(M) = 0$. The proof of (3.37) carries over to this situation if one can prove $\sigma(g_\varepsilon) > 0$ and

$$\frac{\varepsilon}{\sigma(g_\varepsilon)} \mathcal{Z}^2(g_\varepsilon) \to 0, \quad \text{as} \ \varepsilon \to 0. \tag{3.41}$$

In fact, if

$$T = T_\varepsilon = |inf_{\Omega_\varepsilon} \tau_\varepsilon|, \tag{3.42}$$

then this proof shows

$$1 \le T \le (1 + \frac{2\varepsilon}{\sigma(g_\varepsilon)} \mathcal{Z}^2(g_\varepsilon))^{1/2}. \tag{3.43}$$

Whether (3.41) holds or not turns seems to be a delicate question. In case $\sigma(M) < 0$, (3.41) is analogous to

$$\frac{\varepsilon \mathcal{Z}^2(g_\varepsilon)}{\sigma(g_\varepsilon) - |\sigma(M)|} \to 0, \quad \text{as} \ \varepsilon \to 0,$$

which is almost certainly *not* true; we do not discuss this explicitly, but refer to the examples in §6.1 for some indication of the reasons, c.f. in particular (6.7)-(6.10) Thus, any proof of (3.41) probably must use the fact that $\sigma(M) = 0$.

This situation illustrates one of the differences between the cases $\sigma(M) = 0$ and $\sigma(M) < 0$; we will see later other differences, mostly related to collapse, c.f. for example §6. Note that when $\sigma(M) < 0$, $\tau$ and the negative part $s^-$ of the scalar curvature are essentially the same - their ratio is a constant bounded away from 0 and $\infty$. When $\sigma(M) = 0$, the ratio $\tau/s^-$ becomes unbounded as $\varepsilon \to 0$, and thus one would expect it is more difficult to control $\tau$.

We settle with the following analogue of Theorem 3.10 in the case $\sigma(M) = 0$.

**Theorem 3.12.** *Suppose $\sigma(M) = 0$. Then, for all $\varepsilon > 0$,*

$$\sigma(g_\varepsilon) > 0, \tag{3.44}$$

*and, with $\tau_\varepsilon = lim_{\delta \to 0} \tau_{\varepsilon,\delta}$, (3.36) holds. Further, as $\varepsilon \to 0$, we have*

$$inf_{\Omega_\varepsilon} s_\varepsilon \to 0, \tag{3.45}$$

*and*

$$\int_{\Omega_\varepsilon} |\nabla \tau_\varepsilon|^2 dV_\varepsilon \to 0. \tag{3.46}$$

*Proof.* Given (3.44), the proof of (3.36) in Theorem 3.11 holds without any changes here. (Note this argument does not involve letting $\varepsilon \to 0$). Similarly, the proof of (3.37) gives then, as before,

$$inf_{\Omega_\varepsilon} \tau_\varepsilon \ge -\big(1 + \frac{2\varepsilon}{\sigma(g_\varepsilon)} \mathcal{Z}^2(g_\varepsilon)\big)^{1/2},$$

as in (3.43). Since $\tau_\varepsilon = s_\varepsilon^-/\sigma(g_\varepsilon)$, this gives

$$inf_{\Omega_\varepsilon} s_\varepsilon \ge -\sigma(g_\varepsilon)\big(1 + \frac{2\varepsilon}{\sigma(g_\varepsilon)} \mathcal{Z}^2(g_\varepsilon)\big)^{1/2} \to 0, \ \text{as} \ \varepsilon \to 0,$$

by (3.3) and (3.4). Similarly, the proof of (3.46) is the same as that of (3.39).

Thus, it suffices to prove (3.44). By the discussion at the end of Remark 3.10, if $\sigma(g_\varepsilon) = 0$ for some $\varepsilon > 0$, then $\sigma(g_\varepsilon) \equiv 0$, for all $\varepsilon > 0$. As discussed then preceding Remark 3.2, this means that the metrics $g_\varepsilon$ minimize $\mathcal{Z}^2$ among all unit volume metrics on $M$, with no constraint on the scalar curvature. In particular, $(\Omega_\varepsilon, g_\varepsilon)$ satisfies the Euler-Lagrange equations (3.30) with $\tau \equiv 0$ and $s_\varepsilon \geq 0$ everywhere.

Now if $\Omega_\varepsilon = M$, i.e. $\Omega_\varepsilon$ is compact, then a standard result, c.f. [B, 4.49] implies that either $\sigma(M) > 0$, (i.e. $M$ admits a metric of positive scalar curvature), a contradiction to (0.3), or $(M, g_\varepsilon)$ is flat, for all $\varepsilon$. But this implies that $M$ is a graph manifold, which has been ruled out previously.

A similar argument proves that $(\Omega_\varepsilon, g_\varepsilon)$ is flat, even when $\Omega_\varepsilon$ is non-compact, again giving a contradiction, since a complete, open flat manifold has infinite volume. The proof of this, although conceptually similar, is more complicated and relies on work developed later in the paper. Thus, to be conceptually reasonably coherent, we defer the proof of this last situation to Proposition C.3 in Appendix C.

∎

While we certainly conjecture that there is a uniform lower bound for $\tau_\varepsilon$ as $\varepsilon \to 0$ in case $\sigma(M) = 0$, (at least on some sequence $\varepsilon_j \to 0$), we have not been able to prove this. It is not difficult to show that if there is no such bound, then the curvature $\mathcal{Z}^2(g_\varepsilon)$ blows up only very slowly, in that $\mathcal{Z}^2(g_\varepsilon) \leq \varepsilon^{-\delta}$, for any fixed $\delta > 0$ as $\varepsilon \to 0$. Such a slow (global) curvature blow-up would seem unlikely if the family $\{g_\varepsilon\}$ degenerates in the sense of (0.6).

While a uniform upper bound on $T$ in (3.42) in case $\sigma(M) = 0$ is certainly useful, it is not essential. In the work to follow in later sections, we will renormalize $\tau_\varepsilon$ to $\tau_\varepsilon/T_\varepsilon$, and thus, by fiat, have $\tau/T$ bounded below.

We close this section with a discussion on the relation between the full curvature and the potential function $\tau_\varepsilon$ of $g_\varepsilon$. We begin with the following local analogue of Lemma 1.3.

**Proposition 3.13.** *Let $x_\varepsilon \in (\Omega_\varepsilon, g_\varepsilon), \varepsilon = \varepsilon_i \to 0$, be a sequence of points such that*

$$\rho(x_\varepsilon) \geq \rho_o \quad \text{and} \quad \nu(x_\varepsilon) \geq \nu_o \cdot \rho_o, \tag{3.47}$$

*for some arbitrary constants $\rho_o, \nu_o > 0$.*

*If $\sigma(M) < 0$, then a subsequence of $(B_{x_\varepsilon}(\rho_o), g_\varepsilon)$ converges to a constant curvature metric $g_{g_o}$ on $B_x(\rho_o)$ with $s_{g_o} = \sigma(M)$ and $\tau_\varepsilon \to -1$ on $B_x(\rho_o)$. If $\sigma(M) = 0$, and the potential $\tau_\varepsilon$ does not converge to the 0 function in $L^2(B_{x_\varepsilon}(\rho_o))$, then the same conclusion holds for $\{g_\varepsilon\}$ and $\tau_\varepsilon \to const. < 0$.*

*Proof.* The bounds on $\rho(x_\varepsilon)$ and the volume radius $\nu(x_\varepsilon)$ imply that a subsequence converges in the weak $L^{2,2}$ topology to a limit metric $g_o$ on $B_x(\rho_o), x = \lim x_\varepsilon$, (c.f. §2). Theorems 3.11 and 3.12 imply that $\tau$ (sub)-converges to a constant function $\tau_o$ in $L^{1,2}(B_{x_\varepsilon}(\rho_o))$, and in case $\sigma(M) < 0$, (3.38) implies that $\tau_o$ is necessarily $-1$.

The limit metric $g_o$ is a $L^{2,2}$ weak solution of the Euler-Lagrange equations (3.9)-(3.12) with $\varepsilon$ set to 0 and $s_{g_o} = \sigma(M)$. Hence from (3.9) and (3.30), we obtain

$$L^* \tau_o + \sigma(M)(\frac{\tau_o}{4} - \frac{1}{12}) \cdot g = 0. \tag{3.48}$$

If $\sigma(M) < 0$, then the formula (1.6) for $L^*$ and (3.48) give $z = 0$ weakly in $L^2$. Similarly, if $\sigma(M) = 0$ and $\tau_o \neq 0$, then again $z = 0$ in $L^2$. Standard elliptic regularity for $L^{2,2}$ weak solutions to the Einstein equation $z = 0$ implies that $g_o$ is smooth.

∎

It is not known if the assumption that $\tau_o \neq 0$ can be removed in case $\sigma(M) = 0$, although it can be replaced by other assumptions. For instance, we note without proof that if $\lambda_1$, the lowest (non-zero) Laplace eigenvalue of $(\Omega_\varepsilon, g_\varepsilon)$, is bounded away from 0 as $\varepsilon \to 0$, then $\tau_o$ must be $-1$.

**Remark 3.14.** As noted in §0, the main interest in the work to follow is the situation where $(\Omega_\varepsilon, g_\varepsilon)$ degenerates at some sequence $x_\varepsilon \in (\Omega_\varepsilon, g_\varepsilon)$, i.e.

$$\rho(x_\varepsilon) \to 0, \quad \text{as} \ \ \varepsilon \to 0. \tag{3.49}$$

As a kind of converse to Proposition 3.13, suppose $\sigma(M) < 0$, and $q_\varepsilon \in (\Omega_\varepsilon, g_\varepsilon)$ is any sequence of points such that

$$\tau(q_\varepsilon) \geq -1 + \delta, \tag{3.50}$$

for any fixed $\delta > 0$, so that $\tau(q_\varepsilon)$ is bounded away from $\inf_{\Omega_\varepsilon} \tau$. Suppose as above that $\{g_\varepsilon\}$ is non-collapsing at $q_\varepsilon$, i.e. $\nu(x_\varepsilon) \geq \nu_o \rho(x_\varepsilon)$.

Then we claim that the metrics $g_\varepsilon$ *must* degenerate at $q_\varepsilon$ in the sense of (3.49). For if there were $\rho_o > 0$ such that $\rho(q_\varepsilon) \geq \rho_o$, then vol $B_{q_\varepsilon}(\rho_o) \geq \nu_o \cdot \rho_o$. By (3.38)-(3.39), $\tau$ converges to the constant function $-1$ in $L^{1,2}(B_{q_\varepsilon}(\rho_o))$. Further, as we will see later in Theorem 4.2, elliptic regularity for the Euler-Lagrange equations as $\varepsilon \to 0$ implies that then $\tau \to -1$ in $L^{2,2}(B_{q_\varepsilon}(\rho_o))$. Sobolev embedding now implies that $\tau \to -1$ in $C^{1/2}$, which contradicts (3.50).

A similar argument holds in case $\sigma(M) = 0$, at least if (for instance) $\lambda_1$ is bounded away from 0 as above.

Since $\delta$ in (3.50) is arbitrary, one may also choose a sequence $\delta_j \to 0$, sufficiently slowly as $\varepsilon \to 0$, so that there exist points $q_\varepsilon$ with $\tau(q_\varepsilon) \to -1$ and satisfying (3.49), (with $q_\varepsilon$ in place of $x_\varepsilon$), if there is no collapse at $\{q_\varepsilon\}$.

**Remark 3.15.** We make an important conceptual remark here, which will proved however only later in Remark 7.3(ii), since an understanding of blow-up behavior is required for the proof.

Proposition 3.13 implies that a bound on $\rho$ and $\nu$ implies control on the potential $\tau$ and on the metric $g_\varepsilon$. Conversely, we claim that suitable control on the potential $\tau$ also gives control on $\rho$ and hence on the metric $g_\varepsilon$.

Thus, let $D_\varepsilon \subset \Omega_\varepsilon$ be a region in which the oscillation of $\tau$ is small compared with its value at an interior point, i.e.

$$|osc_{D_\varepsilon}\tau| \leq \delta_o \cdot |\tau(y_\varepsilon)|, \tag{3.51}$$

where $\delta_o$ is a fixed sufficiently small constant, (independent of $\varepsilon$) and $t(y_\varepsilon) = dist(y_\varepsilon, \partial D_\varepsilon) > 0$. Assume also that $\tau(y_\varepsilon) < 0$. Then there is a constant $K < \infty$, independent of $\varepsilon$, such that

$$|r_{g_\varepsilon}|(y_\varepsilon) \leq K/t^2(y_\varepsilon). \tag{3.52}$$

The estimate (3.52) implies that $\rho(y_\varepsilon) \geq \kappa \cdot t(y_\varepsilon)$, for some constant $\kappa = \kappa(K) > 0$.

Thus the curvature of $g_\varepsilon$ can become very large in the region $\tau < 0$ only if $\tau$ oscillates by a definite amount in very small balls, (compare with [AnI, Thm. 3.3]).

**Remark 3.16.** It is useful to briefly compare the situation here with that in [AnI]. Recall that in the main result, Theorem A of [AnI], it was necessary to restrict the degeneration of $(M, g_i)$ to the region $U_{\delta_o}$ where the ratio $u/sup\,u$ is bounded below by $\delta_o$, for some $\delta_o > 0$. Degenerations in the complement of $U_\delta$, for any $\delta > 0$, may give rise only to super-trivial solutions of the static vacuum equations, (i.e. $u \equiv 0$), as seen in [AnI, §6]. The need for this restriction is due to the fact that the potential function $u$ in general may not approach its mean value in $L^2$, compare with [AnI, §4]. The function $\tau = s^-/\sigma$ in this paper plays essentially the same role as the function $u$ in [AnI], c.f. also Appendix A.

Thus, in effect, Theorem 3.11 allows one to remove the restriction to $U_{\delta_o}$, when $\sigma(M) < 0$, as does Theorem 3.12 in case $\sigma(M) = 0$, at least when $\lambda_1$ is bounded away from 0. Further the assumption on the $L^2$ bound of $z^T$ in [AnI,Thm. A] plays no role here or in Theorem B.

In fact, we will see below that the precise choice of the functional $I_\varepsilon^-$ with potential $\tau \leq 0$, allows one to avoid super-trivial solutions of the static vacuum equations in all circumstances. Thus the involved descent procedure to construct non-trivial blow-up limits in [An1, Thm.3.10] is not necessary in this paper.

## 4. Blow-up Limits of the Minimizers.

The results of §3 establish the basic geometric properties of the minimizers $(\Omega_\varepsilon, g_\varepsilon)$. In this section, we study the initial geometry of blow-up limits of sequences $\{g_\varepsilon\}$ at points where $\{g_\varepsilon\}$ is degenerating.

**§4.1.** We begin by considering the scaling behavior and the form of the Euler-Lagrange equations (3.9)-(3.12) under blow-up limits. We emphasize that most of the considerations in §4.1 are formal, and will be made more precise in the following subsection.

Thus, suppose there are points $x_\varepsilon \in \Omega_\varepsilon$ such that $\{g_\varepsilon\}$ degenerates at $x_\varepsilon$, i.e. $\rho(x_\varepsilon) = \rho(x_\varepsilon, g_\varepsilon) \to 0$, where $\rho$ is the $L^2$ curvature radius. As in (0.8), we rescale the metrics $g_\varepsilon$ based at $x_\varepsilon$ by the factor $\rho(x_\varepsilon)^{-2}$, i.e. set

$$g_\varepsilon' = \rho(x_\varepsilon)^{-2} \cdot g_\varepsilon. \tag{4.1}$$

The metric $g_\varepsilon'$ has $\rho_\varepsilon'(x_\varepsilon) = 1$, so that from Theorem 2.3, the sequence $(B_{x_\varepsilon}'(1), g_\varepsilon', x_\varepsilon)$ either has a convergent subsequence, modulo diffeomorphisms, or it collapses along a sequence of injective F-structures on $B_{x_\varepsilon}'(1 - \delta_\varepsilon)$. In the latter case, using Remark 2.1 we may pass to sufficiently large covers, e.g. the universal cover, to unwrap the collapse. Thus, we may suppose that $\{g_\varepsilon'\}$ sub-converges in the weak $L^{2,2}$ topology to a limit $(B_x'(1), g', x)$. Since the metrics $g_\varepsilon$, and thus also the metrics $g_\varepsilon'$ satisfy Euler-Lagrange equations, the limit $g'$ must also satisfy some elliptic P.D.E.

To determine the equation for the limit, return to the Euler-Lagrange equation (3.9)-(3.12), (3.30), valid for the metric $g_\varepsilon$. This equation is equally valid for the metric $g_\varepsilon'$ (the functional is scale-invariant). However, the volume term $v = vol_{g_\varepsilon}(M) = 1$ becomes $v' = vol_{g_\varepsilon'}(M) = \rho^{-3}$, while the scalar curvature $s_\varepsilon' = \rho^2 s_\varepsilon$, so that $\tau_\varepsilon' = \rho^2 \tau_\varepsilon$; here $\rho = \rho(x_\varepsilon)$. Thus, by Theorems 3.11-3.12, for most choices of $x_\varepsilon$, and for all choices when $\sigma(M) < 0$, $\tau_\varepsilon' \to 0$, as $\varepsilon \to 0$.

Now it is important to keep the potential function $\tau$ scale invariant; compare with Remark 3.16. Thus, $\tau$ is considered as the $L^2$ normalized (cutoff) scalar curvature *function* on the base $(\Omega_\varepsilon, g_\varepsilon)$ viewed in the scale $g_\varepsilon'$, and not as the normalized scalar curvature of $g_\varepsilon'$. Hence, in the scale $g_\varepsilon'$, we divide the Euler-Lagrange equations by $\rho^2$, and obtain the equations

$$\frac{\varepsilon}{\rho^2} \nabla \mathcal{Z}^2 + L^*(\tau) + (\frac{1}{4} s\tau + c_\varepsilon) \cdot g = 0, \tag{4.2}$$

$$2\Delta(\tau + \frac{\varepsilon s}{12\rho^2}) + \frac{1}{4} s\tau = -\frac{1}{2} \frac{\varepsilon}{\rho^2} |z|^2 + 3c_\varepsilon, \tag{4.3}$$

on $(\Omega_\varepsilon, g_\varepsilon')$. To simplify notation, we have dropped the prime and subscript $\varepsilon$; all metric quantities except $\rho = \rho(x_\varepsilon, g_\varepsilon)$ are taken w.r.t. the metric $g_\varepsilon'$. Thus for example $\tau$ is the same as the unscaled function, but $s = s_\varepsilon' = \rho^2 s_\varepsilon$. The constant $c_\varepsilon$ is given by

$$c_\varepsilon = \rho^2 \left( \frac{1}{12\sigma} \int (s^-)^2 dV + \frac{\varepsilon}{6} \int |z|^2 dV \right), \tag{4.4}$$

where all terms in (4.4) are taken w.r.t. $(\Omega_\varepsilon, g_\varepsilon)$, c.f. (3.33).

Without this renormalization by $\rho^2$, the Euler-Lagrange equations become $0 = 0$ in the limit, so that one has no information; this is analogous, (although only superficially), to the super-trivial solutions in [AnI]. (The factor $\rho^2$ appears essentially since we are dealing with the $L^2$ norm of the scalar curvature, and not its $L^1$ norm).

Hence, in passing to the limit $\varepsilon \to 0$, (in a subsequence), there are three possible behaviors corresponding to the behavior of the term $\varepsilon/\rho^2$, namely:

$$
\begin{aligned}
(i) &\qquad \frac{\varepsilon}{\rho^2} \to \infty, \\
(ii) &\qquad \frac{\varepsilon}{\rho^2} \to \alpha > 0, \\
(iii) &\qquad \frac{\varepsilon}{\rho^2} \to 0.
\end{aligned}
\tag{4.5}
$$

These cases correspond to the possible rates of curvature concentration (or blow-up) near $x_\varepsilon$, as $\varepsilon \to 0$, with (i) corresponding to the fastest concentration and (iii) to the slowest. Note that in general, as $\varepsilon \to 0$, in the scale $g'_\varepsilon, c_\varepsilon \to 0$ and $s_\varepsilon^- \to 0$, (c.f. §3.2), but $s_\varepsilon$, (i.e. the scalar curvature of $g'_\varepsilon$), may or may not go to 0.

Apriori, which of these cases occur depends on several factors. First, it depends on the choice of base points $\{x_\varepsilon\}$ at which the blow-up is centered. Different choices of $\{x_\varepsilon\}$ may lead to different behaviors (i)-(iii). Besides this, the occurence of (i)-(iii) may well depend on the topology of the underlying 3-manifold $M$. Since no assumptions are made on the topology of $M$, one has to allow for all three possibilities.

We discuss each of these cases in turn, formally for the moment; a more precise treatment follows in §4.2. We recall that from §3, the potential $\tau \leq 0$ is uniformly bounded below as $\varepsilon \to 0$ when $\sigma(M) < 0$. If $\sigma(M) = 0$, this remains unknown in general. In this case, the equations (4.2)-(4.3) must be further divided by $T = T_\varepsilon = \inf|\tau_\varepsilon|$, thus replacing the potential $\tau_\varepsilon$ by $\tau_\varepsilon/T_\varepsilon$ and replacing $\varepsilon$ by $\varepsilon/T_\varepsilon$ and $c_\varepsilon$ by $c_\varepsilon/T$. We assume implicitly that this has been done throughout this and the following sections.

**Case (i).** $\varepsilon/\rho^2 \to \infty$.

Renormalize the equations (4.2)-(4.3) by multiplying each term by $\rho^2/\varepsilon$. Then it follows that the limit metric $g'$, defined on a 3-manifold $N$, satisfies

$$
\nabla \mathcal{Z}^2 = 0,
\tag{4.6}
$$

$$
\Delta s = -3|z|^2.
$$

We call these the $\mathcal{Z}^2$ equations. These equations correspond to the case of a limit metric $g$ which is a critical point of $\mathcal{Z}^2 = \int |z|^2$ among all compact perturbations of $g$; there are no constraints on the scalar curvature or volume. Observe that from the discussion above, $s \geq 0$.

**Case (ii).** $\varepsilon/\rho^2 \to \alpha > 0$.

In this case, the limit equations on $(N, g')$ are

$$
\alpha \nabla \mathcal{Z}^2 + L^* \tau = 0,
\tag{4.7}
$$

$$
\Delta(\tau + \frac{1}{12}\alpha s) = -\frac{1}{4}\alpha|z|^2.
$$

We call these the scalar curvature-constrained $\mathcal{Z}^2$ equations, or $\mathcal{Z}_c^2$ equations for short. Observe that on the limit $(N, g')$, $s \geq 0$, $\tau \leq 0$, and

$$
s \cdot \tau \equiv 0.
\tag{4.8}
$$

These equations of course closely resemble the Euler-Lagrange equations (3.9)-(3.12), (3.30) for $g_\varepsilon$, with $\varepsilon$ replaced by the fixed number $\alpha$, but without lower order terms, i.e. $c = 0$ and $\chi = 0$ here.

The equations (4.7) are the equations for a metric minimizing $\mathcal{Z}^2$, (or a critical point of $\mathcal{Z}^2$), subject to a volume and scalar curvature constraint. Thus, $g' = \lim g'_\varepsilon$, and as noted in §3.1, $g'_\varepsilon$ minimizes the $L^2$ norm of $z$ among all metrics with volume and $L^2$ norm of $s^-$ at most that of $g'_\varepsilon$. Consequently, the limit $(N, g')$ minimizes the $L^2$ norm of the trace-free curvature among all compact perturbations $\bar{g}$ of $g'$ on $N$ satisfying $s_{\bar{g}} \geq 0$ and $vol_{\bar{g}} K \leq vol_g K$, for some compact set $K \subset N$ containing the perturbation $\bar{g}$ of $g'$.

In practice, $(N, g')$ will usually admit such comparison metrics, although this need not automatically always be the case. For example, the postive mass theorem [SY1] implies that $\mathbb{R}^3$ with the flat metric admits no compact perturbations with $s \geq 0$. On the other hand, the existence of scalar-flat comparisons would be much harder to establish, since it is not clear if there are any compact scalar-flat perturbations. This is one of the reasons for prefering $I_\varepsilon{}^-$ to $I_\varepsilon$ : compare with the discussion in §1(iii).

If $\tau \equiv 0$ on the limit, or a region in the limit, then the equations (4.7) become the $\mathcal{Z}^2$ equations (4.6), by dividing by $\alpha$. On the other hand, if $s \equiv 0$ on the limit, or a portion of it, then the equations (4.7) become

$$\alpha \nabla \mathcal{Z}^2 + L^* \tau = 0, \tag{4.9}$$

$$\Delta \tau = -\frac{1}{4} \alpha |z|^2.$$

We call the equations (4.9) the 'pure' $\mathcal{Z}_c^2$ equations, or the $\mathcal{Z}_s^2$ equations.

**Case (iii).** $\varepsilon/\rho^2 \to 0$.

In this case, the limit satisfies the equations

$$L^* \tau = 0, \tag{4.10}$$

$$\Delta \tau = 0,$$

i.e. the static vacuum Einstein equations with potential $\tau$.

Similar to the situation in (3.23), there may be a junction $\Sigma$ between different solutions of the types above on different regions of the limit $(N, g')$. As in (3.24), one may decompose $N$ into a disjoint union

$$N = N^- \cup \Sigma \cup N^+, \tag{4.11}$$

where $N^- = \{\tau < 0\}$, $N^+$ is the interior of the set $\{s \geq 0\}$ and $\Sigma = \partial(s = 0)$. Thus on $N^-$, the metric is either a static vacuum solution or $\mathcal{Z}_s^2$ solution, on $N^+$ the metric is a $\mathcal{Z}^2$ solution, and $\Sigma$ is a junction or connection between these two types of solution. Observe that (4.8) holds in all cases.

One might compare such a junction $\Sigma$ with junctions commonly employed in general relativity where vaccum regions of space or space-time are connected or joined to non-vacuum regions containing a non-trivial matter distribution, as for instance the junction at the surface of a star c.f. [Wd, Ch.6.2], or [MTW, Ch. 21.13, 23].

In Case (ii) above, the general $\mathcal{Z}_c^2$ solution is a junction between a $\mathcal{Z}_s^2$ solution on the "left" and a $\mathcal{Z}^2$ solution on the "right", (in terms of $\tau < 0$ on the left and $s > 0$ on the right). The $c$-notation is meant to suggest both the constraint and the existence of a connection, i.e. $\Sigma$, between these two types of equations. In case the junction or connection $\Sigma$ is empty, one has a 'pure' $\mathcal{Z}_s^2$ solution, and the $s$-notation is meant to represent scalar-flat.

It is not surprising that solutions of the static vacuum equations might arise as blow-up limits, given the work in [AnI], c.f. also Appendix A. The $\mathcal{Z}^2$ and $\mathcal{Z}_c^2$ solutions arise from the form of the perturbation of the functional $\mathcal{S}_-^2$. In terms of the two parts $\varepsilon \mathcal{Z}^2$ and $\mathcal{S}_-^2$ of $I_\varepsilon{}^-$, the $\mathcal{Z}^2$ equations arise when the curvature is concentrating so fast as $\varepsilon \to 0$ that the $\varepsilon \mathcal{Z}^2$ term is the dominant term (locally) in $I_\varepsilon{}^-$, the static vacuum equations arise when the curvature is concentrating so slowly that the $\mathcal{S}_-^2$ term is dominant, while the $\mathcal{Z}_s^2$ equations arise when the two terms are roughly balanced.

**§4.2.** We now make the formal arguments above rigorous. The discussion below parallels that in §4.1, and constitutes the proof of Proposition 4.4 and Theorem 4.9, stated later. Theorem 4.9 can be viewed as a preliminary version of Theorem B.

Throughout this section, we suppose the sequence $(\Omega_\varepsilon, g_\varepsilon)$ degenerates, in the sense that

$$\int_{\Omega_\varepsilon} |z|^2 dV_{g_\varepsilon} \to \infty, \ \text{ as } \ \varepsilon \to 0, \tag{4.12}$$

and so the $L^2$ norm of the full curvature $r$ also diverges. It is easy to see, (c.f. [AnI, Lem. 1.4] for a proof), that (4.12) implies the existence of points $x_\varepsilon \in (\Omega_\varepsilon, g_\varepsilon)$ such that

$$\rho(x_\varepsilon) \to 0, \ \text{ as } \ \varepsilon \to 0. \tag{4.13}$$

The set of points $\{x_\varepsilon\}$ for which (4.13) holds could be very complicated in $(\Omega_\varepsilon, g_\varepsilon)$; there seems at this stage to be nothing to prevent $\rho(x_\varepsilon)$ from going to 0 on an $\varepsilon'$-dense set in $\Omega_\varepsilon$, with $\varepsilon' \to 0$ as $\varepsilon \to 0$. Nevertheless, it is simplest, and most natural, to choose the points $x_\varepsilon$ to realize the infimum of $\rho$ on $(\Omega_\varepsilon, g_\varepsilon)$, i.e. choose points $y_\varepsilon$ s.t.

$$\rho(y_\varepsilon) \le \rho(q_\varepsilon), \tag{4.14}$$

for all $q_\varepsilon \in \Omega_\varepsilon$. Such points have the highest concentration of curvature in $L^2$ on $(\Omega_\varepsilon, g_\varepsilon)$. If $\Omega_\varepsilon$ is compact, such points $y_\varepsilon$ clearly exist. More generally, even if $\Omega_\varepsilon$ is non-compact, one may choose points almost satisfying (4.14) locally. This is expressed by the following result.

**Lemma 4.1.** *Let $\{x_\varepsilon\}$ be any sequence of points in $\{(\Omega_\varepsilon, g_\varepsilon)\}$, $\varepsilon = \varepsilon_i$, such that*

$$\rho(x_\varepsilon) \to 0 \ \text{ as } \ \varepsilon \to 0. \tag{4.15}$$

*Then for any fixed $\delta > 0$, there are points $y_\varepsilon \in B_{x_\varepsilon}(\delta)$ such that*

$$\rho(y_\varepsilon) \to 0 \ \text{ as } \ \varepsilon \to 0, \tag{4.16}$$

*and such that for any $K < \infty$ and for all $q_\varepsilon \in B_{y_\varepsilon}(K\rho(y_\varepsilon))$,*

$$\rho(y_\varepsilon) \le 2\rho(q_\varepsilon), \tag{4.17}$$

*provided $\varepsilon$ is sufficiently small, depending only on the choice of $K$.*

*Proof.* Given $B_\varepsilon \equiv B_{x_\varepsilon}(\delta)$, let $t(p) = \text{dist}(p, \partial B_\varepsilon)$ and consider the scale-invariant function

$$\psi = \psi_\varepsilon = \rho/t$$

on $B_\varepsilon$. This has value $+\infty$ on $\partial B_\varepsilon$ and at the center point $x_\varepsilon, \psi(x_\varepsilon) \to 0$. Let $y_\varepsilon$ be any point realizing the minimal value of $\psi_\varepsilon$. Thus $y_\varepsilon \in B_\varepsilon$ and $\psi_\varepsilon(y_\varepsilon) \to 0$ as $\varepsilon \to 0$. Since $t \le \delta$ on $B_\varepsilon$, it follows in particular that

$$\rho(y_\varepsilon) \to 0,$$

and further

$$\rho(y_\varepsilon) << t(y_\varepsilon), \ \text{ as } \ \varepsilon \to 0. \tag{4.18}$$

By the minimizing property of $y_\varepsilon$, for any $q_\varepsilon \in B_\varepsilon$,

$$\rho(y_\varepsilon) \le \rho(q_\varepsilon) \cdot \frac{t(y_\varepsilon)}{t(q_\varepsilon)}.$$

If $q_\varepsilon \in B_{y_\varepsilon}(K\rho(y_\varepsilon))$, and $\varepsilon$ is sufficiently small, depending only on $K$, then (4.18) and (2.2) imply that the ratio $t(y_\varepsilon)/t(q_\varepsilon)$ is close to 1. Thus (4.17) follows. ∎

For the remainder of this section, we assume $y_\varepsilon \in \Omega_\varepsilon$ is chosen to satisfy (4.16)-(4.17); (some remarks on the general case where this is not assumed are made in Remark 4.5). Note that the condition (4.17) implies that the metrics $g_\varepsilon$ are strongly $(\rho, c)$ buffered, for a fixed $c > 0$, c.f. §2. As in (4.1) consider the blow-up sequence

$$g'_\varepsilon = \rho(y_\varepsilon)^{-2} \cdot g_\varepsilon. \tag{4.19}$$

The (scale-invariant) estimate (4.17) implies that for any $q_\varepsilon$ within uniformly bounded $g'_\varepsilon$-distance to $y_\varepsilon$, $\rho(q_\varepsilon, g'_\varepsilon) \geq \frac{1}{2}$, provided $\varepsilon$ is sufficiently small.

If $\{g'_\varepsilon\}$ is not collapsing at $y_\varepsilon$, then the estimates (2.7)-(2.9) together with (4.17), imply that $\{g'_\varepsilon\}$ is not collapsing within any given bounded $g'_\varepsilon$-distance to $y_\varepsilon$. In other words, there is no cusp formation within bounded distance to $y_\varepsilon$. It follows from applications of Theorem 2.3 to $g'_\varepsilon$-balls $B_{q_\varepsilon}(\frac{1}{2})$, within bounded $g'_\varepsilon$-distance to $y_\varepsilon$, that there is a sequence $R_i \to \infty$ and a sequence $\varepsilon_i \to 0$ such that the pointed manifolds $(B_{y_\varepsilon}(R_i), g'_\varepsilon, y_\varepsilon), \varepsilon = \varepsilon_i$, sub-converge, (modulo diffeomorphisms), to a $complete$ limit $(N, g', y)$ with $L^{2,2}$ Riemannian metric $g'$. The convergence is in the weak $L^{2,2}$ topology, and uniform on compact subsets of $N$.

On the other hand, if $\{g'_\varepsilon\}$ is collapsing at $\{y_\varepsilon\}$, then the same estimates (2.7)-(2.9) together with (4.17) imply that $\{g'_\varepsilon\}$ is collapsing everywhere within uniformly bounded $g'_\varepsilon$-distance to $y_\varepsilon$. Further the collapse is along a sequence of injective F-structures so that one may unwrap the collapse by passing to the universal cover $\widetilde{B}_{y_\varepsilon}(R)$ of $(B_{y_\varepsilon}(R), g'_\varepsilon, y_\varepsilon)$ for any given $R < \infty$, c.f. Remark 2.1 and (2.7)-(2.9). The sequence $\{\widetilde{B}_{y_\varepsilon}(R), g'_\varepsilon, y_\varepsilon)\}$, ($y_\varepsilon$ is now some lift to the universal cover), is not collapsing, and choosing a sequence $R_i \to \infty$ as above, it follows as before that the pointed manifolds $(\widetilde{B}_{y_\varepsilon}(R_i), g'_\varepsilon, y_\varepsilon), \varepsilon = \varepsilon_i$, sub-converge in the weak $L^{2,2}$ topology to a $complete$ limit $(\widetilde{N}, g', y)$. In addition, the limit $(\widetilde{N}, g')$ here has (at least) a free isometric $\mathbb{R}$-action. In this case, we will always work on the $\mathbb{Z}$-quotient $N$ of $\widetilde{N}$ so that $(N, g')$ has a free isometric $S^1$-action.

Observe that this dichotomy non-collapse/collapse of the blow-up sequence (4.19) depends only on the behavior at the base point $y_\varepsilon$, satisfying (4.17).

In order to obtain any relation between the geometry of a limit $(N, g', y)$ and that of the approximating sequence $(\Omega_\varepsilon, g'_\varepsilon, y_\varepsilon)$, one needs to prove at least strong $L^{2,2}$ convergence to the limit, compare with §1(I). More generally, the following result improves the initial uniform $L^{2,2}$ regularity of the metric $g_\varepsilon$ and functions $\tau = \tau_\varepsilon$ on balls the size of the $L^2$ curvature radius as $\varepsilon$ varies. We recall from §4.1 that in case $\sigma(M) = 0$, $\varepsilon$, $\tau$, $\omega$ and $c$ are always considered to be divided by $T_\varepsilon = \inf|\tau_\varepsilon|$.

**Theorem 4.2.** *There is a constant $c_1 > 0$, independent of $\varepsilon$ and $x_\varepsilon$, such that for any $x_\varepsilon \in (\Omega_\varepsilon, g_\varepsilon)$,*

$$\rho^{1,2}(x_\varepsilon) \geq c_1 \rho(x_\varepsilon), \tag{4.20}$$

*for $\varepsilon \leq 1$, where $\rho^{1,2}$ is the $L^{1,2}$ curvature radius, c.f. §2. Further, there is a constant $c_2 > 0$, independent of $\varepsilon$ and $x_\varepsilon$, such that for the potential function*

$$\omega = \tau + \frac{\varepsilon}{12} s, \tag{4.21}$$

*one has the estimate,*

$$\rho(x_\varepsilon)^2 ||D^2 \omega||_{\bar{L}^2(B_{x_\varepsilon}(c_1 \rho(x_\varepsilon)))} \leq c_2 \cdot max\big(||\omega||_{\bar{L}^2(B_{x_\varepsilon}(\rho(x_\varepsilon)))}, 1\big), \tag{4.22}$$

*where $\bar{L}^2$ is the $L^2$ average w.r.t. Lebesgue measure.*

*Proof.* This result is a parametrized version of Theorem 3.3. Namely Theorem 3.3 implies that the estimates (4.20) and (4.22) hold for any fixed $\varepsilon > 0$, for $some$ constants $c_1$, $c_2$, which apriori depend on $\varepsilon$. In fact, Theorem 3.3 gives local $L^{3,p}$ regularity of the metric $g_\varepsilon$ and potential $\omega$. Considering the family of metrics $g_\varepsilon$ as $\varepsilon$ varies, the content of Theorem 4.3 is then that the constants $c_1$ and $c_2$ are independent of $\varepsilon$ and the size of $\rho(x_\varepsilon)$, for $\varepsilon$ small.

First, for convenience, we scale the metrics $g_\varepsilon$ at $x_\varepsilon$ as in (4.1) so that $\rho(x_\varepsilon) = 1$, giving the equations (4.2)-(4.3), and pass to the universal cover in case $g_\varepsilon$ is sufficiently collapsed at $x_\varepsilon$, so that $r_h(x_\varepsilon) \sim 1$ where $r_h$ is the $L^{2,2}$ harmonic radius, c.f. Theorem 2.3. Observe that the estimates (4.20)-(4.22) are scale invariant and invariant under coverings. By a further bounded rescaling, we normalize the metric so that $r_h(x_\varepsilon) = 1$, where $r_h$ is the harmonic radius, (c.f. §2). For notational simplicity, we drop the prime from the notation, so that $g'_\varepsilon$ is denoted by $g_\varepsilon$ or just $g$. In particular, observe that the $L^2$ norm of the full curvature $r_{g_\varepsilon}$ is uniformly bounded in $B_{x_\varepsilon}(1)$.

Let $\alpha = \alpha(\varepsilon)$ be the resulting coefficient in (4.2) and set $B = B_{x_\varepsilon}(1)$. Note that $\tau$ is scale-invariant and $\tau \leq 0$. Further, $\tau$ is also uniformly bounded below for $\varepsilon \leq 1$ by Theorem 3.11 and the convention on $T$ above.

Theorem 4.2 is proved in [An4, Thm. 3.6] for $\mathcal{Z}^2_s$ solutions (4.9), i.e. for the limit equation of (4.2)-(4.3) when $\omega = \tau$, (and so $s = 0$). In fact the same bounds, independent of $\alpha$, are obtained in this case for the $L^{k,2}$ curvature radius $\rho^{k,2}$ and $L^{k,2}$ norm of $\omega$, for any $k$. Because the proof of Theorem 4.2 is very similar, and because a full proof with complete details would be very long and repetitious, we give the proof with reasonably complete details and refer to [An4] in case further details are desired.

Before beginning, we recall a key point. By Lemma 3.4, one has an $L^{k,p}$ bound on each summand of $\omega$ in terms of an $L^{k,p}$ bound on $\omega$, provided $k \leq 1$. We also assume familiarity with the proof of Theorem 3.3.

The proof must be separated into several cases, depending on the (relative) behavior of the coefficient $\alpha$ and potential $\omega$ as $\varepsilon \to 0$ on any given sequence $\varepsilon = \varepsilon_i \to 0$, (as in [An4, Thm. 3.6]).

**(i).** Suppose first there is some constant $\alpha_o > 0$ such that $\alpha(\varepsilon) \geq \alpha_o$, as $\varepsilon \to 0$. The main point in this case is that $\alpha \nabla \mathcal{Z}^2$ is the *uniformly* dominant term in the elliptic system (4.2). The uniform lower bound on $\alpha$ here is the same as a uniform lower bound on $\varepsilon$ away from 0 in the proof of Theorem 3.3 and it is clear that the estimates obtained there are independent of $\varepsilon$ when $\varepsilon$ is bounded below. (Note that these estimates are improved when $\alpha$ is large, (or $\varepsilon$ large in Theorem 3.3), as one sees by dividing everywhere by $\alpha$ ). Thus the estimates (4.20)-(4.22) follow exactly as in the proof of Theorem 3.3 where $\varepsilon$ was fixed. The constants $c_1$ and $c_2$ in (4.20)-(4.22) of course then depend on $\alpha_o$.

We may thus suppose that $\alpha$ is bounded, say $\alpha \leq 1$, as $\varepsilon \to 0$. We claim that it then follows as in the proof of Theorem 3.3 that the $L^{2,2}$ norm of $\omega = \omega_\varepsilon$ is uniformly bounded in $B(1 - \delta)$, for any fixed $\delta > 0$, and hence $\omega$ is uniformly controlled in $C^{1/2}$ in $B(1 - \delta)$, i.e.

$$\|\omega\|_{C^{1/2}} \leq K, \tag{4.23}$$

where $K$ depends only on $\delta$. The reasoning here is as follows, referring to (4.2)-(4.3) and the proof of Theorem 3.3. First $\tau$ is bounded in $L^\infty$ by the remarks above, while the bound $\alpha \leq 1$ implies that the term $\varepsilon s$ in the potential $\omega$ in (4.21) is bounded in $L^2$, since in this scale, $\varepsilon s$ becomes $\alpha s$, and $s$ is bounded in $L^2$. Hence, $\omega$ is controlled in $L^2(B)$. The left side of (4.3) is uniformly bounded in $L^1(B)$ and hence $\omega$ is uniformly controlled in $L^{t,2}(B')$, $t < \frac{1}{2}$, where $B(1 - \delta) \subset B' \subset B$. By Lemma 3.4, $\tau$ is then also controlled in $L^{t,2}(B')$ and so $L^*\tau$ is controlled in $L^{t-2,2}(B')$. This gives a uniform bound on $\alpha z$ in $L^{t,2} \subset L^{3-\mu}$, $\mu > 0$, which then gives a uniform bound on $\alpha|z|^2$ in $L^{(6/5)-\mu}$ by the Hölder inequality, on $B(1 - \delta) \subset B'' \subset B'$. In turn and as before, this gives uniform control on $\omega$ in $L^{2,(6/5)-\mu} \subset L^{1,2-\mu}$, and so on. Continuing in this way gives uniform $L^{2,2}$ control on $\omega$ in $B(1 - \delta)$, (c.f. also [An4, Cor. 3.5]).

Observe that these arguments already prove (4.22) in all cases.

**(ii)(a).** Suppose next there is some constant $\tau_o < 0$ such that

$$\tau(x_\varepsilon) \leq \tau_o, \tag{4.24}$$

i.e. $\tau(x_\varepsilon)$ has a definite bound away from 0. Then (4.23) implies that $\tau$ is bounded away from 0 in a ball $B' \subset B$, whose radius is of a definite size, depending only on the value $\tau_o$. Thus in $B'$, the family of equations (4.2)-(4.3) are essentially the $\mathcal{Z}_s^2$ equations (4.9); the difference is only in the 0-order terms $c$ and $s\tau$, which are controlled when $\varepsilon$ is small. In particular, as in Theorem 3.3, the metric $g_\varepsilon$ and potential $\omega$ are $C^\infty$ in $B'$.

The idea now is that when $\tau$ is bounded away from 0, the standard elliptic regularity estimates one would obtain from each term $\alpha \nabla \mathcal{Z}^2$ and $L^*\tau$ in (4.2) alone, act together, and not in conflict, with each other. Thus, the leading order terms in (4.2) are

$$\alpha D^* Dz - \tau \cdot r + D^2 \tau.$$

The Hessian $D^2\tau$ is controlled by $\Delta\tau$, which in turn is controlled by the trace equation (4.3). Now $\alpha D^* Dz$ is a positive operator on $z$ and when $\tau$ is strictly negative in $B'$, $-\tau \cdot r$ is of course also strictly positive in $z$. This implies that elliptic estimates obtained from each of these two terms separately are *reinforced*, in the full equation (4.2). This is explained in full detail in [An4, Lemmas 3.10-3.12], and so we only outline the main idea here. First, pair the equation (4.2) with $z$; the leading order term then becomes

$$-\alpha \Delta |z|^2 + \alpha |Dz|^2 - \tau |z|^2 + < D^2\tau, z >.$$

If one pairs this with a smooth cutoff function $\eta$ supported in $B'$ and integrates by parts, one obtains a uniform bound on $\alpha \int |Dz|^2$, and hence by Sobolev embedding a uniform bound on $\alpha (\int |z|^6)^{1/3}$, where the integrals are over $B'' \subset B'$, c.f. [An4, Lemma 3.10]; here we are using the fact that the $L^2$ norm of $z$ is uniformly bounded on $B$. Similarly, pairing with $\eta |z|$ and integrating by parts gives uniform bounds on $\alpha (\int |z|^9)^{1/3}$ and on $\int |z|^3$ on a slightly smaller ball, c.f. [An4, Lemma 3.11]. This is due to the fact that both terms $\alpha$ and $-\tau$ are *positive* in $B'$. The bound on $\int |z|^3$ depends on the upper bound $\tau_o$ on $\tau$ from (4.24). Note however that since $\alpha$ may be arbitrarily small, these arguments do not give uniform bounds on $\int |Dz|^2$. To achieve this, one takes the covariant derivative of (4.2), pairs (4.2) with $\eta Dz$, and computes in a similar way as above, c.f. [An4, Lemma 3.12]. Here we note that, as remarked before, the metric $g_\varepsilon$ and potential $\omega$ are $C^\infty$ smooth in $B'$, and so satisfy all covariant derivatives of the equation (4.2). An $L^{1,2}$ bound on $s$ then follows from the Bianchi identity and hence one obtains a lower bound for the $L^{1,2}$ curvature radius.

In this way one obtains the estimates (4.20)-(4.22) with constants $c_1$ and $c_2$ depending however (only) on $\tau_o$.

**(ii)(b).** Consider next the opposite case to (4.24) where $\tau \equiv 0$ in some ball $B' \subset B$, so that $s(x) \geq 0$, $\forall x \in B'$. In this case, the family of equations (4.2)-(4.3) are essentially the $\mathcal{Z}^2$ equations (4.6) in $B'$, after dividing by $\alpha$. The *only* difference is in the (low order) constant term $c_\varepsilon$. Further, $c_\varepsilon / \alpha$ is uniformly bounded above, as one sees by pairing the trace equation (4.2) with a smooth cutoff function $\eta$ in $B'$ and applying the self-adjointness of $\Delta$.

Thus, one has a family of solutions $\{g'_\varepsilon\}$ to an essentially fixed equation. Standard elliptic estimates on the $\mathcal{Z}^2$ equations, c.f. [GT, Chs.8.3,9.5], then imply (4.20)-(4.22), as well as their higher order $L^{k,2}$ analogues, for any $k$.

**(iii).** Hence, since this is the only case left, we may assume that $s$ changes sign say in $B(\frac{1}{2})$, (so that the junction $\Sigma$ is non-empty), that $\alpha$ is small, and that $\tau(x_\varepsilon) < 0$ is small, i.e.

$$\tau(x_\varepsilon) \geq \tau_o. \tag{4.25}$$

Hence $\tau$ is small in some neigborhood of $x_\varepsilon$ by (4.23). Observe that this implies that *all* quantities in (4.2) are (weakly) small in that neighborhood.

In this case, we need to renormalize the equations (4.2)-(4.3) again, so that $\tau$ is of roughly unit (negative) size in a small ball of definite size. We carry out this process below, and then show that this case reduces to one of the previous cases.

For the remainder of the proof, we change notation and assume that $B = B_{x_\varepsilon}(\frac{1}{2})$. We apply Proposition 2.4 to the trace equation (4.3), with $u = -\omega = -\tau - \frac{\alpha}{12}s$. From (2.13), this gives

$$\int_B u \le c_3 \int_B (\alpha|r|^2 + \frac{1}{4}s\tau - c) + \int_{B_\delta} u + 2\pi u_{av}(\delta),$$

where $\delta = \frac{1}{10}$, $B_\delta$ is the $\delta$-ball about the center $x = x_\varepsilon$ and $u_{av}(\delta)$ is the average value of $u$ on $S_\delta = \partial B_\delta$. Thus,

$$\int_B -\tau \le c_3 \int_B (\alpha|r|^2 + \frac{1}{4}s\tau - c) + \int_B \frac{\alpha}{12}s + \int_{B_\delta} u + 2\pi u_{av}(\delta), \tag{4.26}$$

Note that the constant term $c = c_\varepsilon$ here is positive. Before proceeding further, we claim that under the assumption (4.25) that $\tau(x_\varepsilon)$ is small, (the required smallness will be defined below),

$$s\tau(q) < c_3|\tau(q)|, \tag{4.27}$$

for all $q \in B$, for all $\varepsilon$ sufficiently small. To see this, if $\sigma(M) = 0$, then by Theorem 3.12, inf $s \to 0$ pointwise, so (4.27) follows, since $s\tau = 0$ when $s > 0$. Thus assume $\sigma(M) < 0$. If $\rho(x_\varepsilon, g_\varepsilon)$ is sufficiently small (as $\varepsilon \to 0$), then $s = s_{g'_\varepsilon} = \rho(x_\varepsilon)^2 \cdot s_{g_\varepsilon}$ and $s_{g_\varepsilon}$ is uniformly bounded below by Theorem 3.11, so that (4.27) again follows. Hence (4.27) holds if $\rho(x_\varepsilon, g_\varepsilon) \le \rho_o$, for some (absolute) $\rho_o > 0$. Finally, suppose $\rho(x_\varepsilon, g_\varepsilon) \ge \rho_o$. By (4.23), the functions $\tau_\varepsilon$ are uniformly bounded in $C^{1/2}$, while $\omega_\varepsilon$ is uniformly bounded in $L^{2,2}$ on $B = B_{x_\varepsilon}(\frac{1}{2})$. Taking a limit as $\varepsilon \to 0$, it follows that the limit function $\tau$ satisfies $2\Delta\tau + \frac{1}{4}s\tau = c$. If $x = \lim x_\varepsilon$ and $\tau(x) = 0 = \sup \tau$, this contradicts the maximum principle, (c.f. [GT, Thm. 3.5]), since $c > 0$ in the limit. Hence, it follows in this case (by the $C^\alpha$ control) that $\tau(x_\varepsilon)$ could not have been too small, depending only on $\rho_o$. This determines a bound for $\tau_o$ away from 0 in (4.25), i.e. the cutoff between the cases (ii)(a) and (iii) above. Thus, if $\tau_o$ in (4.25) is chosen sufficiently small, then (4.27) holds.

The point here is that the estimate (4.27) implies that the positive term $\frac{1}{4}s\tau$ on the right in (4.26) may (and now will) be absorbed into the left side of (4.26).

Next, since $\tau \le 0$, we may apply (2.14) to any ball $\bar B \subset B$, with $u = -\tau - \frac{\alpha}{12}s$ as before to obtain

$$d \cdot \| -\tau \|_{L^2(\bar B)} \le c_4 \left( \int_B \alpha|r|^2 + \int_B \frac{\alpha}{12}s^+ + \frac{1}{2}(-\tau)_{av}(\delta) \right). \tag{4.28}$$

Here $d$ depends only on $dist(\partial \bar B, \partial B)$ and $c_4$ depends only on $c_3$.

The quantity $(-\tau)_{av}(\delta)$ is small, since otherwise by (4.23) one is in case (ii)(a). Observe then that the $L^2$ norm of $-\tau$ on $\bar B$ is small, since all other terms on the right in (4.28) are small, (since $\alpha$ is small). In this situation renormalize $\tau$ and $\alpha$ by setting

$$\bar\tau(x) = \frac{\tau(x)}{|\tau_{av}(\delta)|}, \bar\alpha = \frac{\alpha}{|\tau_{av}(\delta)|}. \tag{4.29}$$

Similarly, the equations (4.2)-(4.3) are renormalized by dividing everywhere by $|\tau_{av}(\delta)|$, as is (4.28). (This is the same renormalization as in [An4,(3.35)-(3.36)]).

Given this renormalization, we now repeat the analysis above in the previous cases. First, if $\bar\alpha$ is bounded away from 0, (4.20)-(4.22) follow as in the initial situation above, when $\alpha$ was bounded away from 0, (Case (i)). Thus, we may assume $\bar\alpha$ is also small and hence by (4.28), $\bar\tau$ is bounded in $L^2(\bar B)$. As above, this data implies that $\bar\tau$ is bounded in $L^{2,2} \subset C^{1/2}$, so that (4.23) holds with $\bar\tau$ in place of $\omega$ inside $\bar B$. Now however $\bar\tau$ is bounded away from 0 inside $\bar B$, and so (4.20)-(4.22) follow as in the arguments in Case (ii)(a). ∎

**Remark 4.3. (i).** The same proof as above shows that the bound (4.20) can be improved to a bound on the $L^{1,p}$ curvature radius for any $p < \infty$, i.e.

$$\rho^{1,p}(x_\varepsilon) \ge c_1(p)\rho(x_\varepsilon). \tag{4.30}$$

However, this estimate cannot be improved beyond this range in neighborhoods of the junction $\Sigma = \partial\{s = 0\}$, for the same reasons as discussed in §3. Namely, one does not expect $g_\varepsilon$ to be smoother than $L^{3,p}$ in a neighborhood of $\Sigma$. Similarly, the bound (4.22) can be improved, (with the same methods), to a $L^{3,p}$ estimate on $\omega$, but not further. As in Lemma 3.4, one then also has the estimate

$$||D\tau||_{L^p(B_{x_\varepsilon}(c_1))} \leq c_2(p) \cdot max\left(||\omega||_{L^2(B_{x_\varepsilon}(1))}, 1\right), \tag{4.31}$$

in balls normalized so that $\rho(x_\varepsilon) \sim r_h(x_\varepsilon) \sim 1$.

On the other hand, away from the junction $\Sigma$, i.e. within $U^+$ and $U^-$ from (3.24), as explained in the proof above, the estimates (4.20)-(4.22) can be improved to similar estimates on $\rho^{k,2}$ and the $L^{k,2}$ norm of $\omega$ as above, for any $k < \infty$; the estimates then depend on the distance to $\Sigma$.

**(ii).** The validity of (4.20) and (4.22) depends strongly on the fact that $\tau \leq 0$. Thus, such estimates will not hold in this generality for the functional $I_\varepsilon$ from (3.6) with $\mathcal{S}^2$ in place of $\mathcal{S}_-^2$. The reason for this, as explained in the proof of Theorem 4.2, is that the two terms $\varepsilon\nabla\mathcal{Z}^2$ and $L^*\tau$ in (4.2) act in unison with each other where $\tau < 0$, (regarding regularity estimates), but may interfere with each other where $\tau > 0$. This sensitivity to the sign of $\tau$ is in contrast to the situation in [AnI], where the sign of the potential function $u$ played no role.

Now return to the blow-up sequence $(\Omega_\varepsilon, g'_\varepsilon, y_\varepsilon)$ from (4.19). Theorem 4.2 implies that the $L^{1,2}$ curvature radius of $g'_\varepsilon$ is uniformly bounded below, within uniformly bounded distance to the base points $y_\varepsilon$. Since the embedding $L^{1,2} \subset L^2$ is compact, this implies that $\{(B_{y_\varepsilon}(R_i), g'_\varepsilon, y_\varepsilon)\}, \varepsilon = \varepsilon_i$, (or sufficiently large covers of it, as discussed prior to Theorem 4.2), converges *strongly* in $L^{2,2}$ to its limit $(N, g', y)$. Under these conditions, the $L^{2,2}$ curvature radius is continuous, (c.f. [AnI, (1.28)] and references therein), i.e.

$$\rho'(y) = lim_{\varepsilon\to 0}\rho'_\varepsilon(y_\varepsilon) = 1, \tag{4.32}$$

so that the limit $(N, g')$ is *not flat*. By (4.17) and (4.30), the limit $(N, g')$ is complete and has uniformly bounded curvature, and $g'$ is uniformly bounded in $L^{3,p}$. It is also clear from Theorems 3.11 and 3.12 that the scalar curvature $s$ of the limit $g'$ is non-negative. By (4.31) the functions $\tau_\varepsilon$ are locally uniformly bounded in $L^{1,p}$ and so (sub)-converge to a limit $L^{1,p}$ function $\tau$ on $N$.

Following the discussion in §4.1, we now consider the form of the Euler-Lagrange equations on the limit $(N, g', y)$. Thus, return to the Euler-Lagrange equations (4.2)-(4.3) in the scale $g'_\varepsilon$, i.e.

$$\frac{\varepsilon}{\rho^2}\nabla\mathcal{Z}^2 + L^*(\tau) + (\frac{1}{4}s\tau + c_\varepsilon) \cdot g = 0, \tag{4.33}$$

$$2\Delta(\tau + \frac{\varepsilon s}{12\rho^2}) + \frac{1}{4}s\tau = -\frac{1}{2}\frac{\varepsilon}{\rho^2}|z|^2 + 3c_\varepsilon, \tag{4.34}$$

where as before we have dropped the prime and $\varepsilon$, and $\rho = \rho(y_\varepsilon, g_\varepsilon)$. These equations are clearly invariant under coverings, and so in the case of collapse, also hold on sufficiently large or the universal cover. The remarks above on replacing $\varepsilon$ by $\varepsilon/T, \tau$ by $\tau/T$ and $c$ by $c/T$ when $\sigma(M) = 0$ also remain in effect.

**Case (I).** $\varepsilon/\rho^2 \to \infty$.

As before, divide the equations (4.33)-(4.34) by $\varepsilon/\rho^2$. Since $\tau_\varepsilon$ is uniformly bounded, $(\rho^2/\varepsilon)\tau_\varepsilon \to 0$ uniformly on compact subsets of $(N, g')$. It follows that the limit $(N, g', y)$ is a complete non-flat solution to the $\mathcal{Z}^2$ equations (4.6). Apriori, the limit metric $g'$ is only $L^{3,p}$ and is a weak $L^{3,p}$ solution of (4.6). However since $s \geq 0$ on $(N, g')$, the maximum principle applied to the trace equation in (4.6) shows that $s > 0$ everywhere on $(N, g')$. By Remark 4.3(i), (or Case (ii)(b) in the proof of Theorem 4.2), in this case $g'$ is $C^\infty$ and $g'$ is a smooth solution of the equations.

**Case (II).** $\varepsilon/\rho^2 \to \alpha > 0$.

In this case, the limit equations take the form (4.7), i.e. the constrained $\mathcal{Z}_c^2$ equations. It is useful to distinguish three further subcases according to the behavior of the limit function $\tau$ on $(N, g')$.

**(a).** $\tau \equiv 0$. In this case, the limit $(N, g')$ is as in Case (I) above, i.e. a complete non-flat smooth solution of the $\mathcal{Z}^2$ equations (4.6).

**(b).** $\tau < 0$ everywhere. In this case, the limit $(N, g')$ is a complete non-flat scalar-flat solution to the (pure) $\mathcal{Z}_s^2$ equations (4.9). Thus there is no junction region $\Sigma = \partial(s = 0)$ in this case. Again by Remark 4.3(i), (or [An4, Thm. 3.6]), the metric $g'$ is smooth, and is a smooth solution to the equations.

**(c).** $\tau < 0$ somewhere and $s > 0$ somewhere. In this case, the junction $\Sigma$ is non-empty. The proof of Theorem 3.3 implies that the limit $(N, g', \omega)$ is a complete $L^{3,p}$ weak solution to the equations (4.7), which may or may not be smooth in a neighborhood of the junction $\Sigma$, but is smooth away from $\Sigma$. To the 'right' of $\Sigma$, i.e. in the region $N^+$ given by the interior of $\{s \geq 0\}$, the metric is a $\mathcal{Z}^2$ solution. To the 'left' of $\Sigma$, i.e. in the region $N^- = \{\tau < 0\}$, the metric is a $\mathcal{Z}_s^2$ solution, as in (4.9). The functions $\tau$ and $s$ are locally Lipschitz, satisfying (4.8), while $\omega$ is a $C^{2,\alpha}$ function on $N$.

**Case (III).** $\varepsilon/\rho^2 \to 0$.

Again we distinguish several possibilities according to the behavior of $\tau$ on the limit $(N, g')$.

**(a).** $\tau < 0$ everywhere. In this case, the limit $(N, g', \tau)$ is a complete $C^\infty$ solution of the static vacuum equations (4.10).

**(b).** $\tau < 0$ somewhere. In the region $N^- = \{\tau < 0\}$, $(g', \tau)$ is a solution of the static vacuum equations with non-empty horizon $\Sigma = \partial\{\tau < 0\}$; (otherwise, if $\Sigma = \emptyset$, one is in Case (a)). Since $\tau \to 0$ somewhere, the maximum principle applied to the limit of the trace equation (4.34) implies that there is a non-empty region $N^+ \subset N$, given as before as the interior of $\{s \geq 0\}$, with $\partial N^+ = \Sigma$. In this region, the metric is a $\mathcal{Z}^2$ solution, as one sees by renormalizing the equations (4.33)-(4.34) by dividing by $\varepsilon/\rho^2$, using the fact that $\tau \equiv 0$ on $N^+$. Thus, there is a junction between these solutions at $\Sigma$. Again the metric $g'$ is $L^{3,p}$ smooth across $\Sigma$.

**(c).** $\tau \equiv 0$. In this case, one would seem to obtain only a super-trivial solution of the static vacuum equations; of course the solution is such a super-trivial solution, i.e. both sides of (4.33)-(4.34) are 0 in the limit. However, there must be some non-trivial equation describing the structure of the limit, since the minimizing properties of the metrics $g'_\varepsilon$ cannot disappear at $g'$. To determine the limit equations in this case, return to (4.33)-(4.34) and, as in Case (I), divide by $\varepsilon/\rho^2$, to obtain the equations

$$\nabla \mathcal{Z}^2 + L^*(\frac{\rho^2}{\varepsilon}\tau) + \frac{\rho^2}{\varepsilon}(\frac{1}{4}s\tau + c_\varepsilon) \cdot g = 0, \tag{4.35}$$

$$2\Delta(\frac{\rho^2}{\varepsilon}\tau + \frac{1}{12}s) + \frac{\rho^2}{4\varepsilon}s\tau = -\frac{1}{2}|z|^2 + \frac{3\rho^2}{\varepsilon}c_\varepsilon, \tag{4.36}$$

Now of course $\rho^2/\varepsilon \to \infty$ while $|z|^2$ and $s$ are bounded. Recall that in this scale, the constant term $c_\varepsilon = \rho^2 c_{g_\varepsilon}$, where $c_{g_\varepsilon} = \frac{1}{6}(\varepsilon \mathcal{Z}^2(g_\varepsilon) + \frac{1}{2}\sigma(g_\varepsilon))$ is the constant term on the original $(\Omega_\varepsilon, g_\varepsilon)$ sequence. Thus,

$$\frac{\rho^2}{\varepsilon}c_\varepsilon = \frac{\rho^4}{\varepsilon}c_{g_\varepsilon}, \tag{4.37}$$

where $c_{g_\varepsilon}$ is bounded. (In fact $c_{g_\varepsilon} \to 0$ as $\varepsilon \to 0$ in case $\sigma(M) = 0$). Apriori, one then has again three possibilities, namely $\rho^4 c_{g_\varepsilon}/\varepsilon$ goes to $\infty$, is bounded, or goes to 0, (in a subsequence).

For the moment, we suppose that

$$c_{g_\varepsilon}\rho^4/\varepsilon \to 0, \tag{4.38}$$

so that the constant term in (4.35)-(4.36) goes to 0 as $\varepsilon \to 0$. Similarly, since $s = s_{g'_\varepsilon}$ in (4.35)-(4.36) is of the form $s = \rho^2 \cdot s_{g_\varepsilon}$, the term $(\rho^2/\varepsilon)s\tau$ in (4.35)-(4.36) also goes to 0 under the estimate (4.38).

We then have three further subcases of this case according to the behavior of the (renormalized) potential function

$$\bar{\tau} = \bar{\tau}_\varepsilon = (\rho^2/\varepsilon)\tau_\varepsilon. \tag{4.39}$$

Note first that $\bar{\tau} \leq 0$, and $s$ is uniformly bounded above, so that $\bar{\omega} = \bar{\tau} + \frac{1}{12}s$ is uniformly bounded above as $\varepsilon \to 0$. Hence, elliptic regularity applied to the trace equation (4.36), c.f. [GT, Thms. 8.17, 8.18], implies that the oscillation of $\bar{\omega}$ is uniformly controlled locally, i.e. $sup|\bar{\omega}| \leq \kappa \cdot (inf|\bar{\omega}| + 1)$, for some constant $\kappa$, where the sup and inf are taken over unit balls in $(\Omega_\varepsilon, g'_\varepsilon)$ of $g'_\varepsilon$-bounded distance to $y_\varepsilon$. This estimate thus also holds for $\bar{\tau}$, by Lemma 3.4. The three subcases are:

$(\mathbf{c_{(I)}})$. $\bar{\tau} = (\rho^2/\varepsilon)\tau \to 0$ everywhere. Taking the limit of the equations (4.35)-(4.36), it follows that the limit is a $\mathcal{Z}^2$ solution, as in Case (I).

$(\mathbf{c_{(II)}})$. $\bar{\tau} = (\rho^2/\varepsilon)\tau$ bounded, away from 0 and $-\infty$, on some sequence $x_\varepsilon \to x \in (N, g')$. From the estimates above on $\bar{\tau}$ and $\bar{\omega}$, it follows that $\bar{\tau}$ is locally bounded everywhere on the limit $(N, g')$. In this case, the limit is either a (pure) $\mathcal{Z}_s^2$ solution, or a $\mathcal{Z}_c^2$ solution with non-empty junction, as in Cases (II)(b) or (c).

$(\mathbf{c_{(III)}})$. $\bar{\tau} = (\rho^2/\varepsilon)\tau \to -\infty$, at the base point $y_\varepsilon$, and hence everywhere within $g'_\varepsilon$-bounded distance to $y_\varepsilon$. In this case, renormalize the equations (4.35)-(4.36) once more, (for the last time), by dividing by $\rho^2|\tau(y_\varepsilon)|/\varepsilon$, so that the (again renormalized) potential function $\tilde{\tau}$ satisfies $\tilde{\tau}(y_\varepsilon) = -1$. Thus the resulting renormalized coefficients of $s$, $z$, $|z|^2$ all go to 0, as does the constant term. It follows that the limit is a complete solution of the static vacuum equations with potential $\tilde{\tau} < 0$, i.e. there is no horizon or junction $\Sigma$, as in Case (III)(a).

We summarize the results obtained above in the following:

**Proposition 4.4.** *Let $y_\varepsilon$ be base points in $(\Omega_\varepsilon, g_\varepsilon), \varepsilon = \varepsilon_i \to 0$, satisfying (4.16)-(4.17) and, (if necessary), (4.38). Then the blow-up limit $(N, g', y)$ based at $y = \lim y_\varepsilon$ is a complete non-flat solution to either the $\mathcal{Z}^2$ equations, the $\mathcal{Z}_c^2$ or $\mathcal{Z}_s^2$ equations, or the static vacuum Einstein equations.*

$\blacksquare$

It is understood here that if the sequence $\{g'_\varepsilon\}$ is collapsing at $y_\varepsilon$, then the collapse is unwrapped by passing to sufficiently large covers. In such situations, the limit $(N, g')$ has in addition a free isometric $S^1$ action.

**Remark 4.5.** We have chosen the base points $y_\varepsilon$ to satisfy (4.16)-(4.17) since these are geometrically the most natural. However, we point out that Proposition 4.4 holds for arbitrary base points $x_\varepsilon \in (\Omega_\varepsilon, g_\varepsilon)$, (satisfying (4.38) if necessary), except that the limits may be flat and/or not complete, i.e. one may have only locally defined limits.

Flat limits will arise if the base points $x_\varepsilon$ are not $(\rho, c)$ buffered, for a fixed $c > 0$, while incomplete limits will arise if the curvature of the blow-up metrics $g'_\varepsilon = \rho(x_\varepsilon)^{-2} \cdot g_\varepsilon$ based at $x_\varepsilon$ becomes unbounded within bounded $g'_\varepsilon$ distance to $x_\varepsilon$.

If (4.38) does not hold, so that $c_{g_\varepsilon}\rho^4/\varepsilon$ is bounded away from 0, then the curvature is blowing up only very slowly at $\{y_\varepsilon\}$, (in comparison with the other cases above). In this situation, one obtains blow-up limits satisfying equations of the same form as above, except that the constant term $c$ and the $s\tau$ term may no longer vanish. Of course from the discussion above, $c_{g_\varepsilon}\rho^4(x_\varepsilon)/\varepsilon$ can be bounded away from 0 only at points where $\tau(x_\varepsilon) \to 0$, provided $x_\varepsilon$ satisfies (4.16)-(4.17).

While it may be possible to extend the results to follow in §5 and §7 to such equations, we prefer to rule them out by restricting somewhat further the choice of the base point sequence $\{y_\varepsilon\}$. This is done in the next two results. Several possibilities then arise for the specification of $\{y_\varepsilon\}$ which are summarized in Definition 4.8 below.

First we prove that the ratio in (4.38) is necessarily bounded above, under weak conditions.

**Lemma 4.6.** *Let $U_\varepsilon$ be a domain with smooth closure in $(\Omega_\varepsilon, g_\varepsilon)$ for which there exist $q_\varepsilon \in U_\varepsilon$ and $\delta_o > 0$ such that $\tau(q_\varepsilon) \to 0$ as $\varepsilon = \varepsilon_i \to 0$, and $B_{q_\varepsilon}(\delta_o) \subset U_\varepsilon$.*
*Then there is a constant $D_o = D_o(\delta_o) < \infty$ such that*

$$c_{g_\varepsilon} \rho^4(p_\varepsilon)/\varepsilon \leq D_o, \tag{4.40}$$

*for some $p_\varepsilon \in B_{q_\varepsilon}(\delta_o)$.*

*Proof.* Note first that the ratio in (4.40) is invariant under the substitutions $\varepsilon \to \varepsilon/T$ and $c \to c/T$ used in case $\sigma(M) = 0$. Thus, we may and will ignore this renormalization here.

Observe that by Theorem 4.2 and Remark 4.3(i), (4.40) is equivalent to

$$\varepsilon |r|^2(p'_\varepsilon)/c_{g_\varepsilon} \geq d_o, \tag{4.41}$$

for some $d_o > 0$, and some $p'_\varepsilon \in B_{q_\varepsilon}(\delta_o)$. Suppose then that (4.41) does not hold, so that $\varepsilon |r|^2 \to 0$ everywhere in $B_{q_\varepsilon}(\delta_o)$ as $\varepsilon \to 0$. This of course implies that $\varepsilon s^2 \to 0$, and so $\varepsilon s \to 0$ everywhere in $B_{q_\varepsilon}(\delta_o)$, as $\varepsilon \to 0$. Together with the assumption on $q_\varepsilon$, it follows that $m_\varepsilon = sup_{B_{q_\varepsilon}(\delta_o)} \omega \to 0$ as $\varepsilon \to 0$, where $\omega$ is the potential $\omega = \tau + \frac{\varepsilon}{12} s$.

Consider the trace equation (3.12) on $(\overline{\Omega}_\varepsilon, g_\varepsilon)$, i.e.

$$2\Delta\omega + \frac{1}{4}s\tau = -\frac{1}{2}\varepsilon|z|^2 + 3c. \tag{4.42}$$

We have $\omega(q_\varepsilon) \to \lim m_\varepsilon = 0$. It follows that $\Delta\omega \leq 0$ either at $q_\varepsilon$, or at a (nearby) point $p_\varepsilon \in B_\varepsilon(\delta_o)$ still satisfying $\tau(p_\varepsilon) \to 0$. In particular, $\frac{s\tau}{c}(p_\varepsilon) \leq \frac{s\tau}{\sigma}(p_\varepsilon) = \tau^2(p_\varepsilon) \to 0$. However, (4.42) then implies that $(\varepsilon|z|^2/c)(p_\varepsilon) \geq 5$. This contradiction thus implies (4.40) holds. ∎

If $\{p_\varepsilon\} \in B_{q_\varepsilon}(\delta_o)$ is any sequence of points satisfying (4.40), we may apply Lemma 4.1 to obtain points $\{y_\varepsilon\} \in B_{q_\varepsilon}(\delta_o)$ satisfying (4.16)-(4.17). The construction of $\{y_\varepsilon\}$ implies that $\rho(y_\varepsilon) \leq 2\rho(p_\varepsilon)$, so that (4.40) also holds for $\{y_\varepsilon\}$, (with $2D_o$ in place of $D_o$).

Now if $\{y_\varepsilon\}$ is any sequence of points in $U_\varepsilon$ as above satisfying (4.16)-(4.17) and (4.40) such that

$$(\rho^2(y_\varepsilon)/\varepsilon)\tau(y_\varepsilon) \to -\infty, \tag{4.43}$$

then the blow-up limit based at $\{y_\varepsilon\}$ is of type (III)$c_{(III)}$ above, i.e. a complete static vacuum solution, for the same reasons as discussed there.

Thus it remains to understand the situation where (4.40) holds at $\{y_\varepsilon\}$ and

$$(\rho^2(y_\varepsilon)/\varepsilon)|\tau(y_\varepsilon)| \leq K, \tag{4.44}$$

for some $K < \infty$. This is addressed by the following result.

**Proposition 4.7.** *Let $U_\varepsilon$ be as in Lemma 4.4, and suppose that $U_\varepsilon$ satisfies in addition the following property as $\varepsilon = \varepsilon_i \to 0$:*
*(i). If $\sigma(M) < 0$, then $vol_{g_\varepsilon} U_\varepsilon \geq \nu_o$, for some arbitrary $\nu_o > 0$, independent of $\varepsilon$.*
*(ii). If $\sigma(M) = 0$, then $vol(\{\tau \leq -\nu_o\} \cap U_\varepsilon) \geq \nu_o$, for some arbitrary $\nu_o > 0$. Here, and only here, $\tau$ is not renormalized by its infimum $T$.*

*Fix an arbitrary $\delta > 0$ small, and let $y_\varepsilon$ be a sequence of points in the $\delta$-neighborhood of $U_\varepsilon$ almost realizing $\rho(U_\varepsilon) = \inf_{U_\varepsilon}\rho$, in the sense that $\rho(y_\varepsilon) \leq 2\rho(U_\varepsilon)$, as in Lemma 4.1. Then, as $\varepsilon = \varepsilon_i \to 0$, (in a subsequence),*

$$c_{g_\varepsilon} \rho^4(y_\varepsilon)/\varepsilon \to 0. \tag{4.45}$$

*Proof.* The proof is by contradiction, so suppose that (4.45) does not hold, i.e. $c_{g_\varepsilon} \rho^4(y_\varepsilon)/\varepsilon \geq 2d_1$, for some $d_1 > 0$, (on some subsequence $\varepsilon = \varepsilon_i \to 0$). Hence, since $\rho(y_\varepsilon)$ (almost) realizes the minimal value of $\rho$ on $U_\varepsilon$,

$$c_{g_\varepsilon} \rho^4(q_\varepsilon)/\varepsilon \geq d_1, \tag{4.46}$$

for all $q_\varepsilon \in U_\varepsilon$. By Lemma 4.6 and the minimality property of $y_\varepsilon$ again,

$$c_{g_\varepsilon} \rho^4(y_\varepsilon)/\varepsilon \leq 2D_o. \tag{4.47}$$

Next we claim that since (4.47) holds at $\{y_\varepsilon\}$, it then holds in all of $U_\varepsilon$, i.e.

$$c_{g_\varepsilon} \rho^4(q_\varepsilon)/\varepsilon \leq D_1, \tag{4.48}$$

for some constant $D_1 < \infty$, and all $q_\varepsilon \in U_\varepsilon$.

Assuming (4.48) for the moment, the proof is completed as follows. If both (4.46) and (4.48) hold, then the ratio $\rho(y_\varepsilon)/\rho(q_\varepsilon)$ is uniformly bounded, away from 0 and $\infty$, for all $q_\varepsilon \in U_\varepsilon$. This means that all sequences $\{q_\varepsilon\}$ in $U_\varepsilon$ are strongly $(\rho, d)$ buffered, c.f. §2, for some uniform constant $d > 0$, (independent of $\varepsilon$ and $\{q_\varepsilon\}$). Hence, from the strong convergence of Theorem 4.2, *all* blow-up limits based in $U_\varepsilon$ are complete and non-flat, (c.f. the discussion concerning (4.32)).

However, if $\sigma(M) < 0$, then $\tau_\varepsilon$ converges to $-1$ on a set of almost full measure in $U_\varepsilon$, by Theorem 3.11, c.f. (3.38) and assumption (i). This means that most blow-up limits in $U_\varepsilon$ have $\tau = const.$ $= -1$; note in particular that $\tau$ is not converging to 0 on such basepoint sequences. It follows of course from the form of the three equations in Proposition 4.4 that such limits are flat. This contradiction thus proves (4.45).

Similarly, if $\sigma(M) = 0$, then by (3.46), and assumption (ii), $\tau_\varepsilon$ approaches constant functions almost everywhere as $\varepsilon \to 0$, (w.r.t. volume), and such constants are bounded away from 0 on a definite percentage of base points. As above, it follows that there exist many blow-up limits with limit potential function $\tau = const. \neq 0$, and one has a contradiction as before.

Thus the main work is to prove (4.48) given (4.47). To do this, consider first the blow-up limit $(N, g', y)$ based at $y = \lim y_\varepsilon$. As noted above, $\{y_\varepsilon\}$ is buffered, so that the limit is non-flat. The minimality property of $\{y_\varepsilon\}$ implies that $(N, g', y)$ is complete. The limit equations for $(N, g')$ are the limit of the equations (4.35)-(4.36) as $\varepsilon = \varepsilon_i \to 0$. These are of the form

$$\nabla \mathcal{Z}^2 + L^*(\bar\tau) + cg = 0. \tag{4.49}$$

$$2\Delta(\bar\tau + \frac{s}{12}) = -\frac{1}{2}|z|^2 + 3c.$$

Note that $c = \lim c_{g_\varepsilon} \rho^4(y_\varepsilon)/\varepsilon < \infty$ by Lemma 4.6 and $c > 0$ by assumption, i.e. (4.46). Further, by (4.44), $\bar\tau = \lim (\rho^2(y_\varepsilon)\tau/\varepsilon)$ is locally bounded on $N$, although possibly identically 0. In this latter case, the equations (4.49) are the $\mathcal{Z}^2$ equations, with a non-zero constant term. We also observe that the term $\rho^2 s\tau/\varepsilon$ in (4.35)-(4.36) goes to 0 as $\varepsilon \to 0$, since $\rho^2 s\tau/\varepsilon = \rho^4 \tau^2 \sigma/\varepsilon$, and $\tau^2 << 1$, $\sigma < c_{g_\varepsilon}$.

Now we claim that there is a constant $\Lambda < \infty$, depending only on a lower bound for the value of $c$ in (4.49), such that

$$\rho(q, g') \leq \Lambda, \tag{4.50}$$

for all $q \in (N, g', y)$. For if (4.50) were false, then $(N, g')$ becomes flat near some divergent sequence $p_i \in N$. Suppose first that $\bar\tau(p_i)$ is bounded. It follows that any limit $(N_\infty, g'_\infty, p_\infty)$ of the pointed sequence $(N, g', p_i)$ is flat $\mathbb{R}^3$, (passing to covers in case of collapse), and the limit potential $\bar\tau_\infty$ is a solution to the equations

$$L^*(\bar\tau_\infty) = -cg, \quad 2\Delta(\bar\tau_\infty) = 3c > 0. \tag{4.51}$$

The only solutions of (4.51) on $\mathbb{R}^3$ are quadratic functions, with leading order term $4c \cdot |x|^2$. Since $\bar\tau \leq 0$ everywhere, this forces $c = 0$, which is impossible, giving (4.50) in this situation.

Suppose instead that $\bar{\tau}(p_i)$ is unbounded. We claim that whenever $\bar{\tau}(q)$ is sufficiently large, the metric $g'$ is necessarily almost flat in large balls about $q$, and $\bar{\tau}/|\bar{\tau}(q)|$ is almost constant in such balls, both w.r.t. the $L^{2,2}$ topology. To see this, renormalize the equations (4.49) by dividing by $|\bar{\tau}(q)|$. If $\bar{\tau}(q_i) \to -\infty$, it follows that any limit is a complete static vacuum solution, with renormalized potential $\bar{\tau}_{\infty} = \lim \bar{\tau}/|\bar{\tau}(q_i)| \leq 0$. By [AnI, Thm. 3.2], (c.f. also Theorem 5.1(III) below), this implies the limit is flat, and $\bar{\tau}_{\infty} \equiv -1$. Since the convergence to the limit is in the strong $L^{2,2}$ topology, this proves the claim. In particular, it follows that $|\bar{\tau}|$ has sublinear growth about any base point.

On the other hand, the arguments above concerning (4.51) then show that $\bar{\tau}$ is quadratic as above to leading order, whenever $|\bar{\tau}(q)|$ is sufficiently large. This again contradicts the fact that $\bar{\tau} \leq 0$ everywhere, and hence proves (4.50).

Now (4.50) implies that $\rho(q_{\varepsilon}, g_{\varepsilon}) \leq \Lambda \cdot \rho(y_{\varepsilon}, g_{\varepsilon})$ for all $q_{\varepsilon}$ within uniformly bounded $g_{\varepsilon}$-distance to $y_{\varepsilon}$, where $\Lambda$ depends only on a lower bound for the limit value of $c$ in (4.49). Observe that by (4.46), $c$ has a *uniform* lower bound on all of $U_{\varepsilon}$. Hence, this argument applies to all blow-up limits in $U_{\varepsilon}$ and shows that $\rho(q_{\varepsilon}, g_{\varepsilon}) \leq \Lambda \rho(y_{\varepsilon}, g_{\varepsilon})$, for all $q_{\varepsilon} \in U_{\varepsilon}$. This proves (4.48) and completes the proof.

∎

Observe that the hypotheses of Lemma 4.6 and Proposition 4.7 are always satisfied, (for example one may take $U_{\varepsilon} = \Omega_{\varepsilon}$), unless $\tau_{\varepsilon}$ is everywhere bounded away from 0 in $U_{\varepsilon}$. In this latter situation, as noted above, (4.45) holds automatically on sequences satisfying (4.16)-(4.17).

We organize the discussion above into the following

**Definition 4.8.** *A sequence of base points $y_{\varepsilon} \in (\Omega_{\varepsilon}, g_{\varepsilon})$, $\varepsilon = \varepsilon_i \to 0$, is called* preferred *if (4.16)-(4.17) hold at $y_{\varepsilon}$ and one of the following hold:*

*(i). $c_{g_{\varepsilon}} \rho^4(y_{\varepsilon})/\varepsilon \to 0$, as $\varepsilon = \varepsilon_i \to 0$.*

*(ii). There exists $D_o < \infty$ such that $c_{g_{\varepsilon}} \rho^4(y_{\varepsilon})/\varepsilon \leq D_o$ and $(\rho^2(y_{\varepsilon})/\varepsilon)\tau(y_{\varepsilon}) \to -\infty$.*

The previous results, Lemma 4.6 through Proposition 4.7 then imply the existence of preferred base points $y_{\varepsilon} \in (\Omega_{\varepsilon}, g_{\varepsilon})$ on a subsequence of any given sequence $\varepsilon = \varepsilon_i \to 0$, (under the standing assumption (4.12)). We do not discuss here the issue of the uniqueness of such base points $y_{\varepsilon}$; apriori, there could be a large variety of them.

Combining the results above in §4.2, it follows that the blow-up limits $(N, g', y)$ of $(\Omega_{\varepsilon}, g_{\varepsilon}, y_{\varepsilon})$, at preferred base points $y_{\varepsilon}$ are of one of the types in Cases (I)-(III).

Note that the decomposition (4.11) of $N$ is valid in all cases, as is (4.8). By (the proof of) Theorem 3.3 applied to $g'$, the limit metric $g'$ is locally in $L^{3,p}$, for any $p < \infty$, and of uniformly bounded curvature. The potential $\omega$ is also locally in $L^{3,p}$, while each summand $\tau, s$ is locally Lipschitz. The pair $(g', \omega)$ is a $L^{3,p}$ weak solution of the corresponding equations and $(g', \omega, s)$ are smooth away from $\Sigma$. By the same proof, the relation (3.25) holds on $\Sigma \subset N$.

Summarizing, with the above understood, we have proved the following result.

**Theorem 4.9.** *Any blow-up limit $(N, g', y)$ of $(\Omega_{\varepsilon}, g_{\varepsilon})$, $\varepsilon = \varepsilon_i \to 0$, with $y = \lim y_{\varepsilon}$ and $y_{\varepsilon}$ a preferred base point sequence, is a complete, non-flat solution either to the $\mathcal{Z}^2$ equations (4.6), the $\mathcal{Z}_c^2$ or $\mathcal{Z}_s^2$ equations (4.7),(4.9), or the static vacuum equations (4.10), possibly with junction to a $\mathcal{Z}^2$ solution.*

*In case the sequence $(\Omega_{\varepsilon}, g'_{\varepsilon}, y_{\varepsilon})$ is collapsing at $y_{\varepsilon}$, the limit $(N, g', y)$ has in addition a free isometric $S^1$ action.*

∎

We also obtain the following result on the structure of the junction set $\Sigma$, compare with Remark 3.6.

**Proposition 4.10.** *Let $(N, g)$ be a complete, non-flat $\mathcal{Z}_c^2$ solution with junction $\Sigma = \partial\{N^-\} = \partial\{N^+\}$, as in (4.11). Then $\Sigma$ has empty interior and is a $C^{1,\beta}$ surface in $N$, $\forall\beta < 1$. In particular, $\Sigma = \{s = 0\} \subset N$.*

*Proof.* The domain $N^- = \{\omega < 0\}$ is non-empty in $N$, with $\partial N^- \subset \Sigma$. If $\Sigma$ has non-empty interior, then there is an open set $U$, with $\omega \leq 0$ on $U$, whose intersection with $N^-$ and $\mathrm{int}(\Sigma)$ is non-empty. The strong maximum principle, c.f. [GT, Thm. 3.5], applied to the trace equation (4.7) implies that $\omega \equiv 0$ in $U$, i.e. $N^- \cap U = \emptyset$, a contradiction.

Similarly, since $\omega \in C^{2,\beta}$, the boundary maximum principle, c.f. [GT, Lem. 3.4], (again applied to the trace equation (4.7)), implies that $|\nabla\omega| \neq 0$ everywhere on $\Sigma$, so that $\Sigma$ is a $C^{1,\beta}$ surface. ∎

**Remark 4.11.** By construction, the $\mathcal{Z}_c^2$ and static vacuum solutions $(N, g')$ arising in Cases II, III(a),(b) above have potential function $\tau \leq 0$ uniformly bounded below on $N$. This may not be the case for the $\mathcal{Z}_c^2$ and static vacuum solutions arising in Cases III(c), since in this case, the construction of the limit requires renormalizing $\tau$ by factors diverging to infinity as $\varepsilon \to 0$.

## 5. Some Non-Existence Results.

In this section, we rule out a number of the apriori possible blow-up limits discussed above in Cases I - III of §4.

To begin, we quote the following results.

**Theorem 5.1. (I).** *Let $(N, g)$ be a complete non-compact $\mathcal{Z}^2$ solution, of non-negative scalar curvature. Then $(N, g)$ is flat.*

**(II).** *Let $(N, g)$ be a complete $\mathcal{Z}_s^2$ solution, (necessarily scalar-flat), with a free isometric $S^1$ action, and potential $\tau$ satisfying $-\lambda \leq \tau \leq 0$ everywhere, for some $\lambda < \infty$. Then $(N, g)$ is flat.*

**(III).** *Let $(N, g)$ be a complete static vacuum solution, with potential $\tau \leq 0$ everywhere. Then $(N, g)$ is flat.*

Parts I and II of Theorem 5.1 are proved in [An4, Thms.0.1,0.2], c.f. also [An4, Prop.5.4], respectively, while Part III is proved in [AnI, Thm.3.2(I)]. (Here we note that any complete $L^{2,2}$ static vacuum solution with $\tau \leq 0$ necessarily satisfies $\tau < 0$, by the weak maximum principle applied to the trace equation (4.10), c.f. [GT, Thm.8.1]).

Theorem 5.1(I) implies that Case (I) in §4.2 above, i.e. where $\varepsilon/\rho^2 \to \infty$, cannot occur. In fact the following Corollary is now essentially obvious.

**Corollary 5.2. (Maximal Scale).** *There is a constant $K_o = K_o(M)$, independent of $\varepsilon$, such that, for any $x_\varepsilon \in (\Omega_\varepsilon, g_\varepsilon)$,*

$$\frac{\varepsilon}{\rho^2(x_\varepsilon)} \leq K_o. \tag{5.1}$$

*Proof.* Suppose (5.1) is false, so that there exists a sequence $\{x_\varepsilon\}$ in $\{(\Omega_\varepsilon, g_\varepsilon)\}$ such that $\rho(x_\varepsilon) \to 0$ invalidating (5.1). By Lemma 4.1, we may choose points $y_\varepsilon$ satisfying (4.16)-(4.17), with in addition $\rho(y_\varepsilon) \leq 2\rho(x_\varepsilon)$. It follows that $\alpha_\varepsilon = \varepsilon/\rho^2(y_\varepsilon) \to \infty$. By Theorem 4.9, a subsequence of $(\Omega_\varepsilon, g_\varepsilon', y_\varepsilon)$ converges, (after passing to the universal cover in case of collapse), to a non-flat limit $(N, g', y)$, which is a complete, non-flat solution of the $\mathcal{Z}^2$ equations. This contradicts Theorem 5.1(I), and thus (5.1) must hold. ∎

We recall that in case $\sigma(M) = 0$, $\varepsilon$ is understood to be $\varepsilon/T$ in (5.1). This result of course shows that there is a maximal rate at which the curvature of $(\Omega_\varepsilon, g_\varepsilon)$ can blow up.

Next, we have the following strengthening of Theorem 5.1(II).

**Theorem 5.3.** *Let* $(N, g, \omega)$ *be a complete* $\mathcal{Z}_c^2$ *solution, with a free isometric* $S^1$ *action. Then* $(N, g)$ *is flat.*

*Proof.* When the potential $\omega = \tau + \frac{\varepsilon}{12}s$ is bounded below, this result can be proved by a minor modification of the proof of this result for $\mathcal{Z}_s^2$ solutions in [An4, Thm. 0.2], i.e. Theorem 5.1(II).

However, when $\omega$ is unbounded below, the proof requires different methods, which are first introduced in the proof of Theorem C. (Briefly, the methods of [An4, Thm.0.2] and those of Theorem C can be combined to prove that the $S^1$-invariant potential $\omega$ is necessarily a proper function onto its range, so that Theorem 5.3 is then a simple consequence of Theorem C).

For these reasons, we place the proof of Theorem 5.3 in Appendix C, where the results of Theorem C can be convieniently referred to.

∎

**Remark 5.4.** It is worth pointing out that Theorem 5.3 is false if the assumption that the $\mathcal{Z}_c^2$ solution admit a *free* isometric $S^1$ action is dropped, c.f. §6.1.

Next we turn to the generalization of Theorem 5.1(III) in this context.

**Theorem 5.5.** *Let* $(N, g', y)$ *be a complete blow up limit of* $(\Omega_\varepsilon, g'_\varepsilon, y_\varepsilon)$ *as in Theorem 4.9, which is either a complete static vacuum solution, or the junction of a static vacuum solution to a* $\mathcal{Z}^2$ *solution, (i.e. of the type in Cases III(a),(b), or* $(c)_{(III)}$). *Then* $(N, g)$ *is flat.*

*Proof.* By Theorem 5.1(III), we may assume that $(N, g')$, although complete, is not a complete static vacuum solution, and so the junction $\Sigma$ is non-empty. Hence, using Proposition 4.10, $N^+ = \{s > 0\}$, i.e. the region where $(N, g)$ is a $\mathcal{Z}^2$ solution, is non-empty. By Proposition 3.7 we have, on the sequence $(\Omega_\varepsilon, g'_\varepsilon, y_\varepsilon)$, at any $x \in \Sigma_\varepsilon$,

$$\nabla^-(\tau + \frac{\alpha}{12}s)(x) = \frac{\alpha}{12}\nabla^+ s(x) = \nabla^+ \omega, \tag{5.2}$$

where the notation is as in Proposition 3.7. By Remark 4.3(i), $\nabla s$ is uniformly bounded as $\varepsilon \to 0$. Since $\alpha \to 0$, it follows that on the blow-up limit $(N, g')$, we have, for all $x \in \Sigma$,

$$\nabla^+ \omega(x) = 0. \tag{5.3}$$

However, exactly as in the proof of Proposition 4.10, the boundary point maximum principle implies that $|\nabla^+ \omega|(x) > 0$, at every point $x \in \Sigma$, whenever $N^+$ is non-empty. Hence $N^+$ must be empty, and the result follows.

∎

The results above, namely Theorems 5.1(I), 5.3 and 5.5, imply in particular that any complete blow-up limit $(N, g', y)$, $y = \lim y_\varepsilon$, with $y_\varepsilon$ a preferred base point sequence, (c.f. Def. 4.8), which has a free isometric $S^1$ action, must be flat. However, if the blow-up sequence $(\Omega_\varepsilon, g'_\varepsilon, y_\varepsilon)$ is collapsing at $y_\varepsilon$, then Theorem 4.9 implies that such blow-up limits are complete and non-flat, and do have a free isometric $S^1$ action. Thus, the following Corollary is immediate.

**Corollary 5.6.** *Let* $y_\varepsilon \in (\Omega_\varepsilon, g_\varepsilon)$, $\varepsilon = \varepsilon_i \to 0$, *be a preferred base point sequence, with* $\rho(y_\varepsilon) \to 0$ *as* $\varepsilon \to 0$. *Then the manifolds* $(\Omega_\varepsilon, g'_\varepsilon, y_\varepsilon)$, $g'_\varepsilon = \rho(y_\varepsilon)^{-2} \cdot g_\varepsilon$ *cannot collapse at* $y_\varepsilon$ *and hence cannot collapse within* $g'_\varepsilon$-*bounded distance to* $y_\varepsilon$.

∎

The results in §4 and above in §5 lead now easily to the proof of Theorem B.

**Proof of Theorem B:**

Suppose the sequence $\{(\Omega_\varepsilon, g_\varepsilon)\}$ degenerates in the sense of (0.8), for some sequence $\varepsilon = \varepsilon_i \to 0$. As noted following Definition 4.8, the results from Proposition 4.4 through Proposition 4.7 imply that there exists, (in a subsequence), a preferred base point sequence $\{y_\varepsilon\}$, $y_\varepsilon \in (\Omega_\varepsilon, g_\varepsilon)$, with $\rho(y_\varepsilon) \to 0$ as $\varepsilon \to 0$.

From the analysis in §4.2, a subsequence of $\{(\Omega_\varepsilon, g'_\varepsilon, y_\varepsilon)\}$, $g'_\varepsilon = \rho(y_\varepsilon)^{-2} \cdot g_\varepsilon$, converges in the strong $L^{2,2}$ topology, in fact in the $L^{3,p}$ topology, to a complete non-flat limit Riemannian manifold $(N, g', y)$, passing to sufficiently large covering spaces in case the sequence $g'_\varepsilon$ collapses at $y_\varepsilon$. By Theorem 4.9, the limit metric is a solution of one of the equations discussed in Cases (I)-(III) in §4.2. By Corollary 5.6, the subsequence cannot collapse at $y_\varepsilon$, and hence must converge to the limit. Further, the blow-up limit cannot be of the type in Case (I) or Case (II)(a), i.e. a $\mathcal{Z}^2$ solution, by Theorem 5.1(I), nor of the types in Cases (III)(a),(b), $(c_{(III)})$, i.e. a static vacuum solution, possibly with junction, by Theorem 5.5, nor of the type Case (III)$(c_{(I)})$, (a $\mathcal{Z}^2$ solution), by Theorem 5.1(I) again.

Hence the limit is a complete $\mathcal{Z}^2_c$ solution, (with junction), or a (pure) $\mathcal{Z}^2_s$ solution, i.e. belongs to one of Case (II)(b),(c) or Case (III)$(c_{(II)})$. The remainder of the proof follows from the discussion in §4.2 on the structure of such limit solutions, summarized preceding Theorem 4.9, together with Proposition 4.10.

∎

## 6. Models for the Degeneration of the Minimizers.

In this section, (which is not required for the later developments), we discuss some examples which illustrate the results obtained in previous sections. Thus we present some explicit sequences of metrics conjectured to model the behavior of the minimizing metrics $g_\varepsilon$ of $I_\varepsilon^-$ on a large class of manifolds $M$ with $\sigma(M) \leq 0$, whose blow-up limits are solutions of the $\mathcal{Z}^2_s$ and $\mathcal{Z}^2_c$ equations. We also discuss the non-existence of both static vacuum limits (with junction), as well as collapsed solutions, in the context of these models.

These models also illustrate the validity Theorem C, in that they satisfy all of its assumptions and conclusions.

**§6.1. Example 6.1.** Let $M_1$ and $M_2$ be closed hyperbolic manifolds, and let

$$M = M_1 \# M_2. \tag{6.1}$$

The simplest possible picture of the degeneration of the minimizing sequence $\{g_\varepsilon\}$ on $M$ is that $g_\varepsilon$ converges smoothly to a hyperbolic metric $g_o$ away from a small geodesic ball $B_x(\mu)$, while in $B_x(\mu)$, the metrics $g_\varepsilon$ are forming a neck $S^2 \times \mathbb{R}$, which is being crushed to a pair of (topological) 3-balls $B^3$ glued at one point, say $x$, as $\varepsilon \to 0$. Thus, the core 2-sphere $S^2 \times \{0\} \subset S^2 \times \mathbb{R}$ is being crushed to the point $x$ as $\varepsilon \to 0$. If one rescales such metrics $g_\varepsilon$ by the size of the curvature at $x$, the metrics $g'_\varepsilon$ should converge to the Schwarzschild metric

$$g_S = (1 - \frac{2m}{r})^{-1} dr^2 + r^2 ds^2_{S^2(1)}, \text{ on } [2m, \infty) \times S^2, \tag{6.2}$$

doubled or reflected across the core $S^2 = t^{-1}(0) = r^{-1}(2m)$, or possibly a metric quasi-isometric to the double of $g_S$. Here $t$ is defined by $(1 - \frac{2m}{r})^{-1} dr = dt$, and the initial condition above. Note that the Schwarzschild metric is the canonical solution of the static vacuum Einstein equations, with potential function $u = (1 - \frac{2m}{r})^{1/2}$. The core $S^2 = \{u = 0\}$ is called the horizon of the solution.

Of course this description implicitly requires that $\sigma(M_i)$ is realized by the hyperbolic metric on $M_i$, which is unknown, c.f. [An2, §4], and similarly that $\sigma(M)$ is realized by the union of the hyperbolic metrics on the one point union (at $x$) of $M_1$ and $M_2$, i.e.

$$\sigma(M) = -(vol_{-1} M_1 + vol_{-1} M_2)^{2/3},$$

where $vol_{-1}$ denotes the hyperbolic volume; this relation is also only conjectured, c.f. [An3, §7]. It is known that the hyperbolic metric on each $M_i$ is a local minimum of $I_\varepsilon^-$, for $\varepsilon \geq 0$, ,c.f. [B,

Ch.4G], and it is likely that there exist metrics $\gamma_\varepsilon$ which are local minima of $I_\varepsilon^-$ on $M$, for which $I_\varepsilon^-(\gamma_\varepsilon)$ converges to the value above. All of the results of the previous sections would hold equally for the family $\{\gamma_\varepsilon\}$. In any case, regardless of these conjectures, we proceed to explicitly construct metrics which model the behavior above and illustrate Theorems B and C.

To start, as in Example 2 of [AnI, §6.2], we may bend the hyperbolic metric $g_o$ on $M_1 \setminus B_x(\delta)$ to a metric $g_\delta$, so that $g_\delta = g_o$ on $M_1 \setminus B_x(\mu)$ while the boundary $S = \partial(M_1 \setminus B_x(\delta)), \delta << \mu$, is a totally geodesic symmetric 2-sphere of radius $\delta$. (The parameter $\delta$ was called $a$ in [AnI, §6.2].)

More precisely, consider first the conformally flat metric

$$h_\delta = dt^2 + f^2(t)ds_{S^2(1)}^2, \tag{6.3}$$

on $[0, \infty) \times S^2$, where $f = f_\delta$ is the solution to the initial value problem

$$(f')^2 = 1 + f^2 - (\delta + \delta^3)f^{-1}, f(0) = \delta, f' \geq 0. \tag{6.4}$$

In this metric, the boundary $S^2 = \{t^{-1}(0)\}$ is a totally geodesic 2-sphere, of constant curvature $\delta^{-2}$, while the metric is asymptotic to the hyperbolic metric of curvature $-1$ for large $t$. In particular, as $\delta \to 0$, these metrics converge smoothly to the hyperbolic metric $H^3(-1)$ at any given distance away from the boundary $S^2$. The scalar curvature $s_\delta$ of $h_\delta$ is $s_\delta = -6$. In fact, the sectional curvatures $K_{ij}$ of $h_\delta$ are given by

$$K_{12} = K_{13} = -\frac{f''}{f} = -1 - \frac{1}{2}\frac{\delta + \delta^3}{f^3}, \tag{6.5}$$

$$K_{23} = \frac{1 - (f')^2}{f^2} = -1 + \frac{\delta + \delta^3}{f^3}.$$

Now let $\mu = \delta^{1/3}$. A simple estimate, using (6.4) and (6.5), shows that the curvature of $h_\delta$ in the annulus $A(\mu/2, \mu) = t^{-1}(\mu/2, \mu)$ is uniformly bounded, as $\delta \to 0$. Thus, the metric $h_\delta$ may be smoothly matched to the hyperbolic metric $g_o$ near the boundary $S^2(\mu)$, by bending the function $f_\delta$ suitably in the interval $(\mu/2, \mu)$. Note however that although $h_\delta$ has bounded curvature in $A(\mu/2, \mu)$, its curvature is not close to -1, i.e. $h_\delta$ is not close to the hyperbolic metric in this region.

With a slight abuse of notation, let $h_\delta$ denote the resulting (conformally flat) metric on $M_1 \setminus B_x(\delta)$. The smoothed metric $h_\delta$ has scalar curvature bounded everywhere, while

$$vol_{h_\delta}A(\mu/2, \mu) \leq c \cdot vol_{-1}A(\mu/2, \mu), \tag{6.6}$$

where $c$ is an absolute constant, independent of $\delta$.

An identical construction may be carried out on $M_2$ and the metrics may be matched smoothly at the common boundary $S$ to give a smooth metric $h_\delta$ on $M = M_1 \# M_2$.

We now fix an $\varepsilon > 0$ and estimate $I_\varepsilon^-(h_\delta)$ on the family $\{h_\delta\}$. First, it is easy to see that

$$\varepsilon \int_M |z|^2 \sim \frac{\varepsilon}{\delta}. \tag{6.7}$$

Here, $x \sim y$ means that their ratio is uniformly bounded, away from 0 and $\infty$ as $\delta \to 0$. For the scalar curvature integral, we have $s + 6 = 0$ off $B_x(\mu)$, while $s + 6$ is bounded in $B_x(\mu)$, but bounded away from 0 in $A(\mu/2, \mu)$, since the metric is not close to the hyperbolic metric. Thus

$$v^{1/3} \int_M (s_\delta^-)^2 - \sigma(M)^2 \sim vol_{h_\delta}B_x(\mu). \tag{6.8}$$

Note that $vol_{g_{-1}}B_x(\mu) \sim \mu^3 = \delta$. We claim that also

$$[vol_{h_\delta}B_x(\mu) - vol_{g_{-1}}B_x(\mu)] << \delta. \tag{6.9}$$

This may either be verified by direct computation, or seen as follows. The blow-up limit of the metrics $h_\delta$ is the Schwarzschild metric. It is easy to verify from (6.2) that if $B(R)$ is the domain $\{r \le R\}$ in the Schwarzschild metric, then $volB(R) = vol_{Eucl}B(R) + cR^2 + o(R^2)$, where $vol_{Eucl}$ denotes the volume of the $R$-ball in $\mathbb{R}^3$. Rescaling the $R$-ball to a ball of size $\mu$, (i.e. rescaling the Schwarzschild metric by $(\frac{\mu}{R})^2$), we obtain

$$vol_{h_\delta}B(\mu) \sim vol_{g_{-1}}B(\mu) + c(\frac{\mu}{R})^3 R^2 \sim vol_{g_{-1}}B(\mu) + c\frac{\delta}{R},$$

which gives (6.9). It follows that

$$v^{1/3} \int_M (s_\delta^-)^2 - \sigma(M)^2 \sim c \cdot \delta,$$

and thus

$$I_\varepsilon{}^-(h_\delta) \sim \frac{\varepsilon}{\delta} + c \cdot \delta + |\sigma(M)|. \tag{6.10}$$

This means that the value of $\delta$ minimizing $I_\varepsilon(h_\delta)$ is given (approximately) by

$$\delta = (c\varepsilon)^{1/2}. \tag{6.11}$$

Thus, abusing notation slightly, we define the metric $h_\varepsilon$ to be $h_\delta$, where $\delta = (c\varepsilon)^{1/2}$ and have thus

$$I_\varepsilon{}^-(h_\varepsilon) - |\sigma(M)| \sim \varepsilon^{1/2}. \tag{6.12}$$

Note that the maximal curvature of $h_\delta$ occurs at the core 2-sphere $S$ and is thus of the order $\delta^{-2}$. It follows that, for $x \in S$,

$$|r_{h_\varepsilon}|(x) \sim \rho_\varepsilon(x)^{-2} \sim \delta^{-2} = c\varepsilon^{-1}. \tag{6.13}$$

Thus, we see that for this choice of metrics,

$$\varepsilon/\rho^2 \to \alpha > 0, \tag{6.14}$$

where $\rho = \rho_\varepsilon(x)$. Further, the blow-ups of $\{h_\varepsilon\}$ at $x$ have limit given by the Schwarzschild metric (6.2), doubled or reflected across the horizon $S$. Note that the scalar curvature of $h_\varepsilon$, viewed as a family of functions on the blow-up limit, is identically $-6$.

For this discussion to be a good model of the possible behavior of the actual minimizers $\{g_\varepsilon\}$ of $I_\varepsilon$, it follows that the limit Schwarzschild metric should also satisfy the pure $\mathcal{Z}_s^2$ equations (4.9), i.e.

$$\alpha \nabla \mathcal{Z}^2 + L^*(\tau) = 0. \tag{6.15}$$

Here $\tau = \lim \tau_\varepsilon$ represents the limit of the scalar curvatures $s_\varepsilon{}^-$ of the minimizers $g_\varepsilon$, viewed as functions on the blow-up limit. It is easy to see that there is no non-flat solution to (6.15) with $\tau = $ const, so that the potential function $\tau$ does not represent the limiting scalar curvature of $h_\varepsilon$. Nevertheless, if this is an accurate model of $\{g_\varepsilon\}$, one expects $\tau$ to have similar features to the scalar curvature of $h_\varepsilon$; in particular, $\tau$ should be negative and asymptotic to a negative constant at each end of the doubled Schwarzschild metric.

The Schwarzschild metric satisfies the static vacuum equations, so that $L^*(u) = 0$, with $u = (1 - \frac{2m}{r})^{1/2}$. Thus, this function does not correspond to $\tau$; the harmonic function $u$ is naturally *odd* w.r.t. reflection across $\Sigma$, and this does not correspond to the expected behavior of $\tau$, which should be even across $\Sigma$.

Nevertheless, it is shown in [An4, Prop.5.1] that the Schwarzschild metric (6.2) satisfies the $\mathcal{Z}_s^2$ equations (6.15). The potential $\tau$ is spherically symmetric, satisfies $\tau < 0$ everywhere, and is even w.r.t. reflection across $\Sigma$. Explicitly $\tau$ is given by

$$\tau = lim_{a \to 2m} \tau_a \tag{6.16}$$

where $a > 2m$ and, for $r \geq a$,

$$\tau_a(r) = \frac{\alpha}{8}(1 - \frac{2m}{r})^{1/2} \left( \int_a^r \frac{1}{s^5(1 - \frac{2m}{s})^{3/2}} ds - \frac{1}{ma^3} \frac{1}{(1 - \frac{2m}{a})^{1/2}} \right). \qquad (6.17)$$

One may verify without difficulty that $\tau(2m) = \frac{\alpha}{8m}(2m)^{-3}$ and $\frac{d\tau}{dr}(2m) = 0$.

Thus the Schwarzschild metric, with potential function $\tau$, gives the simplest non-trivial solution to the $\mathcal{Z}_s^2$ equations (4.9), just as it is the canonical solution of the static vacuum equations (4.10). This discussion shows that the metrics $\{h_\varepsilon\}$ do provide a good model that might actually describe the degeneration of the minimizers $\{g_\varepsilon\}$ on $M = M_1 \# M_2$. In particular blow-up limits of the form described in Theorem B, may actually arise.

We make a number of further remarks on this construction.

**Remark 6.2.** The discussion above shows that the doubled Schwarzschild metric is a critical point of $\mathcal{Z}^2$ among asymptotically flat metrics of non-negative scalar curvature, and non-positive volume deformation, with a given mass. The minimizing property of $\{g_\varepsilon\}$ implies that a limit $g'$ of the blow-ups $g'_\varepsilon$ should actually realize the minimal value of $\mathcal{Z}^2$ among such metrics. We conjecture that in fact the doubled Schwarzschild metric does realize this minimal value, for metrics with two asymptotically flat ends, of equal mass.

**(ii).** Using the construction above, one can also model the degeneration of the minimizers $\{g_\varepsilon\}$ on connected sums of the form

$$M = P \# M_1, \qquad (6.18)$$

where $M_1$ is hyperbolic and $P$ is a closed 3-manifold with $\sigma(P) > 0$.

Namely, as noted in [AnI, Rmk. 6.2(i)], since $P$ admits metrics with $s = -6$, exactly the same construction as above gives a smooth metric $h_\delta$ on $M$ for which all the estimates (6.6)-(6.14) remain valid, as is easily checked.

Alternately, in place of the function $f_\delta$ defined by (6.4), define $f = f_\delta$ to be the solution to the initial value problem

$$(f')^2 = 1 - \delta f^{-1}, f(0) = \delta, f' \geq 0. \qquad (6.19)$$

and form the metric $h_\delta$ as in (6.3). This metric is scalar-flat, and gives the Schwarzschild metric of mass $m = \delta/2$. The same construction can be carried out with respect to this glueing metric on $S^2 \times \mathbb{R}$, giving the same estimates as before.

In either case, we have $\varepsilon/\rho^2 \to \alpha > 0$, and the blow-up limit is the doubled Schwarzschild metric, so that (6.15) holds. Here however, since $M_1$, on the 'right', is hyperbolic, while $P$, on the 'left', has $\sigma(P) > 0$, one might not expect $\tau$ to be even or symmetric across the horizon $S$. It is possible that $\tau$ is asymptotic to a negative constant on the right end, and asymptotic to a possibly different constant, perhaps $\tau \to 0$, on the left end of the doubled Schwarzschild metric. (Note that since $\sigma(P) > 0$, $N$ admits metrics of arbitrarily small scalar curvature and volume, so that the minimizer $g_\varepsilon$ may be (almost) of non-negative scalar curvature, with very small volume on $P$); c.f. again [AnI, Rmk.6.2(i)].

In fact, there are such solutions $\tau'$ to (6.15), obtained from $\tau$ of (6.16)-(6.17) by just readjusting the constant term after the integral in (6.17) suitably. This corresponds to adding suitable multiples of the potential function for the Schwarzschild metric $u = (1 - 2m/r)^{1/2} \in \text{Ker } L^*$, which is of course an odd function under reflection, c.f. [An4, Prop. 5.1].

**(iii).** Continuing the discussion on $M$ as in (6.18), there are other models of the blow-up limit which might better describe the degeneration of the minimizers $\{g_\varepsilon\}$ in this case. Namely, instead of the doubled Schwarzschild metric, the blow-up limit $(N, g')$ may be a complete manifold with *one* asymptotically flat end, topologically of the form $P \setminus \{pt\}$. The end $S^2 \times \mathbb{R}^+$ of $N = P \setminus \{pt\}$ should be glued onto a neighborhood of the end of $M_1 \setminus \{pt\}$.

Consider for instance taking the Schwarzschild metric (6.19) and, instead of doubling, glueing a metric of non-negative scalar curvature onto the horizon $S$. Thus, the blow-up limit $(N, g')$ may be a complete, asymptotically flat $\mathcal{Z}_c^2$ solution, given by the Schwarzschild metric outside the horizon $S$ and a $\mathcal{Z}^2$ solution 'inside' the horizon, with junction equal to $S$.

It is possible, (and even probable although unknown), that such a configuration may have smaller $L^2$ norm of curvature, (among asymptotically flat manifolds of non-negative scalar curvature), than the configuration in (ii) above.

**(iv).** In the converse direction, let $(N', g')$ be a complete manifold of non-negative scalar curvature, with a finite number of asymptotically flat ends $E_i, i = 1,....,k$. Then, as in [AnI,§6.1], one may blow this metric down, i.e. rescale large balls to very small size, and glue it smoothly onto hyperbolic manifolds $M_i' = M_i \setminus B^3$, $i = 1,...,k$, to form the closed manifold

$$M = N' \cup M_i, \tag{6.20}$$

where the union is along essential 2-spheres.

The same techniques as above give a family of metrics $h_\varepsilon$ on $M$, which converge to the one-point union, say at $x$, of the hyperbolic manifolds $M_i$, and which satisfy $\varepsilon/\rho^2 \to \alpha > 0$, $\rho = \rho_\varepsilon(x)$. The blow-ups at $x$ converge to the manifold $(N', g_o)$. As a concrete example, $(N', g')$ may be $S^3$ deleted at $k$ distinct points, $k \geq 2$ with a metric of non-negative scalar curvature, asymptotically flat at each end. We note that when $k \geq 3$, there are no static vacuum solutions of this form. We conjecture that there are complete $\mathcal{Z}_s^2$ or $\mathcal{Z}_c^2$ solutions of this form.

**Remark 6.3.** Suppose $(N, g', \omega)$ is a complete $\mathcal{Z}_c^2$ or $\mathcal{Z}_s^2$ solution, arising as a blow-up limit of $(\Omega_\varepsilon, g_\varepsilon, y_\varepsilon)$ as in Theorem B. Suppose further that

$$\tau \to -1, \tag{6.21}$$

uniformly at infinity in $N$. We will see in §7 that this implies that $\omega$ is a proper function (onto its image) on $N$, so that Theorem C implies that $N$ has a finite number of ends, each asymptotically flat.

To analyse this a bit further, consider the conformally equivalent metric on $N$ given by

$$\bar{g} = (1 + \tau)^4 \cdot g'. \tag{6.22}$$

Using standard formulas for the behavior of scalar curvature under conformal changes, c.f. [B, Ch.1J], and the fact that $g'$ is of non-negative scalar curvature, it is easily seen from the trace equation (4.9) that the scalar curvature $\bar{s}$ of $\bar{g}$ satisfies

$$\bar{s} \geq 2\alpha(1 + \tau)^{-5}|z|^2 \geq 0. \tag{6.23}$$

Further, since $\tau \sim -1 + \frac{m}{t}$, for some $m > 0$, (depending on the end), as in Theorem C, one has $\bar{g} \sim t^{-4} \cdot g'$ asymptotically, where $t$ is the distance to some base point. This implies that the manifold $(N, \bar{g})$ may be compactified by adding one point $p_k$ to each end $E_k$, as in the stereographic projection of $\mathbb{R}^3$ into $S^3$, c.f. [S],[LP]. The compactification $(\bar{N}, \bar{g})$ is a closed Riemannian manifold, and a computation shows that $\bar{g}$ is $L^{2,q}$ smooth across each $\{p_k\}$, for any $q < \infty$. The scalar curvature of $(\bar{N}, \bar{g})$ is non-negative by (6.23), and one has

$$\int_{\bar{N}} \bar{s}d\bar{V} \geq 8 \int_N |\nabla \tau|^2 dV'. \tag{6.24}$$

In particular, this implies that $\sigma(\bar{N}) > 0$.

It follows that when the limit $(N, g')$ is a complete $\mathcal{Z}_c^2$ solution satisfying (6.21), it is obtained from a closed 3-manifold $(\bar{N}, \bar{g})$ with $\sigma(\bar{N}) > 0$ by conformally multiplying $\bar{g}$ by $(\sum a_k G(x, p_k))^4$, where $G(x, p_k)$ is the Green's function for the conformal Laplacian and $a_k > 0$, c.f. again [S],[LP].

**§6.2.** Next we make a few remarks on the non-existence of static vacuum limits. Of course the main reason there are no such limits is that the potential function $\tau$ is constructed so that $\tau \leq 0$ everywhere and by Theorem 5.1(III), there are no non-flat complete static vacuum solutions with this property. In turn, the fact that $\tau \leq 0$ is one of the main reasons for choosing the functional $\mathcal{S}_-^2$, c.f. §1(i)-(ii).

Recall that in Example 6.1, the fact that $\varepsilon/\rho^2 \to \alpha > 0$ rests on (6.10), and thus on the fact that the difference $\mathcal{S}_-^2(h_\delta) - |\sigma(M)|$ is on the order of $c \cdot \delta$. Now before the gluing, the scalar curvature of $h_\delta$ is $-6$, while the volume of $B(\mu)$ is very close to the hyperbolic volume - the difference is *small* compared with $\delta$, c.f. (6.9). Thus, the estimate (6.10) arises from the glueing of the metrics in the band $A(\mu/2, \mu)$, and from the fact that the curvature of $h_\delta$ is merely bounded, and not close to the curvature $-1$ of the hyperbolic metric. (The metric (6.3) is close to the hyperbolic metric only on scales much larger than $\mu = \delta^{1/3}$).

If this construction could be improved, to give a glueing on $A(\mu/2, \mu)$ with the difference of $\mathcal{S}_-^2(h_\delta)$ and $|\sigma(M)|$ on $B(\mu)$ *small* compared with $\delta$, then one would obtain static vacuum limits, having a non-empty horizon. Of course the work in §4 and §5 indicates that this probably cannot be done in such a way as to minimize the curvature in $L^2$ among admissible comparisons; compare also with the examples in [AnI, §6.5].

In effect, we know of no realistic constructions accurately *modeling* the behavior of minimizers $\{g_\varepsilon\}$ whose blow-up limits are given by static vacuum solutions. This may be somewhat surprising, in view of the results in [AnI]. In fact, by the results of [AnI, §7] one would expect that relatively generic, (say $\delta$-dense in a suitable topology on $\mathbb{M}$, for any $\delta > 0$), minimizing sequences for $\mathcal{S}_-^2$ have all (non-collapsing) blow-up limits satisfying the static vacuum equations. Of course, the sequence $\{g_\varepsilon\}$ is not at all generic.

**§6.3** Finally, we discuss some examples which illustrate the differences between the collapse and non-collapse cases, and thus also between the cases $\sigma(M) = 0$ and $\sigma(M) < 0$.

Suppose first $G$ is any closed, irreducible graph manifold. Any such manifold has $\sigma(G) \geq 0$, with $\sigma(G) = 0$ if and only if $|\pi_1(G)| = \infty$, c.f. [An2, §4].

Somewhat surprisingly, it is known that connected sums of graph manifolds are still graph manifolds, c.f. [So], [Wa]. Further, if $G_1$ and $G_2$ are graph manifolds, then $\sigma(G_1 \# G_2) = 0$ exactly when either $\sigma(G_1) = 0$ or $\sigma(G_2) = 0$, c.f. [K].

Now graph manifolds are exactly the class of 3-manifolds for which

$$inf I_\varepsilon^- = inf I_\varepsilon = 0, \tag{6.25}$$

$\forall \varepsilon > 0$, c.f. §3.1. Thus, on an arbitrary closed graph manifold $M$, there is a sequence of unit volume metrics $\{g_i\}$ such that $\mathcal{S}_-^2(g_i) \to 0$ and $\mathcal{Z}^2(g_i) \to 0$; in fact, one may also choose $\{g_i\}$ so that the curvature goes to 0 in $L^\infty$. The metrics $g_i$ collapse $M$ completely along an F-structure, (or possibly along a sequence of F-structures), on $M$.

It follows that on any, possibly reducible, graph manifold $M$ with $\sigma(M) = 0$, there exist unit volume minimizing sequences $g_\varepsilon$, with $\mathcal{S}_-^2(g_\varepsilon) \to 0 = \sigma(M)$, with uniformly bounded curvature, i.e. the sequence $\{g_\varepsilon\}$ does not degenerate anywhere in the sense of (0.6). This despite the fact that $M$ may have essential 2-spheres, e.g. $M = G_k = \#_1^k G$, where $G$ is an irreducible graph manifold of infinite $\pi_1$. In particular, the minimizing sequence $\{g_\varepsilon\}$ does not crush essential 2-spheres to points, as in the examples in §6.1 above; the geometry of collapse, along circles and tori, outweighs any geometry associated with crushing of 2-spheres.

On the other hand, one can also find explicit sequences of metrics $\{g_\varepsilon\}$ on $M$ with $\mathcal{S}_-^2(g_\varepsilon) \to 0$, (or $\mathcal{S}^2(g_\varepsilon) \to 0$), which do crush the essential 2-spheres in $M$ to points, and collapse the irreducible components of $M$ along F-structures on these components, on a larger scale than the collapse of the 2-spheres. For example, on $M = G_2 = G \# G$, take the Schwarzschild metric $g_S(\delta)$ of mass $\delta$, doubled across the horizon $S$ to a complete scalar-flat metric. This metric has horizon $S$ a totally geodesic 2-sphere of curvature $\delta^{-2}$, and outside a small neighborhood of $S$, the metric is very close to two copies of $\mathbb{R}^3$ with the flat metric. As in Example 6.1 or Remark 6.2(ii), one may glue $g_S(\delta)$ onto two copies of $G \setminus B$, where $B$ is a 3-ball. On a larger scale than $\delta$, one may then collapse each copy of $G \setminus B$ along an F-structure, to obtain $\mathcal{S}^2_- \to 0$, as $\delta \to 0$.

The discussion above reflects the fact that the geometrization of (reducible) graph manifolds has some possible ambiguities. It is in fact not necessary to perform the sphere decomposition first in order to geometrize a graph manifold. Reducible graph manifolds can be cut solely along tori so that the components admit geometric structures; however, to do this, one must allow compressible tori. A minimizing sequence for $I_\varepsilon^-$ or $I_\varepsilon$, for any fixed $\varepsilon > 0$, will perform the geometrization of a graph manifold $M$ in this way, c.f. also [An6, §2].

This discussion illustrates, (but of course does not prove), the non-existence of collapsing blow-up sequences for $\{g_\varepsilon\}$ in Corollary 5.8.

**Remark 6.4.** Continuing the discussion in Remark 6.2, consider closed manifolds $M$ on which one expects to have a mixture of non-collapse and collapse behavior; for example

$$M = G \# M_1, \tag{6.26}$$

where $M_1$ is a closed hyperbolic 3-manifold and $G$ is a graph manifold with $\sigma(G) = 0$, (possibly reducible). Here, one expects the minimizing sequence $g_\varepsilon$ to converge to the hyperbolic metric $g_o$ on $M_1$ as $\varepsilon \to 0$ and collapse $G$ to a lower dimensional space; the curvature of $\{g_\varepsilon\}$ should remain bounded away from an essential $S^2$ connecting $G$ to $M_1$.

In this case, following the reasoning in §6.1 and §6.3 above, the blow-up limit $(N, g', x)$ topologically may be of the form

$$N = (D^2 \times S^1) \# \mathbb{R}^3, \tag{6.27}$$

where $D^2 \times S^1$ connects to $G$ with a solid torus removed and the $\mathbb{R}^3$ factor connects to $M_1 \setminus B^3$. On the $\mathbb{R}^3$ factor, the metric $g'$ should be asymptotically flat, i.e. Schwarzschild-like, as in §6.1, while on the $D^2 \times S^1$ factor, $g'$ should collapse in the large, i.e. the volume radius $\nu$ should satisfy $\nu(x) << t(x)$ as $t(x) \to \infty$ in $D^2 \times S^1$, where $t$ is the distance to a base point.

In analogy to the second alternative in §6.3 above, another model for the blow-up limit of $M$ in (6.26) is

$$N = \mathbb{R}^3 \# \mathbb{R}^3, \tag{6.28}$$

in place of (6.27), with each end asymptotically flat. In this case, the collapse of $G$ takes place on a much larger scale than the crushing of the essential $S^2$.

As we will see next, Theorem C confirms that (6.28) is the correct model for the behavior of the blow-up limit.

## 7. Results on the Asymptotics of the Limits.

In this section, we begin the study of the topology and geometry of complete $\mathcal{Z}^2_c$ or $\mathcal{Z}^2_s$ solutions.

We start with a basically known result on the topology of complete open 3-manifolds with non-negative scalar curvature. Recall that $\mathcal{Z}^2_c$ solutions have non-negative scalar curvature; also recall from §2 that all manifolds are assumed to be connected and oriented.

**Proposition 7.1.** *Let $(N, g)$ be a complete, non-compact Riemannian 3-manifold with $s \geq 0$. Suppose that $N$ has finite topological type, i.e. $N$ is diffeomorphic to the interior of a compact manifold with boundary. Then $N$ is diffeomorphic to a connected sum of 3-balls, handlebodies $H_j$, and a closed 3-manifold $P$,*

$$N = (\#^q_1 B^3) \# P \# (\#^p_1 H_j). \tag{7.1}$$

*In particular, each end $E$ of $N$ is of the form $S^2 \times \mathbb{R}^+$ or $\Sigma_g \times \mathbb{R}^+$, where $\Sigma_g$ is a closed surface of genus $g \geq 1$, compressible in $N$.*

*Proof.* By results of Schoen-Yau [SY3, Thm.4], (c.f. also [GL2, Thm.8.4]), there are no closed oriented incompressible surfaces in $N$ of genus $g \geq 1$. Suppose first that $N$ has at least 2 ends. A standard result in 3-manifold topology, c.f. [J, Thm.III.11] then implies that $N$ has a closed oriented, (possibly disconnected), incompressible surface $\Sigma$, disconnecting $N$. Thus, each component of $\Sigma$ must be a 2-sphere $S^2$. If $\Sigma$ is connected, then performing surgery on $\Sigma$ thus implies that $N = N_1 \# N_2$, where each $N_i$ has at least one end. If $\Sigma$ is disconnected, then again by standard 3-manifold topology, c.f. [H, Lemma 3.8], $N = N_1 \# (S^2 \times S^1)$, so that $N_1$ still has at least 2 ends. If $N$ is of finite topological type, it follows by induction that

$$N = N_1 \# ... \# N_k \# P',$$

where each $N_i$ has exactly one end and $P'$ is a sum, (possibly empty), of $S^2 \times S^1$'s.

Assume now that $N$ has one end, and that $N$ is of finite topological type, so that $\partial N$ is a connected surface (at infinity). If $\partial N = S^2$, then $N = B^3 \# Q$, where $Q$ is a closed 3-manifold and $\partial N$ is incompressible exactly when $Q$ is non-empty, (or $Q \neq S^3$). If $\partial N \neq S^2$, then $\partial N$ must be compressible in $N$. By the sphere decomposition theorem, c.f. [H, Thm.3.15], we may write $N = N' \# Q'$, where $N'$ is irreducible and has one end, diffeomorphic to the end of $N$, while $Q'$ is closed and possibly reducible. Thus $N'$ is irreducible, with connected boundary, and has no closed incompressible surfaces. A standard result in 3-manifold topology, c.f. [J, Ex.IV.15], implies that $N'$ is a handlebody on a surface of genus $g \geq 1$. Reassembling this decomposition of $N$ proves (7.1). ∎

It seems likely that the closed factor $P$ in (7.1) is always a 3-manifold which admits a metric of positive scalar curvature, i.e. $\sigma(P) > 0$. We note that the only known closed 3-manifolds with positive scalar curvature metrics are connected sums of $S^2 \times S^1$ and space forms $S^3/\Gamma$.

Conversely, it is not difficult to prove that any manifold of the form (7.1) with $\sigma(P) > 0$ admits complete metrics of positive scalar curvature, c.f. [GL1], [SY2].

The remainder of this section is concerned with the geometry of $\mathcal{Z}_c^2$ or $\mathcal{Z}_s^2$ solutions. We start with the following elementary result, which shows that the geometry in a domain $D$ is controlled by the behavior of the potential function $\tau$, provided the oscillation osc $\tau$ of $\tau$ is sufficiently small.

**Lemma 7.2.** *Let $(N, g, \omega)$ be a $\mathcal{Z}_c^2$ or $\mathcal{Z}_s^2$ solution, and let $D$ be a smooth domain in $N$. Let $t_D(x) = dist(x, \partial D)$ and suppose there is a point $z \in D$ such that $\tau(z) < 0$. Then there is a small constant $\delta_o > 0$, independent of the solution $(N, g, \omega)$, and a constant $c > 0$, independent of $(N, g, \omega)$ and $\delta_o$, such that if*

$$osc_D \tau \leq \delta_o \cdot |\tau(z)|, \tag{7.2}$$

*then, for all $x \in D$,*

$$\rho(x) \geq c \cdot t_D(x), \quad \text{and} \quad |\nabla log|\tau||(x) \leq c^{-1}/t_D(x). \tag{7.3}$$

*Proof.* Under the assumption (7.2), suppose that (7.3) were false; then there must be a sequence $(D_i, g_i, \tau_i)$ of domains in $\mathcal{Z}_c^2$ or $\mathcal{Z}_s^2$ solutions, and points $z_i \in D_i$ with $\tau_i(z_i) = \chi_i < 0$, such that

$$osc_{D_i} \tau_i / |\chi_i| \to 0, \tag{7.4}$$

but

$$\rho(z_i)/t_D(z_i) \to 0. \tag{7.5}$$

Note that the ratio in (7.5) is scale invariant. Note further that the constant $\alpha = \lim \varepsilon/\rho(y_\varepsilon)^2$ in the $\mathcal{Z}_c^2$ equation (4.7) scales as the square of the distance, c.f. §4.1.

We now argue exactly as in the proof of Lemma 4.1. Choose points $y_i \in D_i$ realizing the minimum of the ratio (7.5), and rescale the metrics so that $\rho(y_i) = 1$. In this rescaling, (7.5) implies that $t_D(y_i) \to \infty$ and the metrics have bounded curvature in arbitrarily large neighborhoods of $y_i$. Hence, passing to universal covers in the case of collapse, (as discussed in §4), there is a subsequence converging smoothly to a non-flat limit $(N_\infty, g_\infty, \tau_\infty, y_\infty)$, where $\tau_\infty = \lim \tau_i / |\tau_i(y_i)|$. The smoothness of the convergence follows from Theorem 4.2/Remark 4.3(i). It follows from (7.4) that $\tau_\infty$ is a non-zero constant function.

By the discussion in §4.2, the limit is a solution either of the static vacuum, the $\mathcal{Z}_c^2$, or the $\mathcal{Z}^2$ equations; the last case occurs if $\alpha_i = \alpha \rho(y_i)^2 / |\tau_i(y_i)| \to \infty$. Since $\tau_\infty$ is constant, in the first two cases, the equations show that the limit is flat, while Theorem 5.1(I) gives the same conclusion in the last case. This is hence a contradiction in all cases, and so proves the first estimate in (7.3).

The regularity estimates from Theorem 4.2 or Remark 4.3(i) show that this estimate may be improved to

$$|r|(x) \leq C/t_D(x)^2, \tag{7.6}$$

for some constant $C < \infty$. The second estimate in (7.3) follows in the same way from these regularity estimates. ∎

**Remark 7.3.** If $(N, g, \omega)$ is a complete $\mathcal{Z}_s^2$ solution satisfying $\tau(z) < 0$ and

$$osc_N \tau \leq \delta_o \cdot |\tau(z)|, \tag{7.7}$$

then of course Lemma 7.2 implies that $(N, g)$ is flat, (by letting $t \to \infty$ in (7.3)). Thus, on any non-flat $\mathcal{Z}_s^2$ solution, the potential $\tau$ cannot be too close to a constant globally.

**(ii).** The estimate (7.3) also holds on the "original" family $(\Omega_\varepsilon, g_\varepsilon)$ of minimizers of $I_\varepsilon^-$, c.f. Remark 3.15. Thus, let $D = D_\varepsilon$ be a domain in $(\Omega_\varepsilon, g_\varepsilon)$ and $z_\varepsilon \in D$ a point such that $\tau(z_\varepsilon) < 0$. Then there is a small constant $\delta > 0$ and a constant $c < \infty$, independent of $\delta_o$ and $\varepsilon$, such that if $\varepsilon$ is sufficiently small and

$$osc_D \tau \leq \delta_o \cdot |\tau(z_\varepsilon)|,$$

then

$$\rho(x_\varepsilon) \geq c \cdot min(t_D(x_\varepsilon), 1), \text{ and } |\nabla log|\tau||(x_\varepsilon) \leq c^{-1} \cdot max(t_D(x_\varepsilon)^{-1}, 1) \tag{7.8}$$

for all $x_\varepsilon \in D$, where $t_D(x) = dist(x_\varepsilon, \partial D)$. The proof is identical to that of Lemma 7.2.

We now address the proof of Theorem C. For a first reading of the proof, it may be advisable to understand the statements of the collection of Lemmas and Propositions within the proof, but leave a detailed reading of their proofs until later, after an overall picture of the structure of the proof has been obtained. Generally speaking, Theorem C concerns the (unique) classification of the asymptotic behavior of $\mathcal{Z}_c^2$ solutions, under the hypotheses of the theorem.

**Proof of Theorem C.**

The proof is rather long, and will be separated later into several cases. We assume that $(N, g, \omega)$ is a complete, non-flat $\mathcal{Z}_c^2$ or $\mathcal{Z}_s^2$ solution satisfying

$$\omega \leq \omega_o < 0, \tag{7.9}$$

on $N \setminus K$, for some compact set $K \subset N$ and (arbitrary) constant $\omega_o < 0$, and that the levels $L_a = \{x \in N : \omega(x) = a\}$ of $\omega$ are compact in $N$. In particular $\omega = \tau$ and $g$ is scalar-flat on $N \setminus K$.

We first prove several general estimates on such $\mathcal{Z}_c^2$ solutions, expressed in Lemmas 7.4-7.8 below. Following this, we need to divide into cases, according to the specific possible asymptotic behavior of $(N, g)$.

To begin, observe that by the minimum principle applied to the trace equation (4.7), on any bounded domain $D \subset N$, $inf_D \omega$ occurs on $\partial D$. In particular, $inf_{S(s)} \omega$ is a monotone decreasing function of $s$.

The next result shows that the curvature of $(N, g)$ decays to 0 at infinity in $(N, g)$, (so that essentially $|r|$ is a proper function onto its image in $\mathbb{R}^+$).

**Lemma 7.4.** *Let $(N, g, \omega)$ be a complete $\mathcal{Z}_c^2$ solution satisfying (7.9), such that all level sets of $\omega$ are compact. Let $t(x) = dist(x, y)$, for a fixed base point $y \in N$. Then*

$$limsup_{t \to \infty} |r| = 0. \tag{7.10}$$

*Proof.* Suppose there is a divergent sequence $\{x_i\}$ in $N$ which violates (7.10), so that

$$|r|(x_i) \geq \delta_o, \tag{7.11}$$

for some $\delta_o > 0$. By passing to a subsequence, we may assume w.l.o.g. that $\{\omega(x_i)\}$ is either monotone decreasing to $-\infty$, monotone decreasing to a value $\omega_- \leq \omega_o$, or monotone increasing to a value $\omega_+ \leq \omega_o$. We treat each of these cases separately.

(i). Suppose $\omega(x_i)$ decreases to $-\infty$. Suppose further for the moment that the curvature of $(N, g)$ is uniformly bounded. Renormalize $\omega$ by setting $\omega_i(x) = \omega(x)/|\omega(x_i)|$, and renormalize the $\mathcal{Z}_c^2$ equations (4.7) by dividing by $|\omega(x_i)|$. Thus, $\alpha$ is replaced by $\alpha_i = \alpha/|\omega(x_i)| \to 0$, as $i \to \infty$, and hence the metric approaches a static vacuum solution in neighborhoods of $\{x_i\}$. By the same arguments used in §4.2 to prove convergence, the pointed sequence $(N, g, x_i)$ thus (sub)-converges, (passing to covers in case of collapse), to a complete static vacuum solution $(N_\infty, g_\infty, x_\infty)$, $x_\infty = \lim x_i$, with $\omega \leq 0$ everywhere. Theorem 5.1(III) implies the limit is flat. The convergence to the limit is smooth by Remark 4.3(i) and hence (7.11) passes to the limit to give $|r_{g_\infty}|(x_\infty) \geq \delta_o$. This contradiction implies that (7.10) must hold.

If $|r|$ is not uniformly bounded on $N$, (near $x_i$), then rescale to make it bounded, exactly as in the proof of Lemma 4.1 and apply the same argument as above. (Of course $\mathcal{Z}_c^2$ solutions arising from Theorem B have uniformly bounded curvature.)

Before dealing with the remaining two cases, we discuss some general features of $(N, \omega)$. Each compact level $L_a$ separates $N$ into the sub and super levels of $\omega$, so that

$$N = U_a \cup L_a \cup U^a,$$

where $U_a = \{x \in N : \omega(x) \leq a\}$, $U^a = \{x \in N : \omega(x) \geq a\}$. By the minimum principle applied to the trace equation (4.7), for any given $a$, every component of $U_a$ is non-compact, and thus $U_a$ determines a non-empty collection of ends $E_a^-$ of $N$, whose boundary $\partial E_a^-$ is contained in $L_a$. The remaining ends $E_a^+$ of $N$ (w.r.t. $L_a$), if any, are given by non-compact components of $U^a$. Since $L_a$ is compact, there are only a finite number of components of $L_a$ which are contained in $E_a^-$ but are not contained in $\partial E_a^-$. Such components bound a finite number of compact domains in $E_a^-$. Hence, sufficiently far out in $E_a^-$, all points in $E_a^-$ are in $U_a$. Since this holds for any $a$, we see that, roughly speaking, $\omega$ is decreasing in the ends $E_{a_o}^-$, modulo compact sets, for any fixed $a_o$.

(ii). Suppose $\{\omega(x_i)\}$ decreases to $\omega_- > -\infty$. Choose $a_o > \omega_-$, but sufficiently close to $\omega_-$, depending on $\delta_o$ in (7.11). Then a tail end of $\{x_i\}$ lies in $E_{a_o}^-$, (for some subsequence), and so, outside a compact subset of $E_{a_o}^-$, $\omega(x) < a_o$. If there were any points $y_i$ in $B_{x_i}(s_o)$, for a fixed $s_o < \infty$, such that $\omega(y_i) < \omega_-$, then applying the arguments above to $L_{\omega_-}$ implies that a tail end of $\{x_i\}$ must belong to $E_{\omega_-}^-$, which contradicts the assumption on $\omega(x_i)$. It follows that the oscillation of $\omega$ in $B_{x_i}(s_o)$ is at most $\kappa = |a_o - \omega_-|$, for $i$ sufficiently large. Choosing $\kappa$ sufficiently small and and $s_o$ sufficiently large, and applying (7.6) contradicts (7.11).

(iii). Suppose $\{\omega(x_i)\}$ increases to $\omega_+ \leq \omega_o < 0$, so that a tail end of $\{x_i\}$ lies in $E_{a_o}^+$, for any $a_o < \omega_+$. Hence, for the same reasons as in (ii), $\omega(x) > a_o$ outside a compact subset of $E_{a_o}^+$ and the oscillation of $\omega$ in $B_{x_i}(s_o)$ may be made arbitrarily small by choosing $a_o$ sufficiently close to $\omega_+$. As above, (7.6) then gives a contradiction to (7.11).

■

We assume for the remainder of the proof that (7.9) and (7.10) hold, but, with the single exception of Lemma 7.8 below, do not use any further the assumption that the levels of $\omega$ are compact in $N$.

An argument similar to the proof of Lemma 7.2 now leads easily to the following strengthening of (7.10).

**Lemma 7.5.** *Let* $(N, g, \omega)$ *be a complete* $\mathcal{Z}_c^2$ *solution, satisfying (7.9) and (7.10), with base point* $y$. *Then there is a constant* $\rho_o > 0$ *such that*

$$\rho(x) \geq \rho_o \cdot t(x). \tag{7.12}$$

*where* $t(x) = dist(x, y)$, *as in Lemma 7.4. Consequently, there is a constant* $C < \infty$ *such that for* $t$ *large,*

$$|r|(x) \leq \frac{C}{t^2(x)}, |\nabla log||\omega||(x) \leq \frac{C}{t(x)}. \tag{7.13}$$

*Proof.* The proof is a simple combination of Lemma 7.4 and the proof of Lemma 7.2. Thus, if (7.12) were false, then there exists a sequence $x_i$ in $N$ such that $\rho(x_i)/t(x_i) \to 0$ and $t_i = t(x_i) \to \infty$. Note that by Lemma 7.4, $\rho(x_i) \to \infty$ also. As in Lemma 4.1, we may choose $x_i$ so that it realizes approximately the minimal value of $\rho/t$ on, for instance, the geodesic annulus $A(\frac{1}{2}t_i, 2t_i)$. In other words, (re)-choose $x_i$, if necessary, so that it realizes the minimal value of the ratio

$$\rho(x)/dist(x, \partial A(\frac{1}{2}t_i, 2t_i)),$$

for $x \in A(\frac{1}{2}t_i, 2t_i)$. Again as in Lemma 4.1, such a choice of $x_i$ implies that $x_i$ is strongly $(\rho, \frac{1}{2})$ buffered, c.f. (2.3). Now consider the pointed sequence $(N, g_i, x_i)$, where $g_i = \rho(x_i)^{-2} \cdot g$. In this scale, $\rho_i(x_i) = 1$ and the curvature of $(N, g_i)$ is uniformly bounded in arbitrarily large $g_i$-distances to $x_i$. Hence, as in the proof of Lemma 7.2, a subsequence converges smoothly, (passing to covers in case of collapse), to a limit $(N_\infty, g_\infty, x_\infty)$, with $\rho_\infty(x_\infty) = 1$. Further, $\alpha_i = \rho(y)^2/\rho(x_i)^2 = \rho(x_i)^{-2} \to 0$, and so the limit is a static vacuum solution near $x_\infty$. The limit is thus both complete and a complete static vacuum solution, (i.e. there is no junction $\Sigma$), since $t_i(x_i) \to \infty$, and so by (7.9) $\omega < 0$ in arbitrarily large distances to $x_i$. Renormalizing $\omega$ as in the proof of Lemma 7.2 if $\omega(x_i) \to -\infty$ gives a limit potential $\omega_\infty \leq 0$ everywhere. Hence by Theorem 5.1(III), the limit $(N_\infty, g_\infty)$ is flat. This contradicts the fact that $\rho_\infty(x_\infty) = 1$.

The estimates (7.13) follow from strong convergence, (Theorem 4.2/Remark 4.3(i)), as in Lemma 7.2.

■

The estimates (7.13) together with the proof of Lemma 7.5 above imply that the $\mathcal{Z}_c^2$ solution $(N, g)$ is asymptotically a, possibly flat, static vacuum solution, (although the asymptotic behavior apriori may be far from unique). In [An5, Thm. 0.3], it is proved that an end of a static vacuum solution is either asymptotically flat in the sense of (0.17)-(0.18) or parabolic, provided for instance the horizon is compact; (the parabolic ends roughly have at most quadratic volume growth, so that such ends are 'small' in a natural sense). Hence, it is not surprising that the remainder of the proof of Theorem C is closely related to that of [An5, Thm.0.3].

It is convenient to change sign and set

$$u = -\omega = -\tau > 0, \quad \text{on} \quad N \setminus K. \tag{7.14}$$

Consider the 4-manifold $X = (N \setminus K) \times S^1$, with warped product metric

$$g_X = g + u^2 d\theta^2. \tag{7.15}$$

From standard formulas, c.f. [B,Ch.9J], the Ricci curvature of $g_X$ in horizontal directions, i.e. tangent to $N$, is given by

$$r_X = r_N - u^{-1}D^2u,$$

while in the fiber or $S^1$ direction of unit length, it is given by

$$r_X = -u^{-1}\Delta u.$$

The $\mathcal{Z}_c^2$ equations (4.7) and the decay estimate (7.13), thus imply that

$$|r_X| \leq c/t_X^4, \tag{7.16}$$

where $t_X(x) = dist_X(x,y)$, for a given base point $y \in X$. This is a quantitative version of the statement above that $(N, g)$ is asymptotically a static vacuum solution.

By [Lu], (c.f. also [An5, Lemma 3.1]), the estimate (7.16) implies that $X$, and hence $N$, has *a finite number of ends $E_i$*. Further, the intrinsic diameter of any component $A_c(s, 2s)$ of a geodesic annulus $A(s, 2s) = t^{-1}(s, 2s)$ in $N$, (or in $X$), satisfies

$$diam_{A_c}(s, 2s) \leq d_o \cdot s, \tag{7.17}$$

for some constant $d_o < \infty$, independent of $s$ and the component $A_c$. In addition, there is a uniform bound $D_o$ on the number of components $A_c(s, 2s)$ of $A(s, 2s)$, as well as a uniform bound on the number of components $S_c(s)$ of geodesic spheres $S(s)$ in $N$ or $X$.

Note that (7.13) and (7.17) imply that

$$area_N S(s) \leq d_1 \cdot s^2 \text{ and } vol_N B(s) \leq d_1 \cdot s^3, \tag{7.18}$$

for some $d_1 < \infty$. It is easiest to verify the scale-invariant estimate (7.18) w.r.t. the rescaled metrics $g_s' = s^{-2} \cdot g'$ on any $A_c(\frac{1}{2}s, 2s) \subset (N, g')$. Namely the scale-invariant estimates (7.13) and (7.17) imply that these annuli, in the $g_s$ metric, have uniformly bounded curvature and diameter. Hence the volume and the area of the cross-section $S(s)$ is uniformly bounded in the $g_s$ metric, which gives (7.18).

The scale-invariant estimates (7.13) and (7.17) give quite strong control over the *tangent cones at infinity* of any given end $E \subset N$. We recall the definition of such tangent cones. Consider any divergent sequence $\{x_i\}$ in a given end $E$, and the associated pointed manifolds $(E, g_i, x_i)$, where $g_i = t(x_i)^{-2} \cdot g$. By (7.13) and scale-invariance, this sequence has uniformly bounded curvature on $A_i = A_c(\kappa_i^{-1}t_i, \kappa_i t_i)$, for $t_i = t(x_i)$, where $A_c$ denotes the component of the geodesic annulus containing the base point $x_i$, and where $\kappa_i$ diverges to $\infty$ sufficiently slowly as $i \to \infty$. Note that (7.17) and (7.18) imply that the diameter and volume of the geodesic spheres $(S_i(s), g_i) = t_i^{-2}(S_y(st_i), g)$ are uniformly bounded as $i \to \infty$, for any bounded $s$, and $diam_{g_i}S_i(s) \leq d_o \cdot s$. If this sequence is collapsing, we may unwrap this collapse, as discussed in §4.2 and Remark 2.1. Hence, a subsequence converges to a (maximal) limit $(T_\infty, g_\infty, x_\infty)$ called the tangent cone at infinity of $(E, g)$ associated with the subsequence $\{x_i\}$. By Theorem 4.2/Remark 4.3(i) and (7.9), the convergence to $T_\infty$ is in the $C^\infty$ topology, and is uniform on compact subsets. Observe also that $(T_\infty, g_\infty, x_\infty)$ is the Gromov-Hausdorff limit of $\{(A_i, g_i, x_i)\}$, c.f. [Gr1, Ch.3], and that $T_\infty$ has a distinguished distance function $t = t_\infty$, namely the limit of the renormalized distances $t/t(x_i)$. While this construction and terminology are standard, note that the spaces $(T_\infty, g_\infty)$ are not necessarily metric cones.

The following corollary of the discussion above is now essentially obvious, but expresses significant initial control on the asymptotic behavior of $N$.

**Corollary 7.6.** *On $(N, g, \omega)$ as in Lemma 7.5, there is a constant $C < \infty$ such that*

$$\int_N |r|^2 dV \leq C. \tag{7.19}$$

*Proof.* This is an immediate consequence of the decay estimate (7.13) and the volume growth estimate (7.18), by integrating in exponential normal coordinates and applying Fubini's theorem. ∎

By the trace equation (4.9) and (7.13) again,

$$0 \leq \Delta_N u = \frac{\alpha}{4}|z|^2 \leq c/t^4,$$

on $N \setminus K$. Thus $u$ is subharmonic on $N \setminus K$ and asymptotically harmonic in a strong sense. In particular, (7.19) implies that

$$\Delta_N u \in L^1(N \setminus K). \tag{7.20}$$

Next we observe that $|\nabla u|^2$ satisfies a sub-mean value inequality on $(N, g)$.

**Lemma 7.7.** *There is a constant $C < \infty$ such that*

$$sup_{S_c(s)}|\nabla u|^2 \leq \frac{C}{vol A_c(\frac{1}{2}s, 2s)} \int_{A_c(\frac{1}{2}s, 2s)} |\nabla u|^2 + \frac{C}{s^4}, \tag{7.21}$$

*for all $s$ large, where $S_c(s)$ is any component of the geodesic sphere in $(N, g)$ about $y$; $C$ is independent of $s$ and the choice of component.*

*Proof.* As in the proof of (7.18), it is convenient to work in the scale $g_s = s^{-2} \cdot g$, so that the geodesic annulus $A_c(\frac{1}{2}s, 2s)$, (i.e. the component containing $S_c(s)$), becomes the annulus $A_c(\frac{1}{2}, 2)$ w.r.t. $g_s$. By passing to suitable covering spaces, we may also assume that $A_c(\frac{1}{2}, 2)$ has uniformly bounded geometry, independent of $s$. Note that the estimate (7.21) is invariant under coverings. The Bochner-Lichnerowicz formula, ([B, 1.155]), and the Cauchy-Schwartz inequality imply that

$$\Delta|\nabla u| \geq -|\nabla \Delta u| - |r||\nabla u|,$$

in the $g_s$ metric. By the trace equation (4.9) and the fact that $\alpha$ scales as the square of the distance, it follows then from the decay estimate (7.13) that

$$\Delta|\nabla u| \geq -c|\nabla u| - \frac{c}{t^2}.$$

The sup estimate (7.21) then follows from the DeGiorgi-Nash-Moser sup estimate applied to the elliptic inequality above, c.f. [GT,Thm.8.17]. ∎

This result leads to the following estimate on $(N, g)$:

**Lemma 7.8.** *There is a constant $C$, independent of $s$ and the choice of component, such that*

$$sup_{S_c(s)}|\nabla u| \leq C \cdot area S_c(s)^{-1} + C \cdot s^{-2}. \tag{7.22}$$

*Proof.* By the coarea formula, setting $A_c = A_c(\frac{1}{2}s, 2s)$, we have

$$\int_{A_c} |\nabla u|^2 \leq osc_A u \cdot max_v \left(\int_{A_c \cap L_v} |\nabla u|\right),$$

where $L_v$ is the $v$-level of $u$, and $v$ ranges over the values of $u$ in $A_c$. By (7.17), we have $osc_{A_c} u \leq s \cdot sup_{A_c}|\nabla u|$. Further, there is a constant $k < \infty$, independent of $A_c$, such that $sup_{A_c}|\nabla u| \leq k \cdot sup_{S_c}|\nabla u|$. This follows exactly as in the proof of Lemma 7.7, working in the scale $g_s$ where $(A_c, g_s)$ has uniformly bounded geometry. For the same reasons $s/vol A_c \leq k \cdot (vol S_c(s))^{-1}$.

Hence, (7.22) follows from (7.21), given one has a uniform bound

$$\int_{L_v \cap A_c} |\nabla u| \leq k,$$

for some constant $k < \infty$, independent of $A_c$ and $v$. However, since the levels of $u$ are compact in $N$, this follows easily from (7.20) by the divergence theorem. ∎

Given the preceding general results on the structure of $(N, g, u)$, we now separate the proof into five cases, only the first of which is proved to be possible. (This division parallels closely the division in the proof of [An5, Thm.0.3]).

We work on each one of the finite number of ends separately. Thus, fix an end $E$ of $N$ and let $S_E(s) = \cup S_c(s)$ for $S_c(s) \subset E$ defined as above. Since $vol S_E(s) = \Sigma vol S_c(s)$, we have

$$\int^{\infty} area S_E(s)^{-1} ds = \infty \Leftrightarrow \int^{\infty} area S_{\sigma}(s)^{-1} ds = \infty,$$

where $\sigma$ is some geodesic ray in $E$ and $S_{\sigma}(s)$ is the component of $S_E(s)$ containing $\sigma$. Let $E_{\sigma} = \cup S_{\sigma}(s)$, for $s \geq s_o$, for some initial $S_{\sigma}(s_o) \subset E$.

**Case I.** Assume

$$\int^{\infty} area S_E(s)^{-1} ds < \infty. \tag{7.23}$$

It follows from (7.22) that $u$ is bounded on $E_{\sigma}$, for some geodesic ray $\sigma \subset E$, and in fact $lim_{t \to \infty} u$ exists, where the limit is taken in $E_{\sigma}$.

Since $u$ is asymptotically constant in $E_{\sigma}$, the proof of Lemma 7.5 implies that the curvature $|r|$ of $(E, g)$ decays faster than quadratically, i.e.

$$|r|(x) << t^{-2}(x) \tag{7.24}$$

as $t(x) \to \infty$ in $E_{\sigma}$. Equivalently, since all tangent cones at infinity $T_{\infty}$ in $E_{\sigma}$ are static vacuum solutions with constant potential $u$, they must be flat. This, together with the global control (7.17)-(7.18) implies that $E = E_{\sigma}$, (c.f. [An5, Lem. 3.2] for further details if needed). In particular, there is only one component of $S(s) = S_E(s)$ in $E$, for $s$ large and

$$lim_{t \to \infty} u = u_o > 0, \tag{7.25}$$

where the limit is taken in $E$; the inequality follows from (7.9).

Next, we need to divide the discussion into two further cases, according to whether $E$ is non-collapsing or collapsing at infinity.

**Case I(A). (Non-Collapse).** Under the assumption (7.23), suppose

$$liminf_{t \to \infty} \frac{v_E(s)}{s^3} > 0, \tag{7.26}$$

where $v_E(s) = \mathrm{vol} B_E(s)$. Via (7.22), (and the curvature bound (7.13)), this of course implies

$$sup_{S_E(s)} |\nabla u| \leq c \cdot s^{-2}. \tag{7.27}$$

Now the estimate (7.24), together with (7.17)-(7.18) and the (non-collapse) volume condition (7.26) implies that that $(E, g)$ has a unique tangent cone at infinity $T_{\infty}$, given by $\mathbb{R}^3$ with the flat metric. Hence, there is compact set $K$ such that $E \setminus K$ is diffeomorphic to $\mathbb{R}^3 \setminus B(R)$, for some $R < \infty$ and the metric $g$ is (weakly) asymptotically flat in the sense that

$$g_{ij} = \delta_{ij} + \gamma_{ij},$$

in a chart for $E \setminus K$, and $|\gamma_{ij}| \to 0$ as $|x| \to \infty$ in $\mathbb{R}^3$.

To obtain better control on the behavior of the metric at infinity, consider the metric

$$\widetilde{g} = u^2 \cdot g. \tag{7.28}$$

By (7.25), $\widetilde{g}$ is quasi-isometric, and asymptotically homothetic to $g$ in $E \setminus K$. A standard computation, c.f. [B, Ch.1J], shows that $\widetilde{g}$ has Ricci curvature satisfying

$$\widetilde{r} = r - u^{-1} D^2 u + 2(d log u)^2 - u^{-1} \Delta u \cdot g. \tag{7.29}$$

By the $\mathcal{Z}_s^2$ equation (4.9), (7.13) and (7.25), it follows that

$$|\widetilde{r}| \leq c/t^4, \tag{7.30}$$

compare with (7.16). Since the Ricci curvature determines the full curvature in dimension 3, and the curvature is computed from the $2^{nd}$ derivatives of the metric, this implies that the metric $\widetilde{g}$ differs from the flat metric $\delta_{ij}$ by terms of order at most $t^{-2}$, i.e. in a suitable (harmonic) coordinate chart for $E \setminus K$,

$$\widetilde{g}_{ij} = \delta_{ij} + O(t^{-2}), \tag{7.31}$$

c.f. [BKN] for further details. Elliptic regularity for the $\mathcal{Z}_s^2$ equations implies that the $C^0$ estimate (7.31) holds also in the $C^k$ topology, i.e. $|\partial^k \widetilde{g}_{ij}| = O(t^{-2-k})$, for any $k < \infty$.

Similarly, from standard formulas, c.f. [B, Ch.1J],

$$\widetilde{\Delta} log u = u^{-3} \Delta u = u^{-3} \frac{\alpha}{4} |z|^2 = O(t^{-4}),$$

where the latter estimates follow from the trace equation (4.9) and from (7.13). Let $\Delta_o$ denote the Laplacian w.r.t. the flat metric on $\mathbb{R}^3 \setminus B$. Then we obtain from this and (7.31) that

$$|\Delta_o log u| \leq c \cdot (t_o^{-2})|D^2 log u| \leq c \cdot t_o^{-4},$$

where $t_o(x) = |x|$. Since $log u \to log u_o$ at infinity in $\mathbb{R}^3$, it is then standard that $log u$ has the expansion

$$log \frac{u}{u_o} = 1 - \frac{m}{t_o} + O(t_o^{-2}),$$

and hence, on $E \setminus K$,

$$u = u_o - \frac{m u_o}{t} + O(t^{-2}), \tag{7.32}$$

where $m \in \mathbb{R}$. Thus $\omega = -u$ has the expansion (0.18) in $E$. Further, since $\widetilde{g} = u^{-2} \cdot g$, and $u^{-2} \sim u_o^{-2}(1 + \frac{2m}{t})$, it follows that $E$ is asymptotically flat in the sense of (0.17), (after applying a suitable dilation of $\mathbb{R}^3$), i.e.

$$g_{ij} = (1 + \frac{2m}{t})\delta_{ij} + O(t^{-2}). \tag{7.33}$$

As before, elliptic regularity implies that $|\partial^k g_{ij}| = O(r^{-2-k})$.

The mass $m = m_E$ above will depend on the end $E$, and from (7.32) and (7.33), can be computed as

$$4\pi u_o \cdot m_E = lim_{s \to \infty} \int_{S_E(s)} < \nabla u, \nu > dA, \tag{7.34}$$

where $\nu$ is the outward unit normal. From (7.33), $m$ agrees with the mass of the end $E$ in the sense of general relativity, c.f. [Wd, Ch. 6.2].

Suppose that all ends $E_i$ of $(N, g)$ are asymptotically flat in the sense of (7.33). Observe that by applying the divergence theorem to the trace equation (4.7) and using (7.32), (7.34), one has

$$\frac{\alpha}{4} \int_N |z|^2 = 4\pi \sum_i u_i \cdot m_{E_i} > 0, \tag{7.35}$$

where $u_i$ is the limiting value of $u$ in the end $E_i$, as in (7.25). In particular, the total mass of the ends of $N$ is positive. In fact, since $(N, g)$ is complete and of non-negative scalar curvature, the positive mass theorem [SY1] implies in this situation that $m_E > 0$ for each end $E$, so that $u$ is increasing to its supremum at infinity in each end.

When all ends of $N$ are asymptotically flat, then by Proposition 7.1, $N$ topologically has the form

$$N = P \# (\#_1^q \mathbb{R}^3).$$

It is clear from the asymptotics (7.33) that each end of $N$ may be metrically capped off with a 3-ball to form a closed 3-manifold $(\bar{N}, \bar{g})$ in such a way that $\bar{g}$ has scalar curvature $s_{\bar{g}} \geq 0$, and $s_{\bar{g}}$ is not identically 0. Since $\bar{N}$ is diffeomorphic to $P$, it follows that $\sigma(P) > 0$; compare with Remark 6.3.

Next, we proceed with examination of the other possible ends. As we will see in the following, all other ends are collapsing at infinity, on the scale of the $L^2$ curvature radius.

**Case I(B)(Collapse).** Under the standing assumption (7.23), suppose that the end $E$ is collapsing at infinity, in the sense that

$$liminf_{t \to \infty} \frac{v_E(s)}{s^3} = 0. \tag{7.36}$$

The results above preceding Case I(A), in particular (7.25), still hold, and hence from (7.29), $\widetilde{r}_- \geq -ct^{-4}$, where $\widetilde{r}_-$ denotes the negative part of the Ricci curvature. Thus, $\widetilde{g}$ is of almost non-negative Ricci curvature, in a strong sense. From standard volume comparison theory for such manifolds, (7.36) holds with limsup in place of liminf, so that $E$ is collapsing everywhere, (on the scale of the curvature radius $\rho$).

The main point here is to prove that the curvature of $\widetilde{g}$ decays sufficiently fast, i.e.

$$\int^\infty s|\widetilde{K}|(s)ds < \infty, \tag{7.37}$$

where $|\widetilde{K}|(s)$ is the maximum of $|\widetilde{K}|$ on $\widetilde{S}(s) \subset N$; here all quantities are w.r.t the $\widetilde{g}$ metric (7.28). To prove (7.37), from (7.29), as in (7.30), observe that

$$|\widetilde{K}|(s) \sim |\nabla logu|^2(s) + O(s^{-4}). \tag{7.38}$$

Hence, modulo terms with finite integral,

$$s|\widetilde{K}|(s) \leq c \cdot s|\nabla logu|^2(s) \leq c \cdot s|\nabla logu|(s) \cdot (areaS_E(s))^{-1} \leq c \cdot (areaS_E(s))^{-1},$$

where the second inequality follows from (7.22) and (7.25) while the last inequality follows from (7.13). Thus, the estimate (7.37) follows from the assumption (7.23).

We now use the following Lemma.

**Lemma 7.9.** *Let $(E, \widetilde{g})$ be an end of $(N, \widetilde{g})$ satisfying (7.36) and (7.37). Then outside a sufficiently large compact set $K \subset E$, the manifold $(E \setminus K, \widetilde{g})$ is quasi-isometric to either $(\mathbb{R}^2 \times S^1) \setminus K'$ or $(\mathbb{R}^+ \times S^1 \times S^1) \setminus K'$, for some compact set $K'$.*

*Proof.* This is a relatively standard result in comparison geometry. The condition (7.37) implies that the exponential map into $E$, (w.r.t. $\widetilde{g}$), is a quasi-isometry in the universal cover. The condition (7.36) means that $(E, \widetilde{g})$ cannot be quasi-isometric to $\mathbb{R}^3$ near infinity, and hence it must be quasi-isometric to a quotient of $\mathbb{R}^3$, i.e. of the form stated in the Lemma, near infinity. We refer to [An5,Lemma 3.7], c.f. also [Ab], for further details.

∎

By (7.25), $\widetilde{g}$ is quasi-isometric to $g$ outside a compact set in $E$, so that Lemma 7.9 also holds w.r.t. $g$. In either case of Lemma 7.9, one has at most quadratic volume growth of balls, or at most linear growth for the area of geodesic spheres, i.e.

$$area_g S(s) \leq v_1 \cdot s, \tag{7.39}$$

for some $v_1 < \infty$. However, the estimate (7.39) contradicts the standing assumption (7.23). Thus, Case I(B) cannot occur.

It follows that among the finite number of ends $E_i$ of $N$, each end satisfying (7.23) is asymptotically flat. Next, we consider ends which may not satisfy (7.23).

**Case II.** Suppose

$$\int^{\infty} areaS_E(s)^{-1}ds = \infty. \tag{7.40}$$

As in Case (I), we first prove some general results on ends satisfying (7.40), and will then specialize to three further subcases II(A)-(C) according to the precise behavior of the tangent cones at infinity.

The condition (7.40) implies that the end $E$ is collapsing at infinity, i.e. (7.36) holds, at least for some sequence $s_i \to \infty$. The blow-downs near infinity, i.e. the metrics $g_s = s^{-2} \cdot g$ in $A(\kappa^{-1}s, \kappa s) \subset E, \kappa$ large, are of bounded curvature, by (7.13), and of bounded diameter, by (7.17). Further, the oscillation of $u$ on each component of $A(\kappa^{-1}s, \kappa s)$ is uniformly bounded, depending only on $\kappa$, again by (7.13).

First, we claim these statements imply that $E$ is collapsing everywhere at infinity, i.e.

$$areaS_E(s)/s^2 \to 0 \quad \text{as} \quad s \to \infty. \tag{7.41}$$

To see this, consider the conformally equivalent metric $\widetilde{g} = u^2 \cdot g$ from (7.28). By (7.29)-(7.30), $\widetilde{g}$ has almost non-negative Ricci curvature, i.e. as noted in Case I(B), (7.30) holds with $\widetilde{r}_-$ in place of $\widetilde{r}$. Standard comparison geometry then implies that (7.41) holds w.r.t. $\widetilde{g}$-geodesic spheres for all $s$, if it holds for some sequence $s_i \to \infty$. This, together with the control on $u$ from (7.13) implies the same for $(E, g)$.

Thus, outside a sufficiently large compact set, $E$ admits a sequence of F-structures along which the blow-downs at infinity, (as in the construction of tangent cones at infinity), are collapsing.

Next, we have a global bound on the $L^2$ norm of $\nabla logu$ on $E \setminus K$.

**Lemma 7.10.** *There is a compact set $K$ in $E$ such that*

$$\int_{E \setminus K} |\nabla logu|^2 dV < \infty. \tag{7.42}$$

*Proof.* To obtain this, muliply the trace equation (4.9) by $u^{-1}$ on $B(s) \setminus K = B_E(s) \setminus K \subset E$ and apply the divergence theorem to deduce

$$\int_{B(s) \setminus K} |\nabla logu|^2 + \int_{\partial(B(s))} <\nabla logu, \nu> = \frac{\alpha}{4} \int_{B(s) \setminus K} u^{-1}|r|^2 + C, \tag{7.43}$$

where $C$ is a bound for boundary integral over $B(s) \cap \partial K$. From Corollary 7.6 and (7.9), the right side of (7.43) is uniformly bounded, (independent of $s$). Set

$$F(s) = \int_{B(s) \setminus K} |\nabla logu|^2,$$

and let $\lambda(s) = areaS(s)$ be the area of $\partial B(s)$. From the Hölder inequality, it then follows from (7.43) that there is a constant $C < \infty$ such that

$$F(s) \leq \left(F'(s) \cdot \lambda(s)\right)^{1/2} + C \leq 2max(C, \left(F'(s) \cdot \lambda(s)\right)^{1/2}).$$

Separating variables, i.e. $F$ and $\lambda$, this differential inequality may then be integrated from $r$ to $\infty$. Using the assumption (7.40), i.e. $\int^{\infty} \lambda^{-1}(s)ds = \infty$, one easily obtains (7.42); further details here, (if needed), are in the proof of [An4, Prop.4.1]. ∎

Now again we break the discussion into cases, according to the structure of the tangent cones at infinity.

Thus let $(T_{\infty}, g_{\infty}, x_{\infty})$ be a tangent cone at infinity of $(E, g)$, $x_{\infty} = \lim x_i$, $g_{\infty} = \lim g_i$, where $g_i = t_i^{-2} \cdot g$, $t_i = t(x_i) \to \infty$, as $i \to \infty$. Since the metrics $g_i$ are collapsing at $x_i$, $(T_{\infty}, g_{\infty})$ is obtained by unwrapping the collapse by taking sufficiently large finite covers. By Remark 2.1, $(T_\infty, g_\infty)$ has either a free isometric $S^1$ or $S^1 \times S^1$ action, corresponding (basically) to a rank 1 or rank 2 collapse respectively, c.f. Remark 2.1. These cases are distinguished by the growth of the diameter of geodesic spheres $S_c(t_i)$, where $S_c$ denotes the component containing the base point $x_i$. Namely, if there exists $d_o > 0$ such that

$$diam S_c(t_i) \geq d_o \cdot t_i, \tag{7.44}$$

then $(T_\infty, g_\infty)$ is formed from a rank 1 collapse and has a free isometric $S^1$ action. On the other hand, if

$$diam S_c(t_i) << t_i, \tag{7.45}$$

then $(T_\infty, g_\infty)$ is formed from a rank 2 collapse and has a free isometric $S^1 \times S^1$ action.

**Case II(A).(Rank 2 Collapse).** Suppose

$$diam S_c(t_i) << t_i, \text{ as } i \to \infty. \tag{7.46}$$

We first classify the possible tangent cones based on sequences satisfying (7.46).

**Proposition 7.11.** *Any tangent cone $(T_\infty, g_\infty, x_\infty)$, based on $x_i$ satisfying (7.46) is a flat product of the form $\mathbb{R}^+ \times S^1 \times S^1$.*

*Proof.* It suffices to prove that the tangent cone is flat, since the diameter condition (7.46) then implies that $(T_\infty, g_\infty)$ has sublinear diameter growth, and hence is of the form $\mathbb{R}^+ \times S^1 \times S^1$.

The tangent cone $T_\infty$ is flat exactly when the curvature of $g$ decays faster than quadratically in the annuli $A_i = A_c(\kappa^{-1} t_i, \kappa t_i)$ about the base points $x_i$, as in (7.24). We will suppose this is not the case, and derive a contradiction.

Thus, suppose that there are points $p_i \to \infty$, within bounded $g_i$-distance to $x_i$, such that the opposite inequality to (7.13) holds, i.e. for some $\kappa > 0$,

$$|r|(p_i) \geq \frac{\kappa}{t^2(p_i)}. \tag{7.47}$$

We may assume that $\{p_i\}$ converges to a point $p$ in the tangent cone $(T_\infty, g_\infty, x_\infty)$, which is then non-flat. Further, as in the proof of Lemma 7.5, since $\alpha$ as defined in §4 scales as the square of the distance, so $\alpha_i = \alpha t_i^{-2} \to 0$, it follows that the limit $(T_\infty, g_\infty)$ is a non-flat solution of the static vacuum Einstein equations (4.10). Here the limit potential $\bar{u}$ is renormalized so that $\bar{u}(x)$ $= \lim u(x)/|u(x_i)|$, giving $\bar{u}(x_\infty) = 1$. Since $u$ is bounded away from 0 outside a compact set, this renormalization does not affect the preceding arguments.

As discussed above, the limit $(T_\infty, g_\infty)$ has a free isometric $S^1 \times S^1$-action, arising from the rank 2 collapse of $E$ based at $\{x_i\}$.

Now the only non-flat static vacuum solution with a free $S^1 \times S^1$ action is the *Kasner metric*

$$g_K = dr^2 + r^{2a} d\theta_1^2 + r^{2b} d\theta_2^2, \tag{7.48}$$

on $(\mathbb{R}^+ \times S^1) \times S^1$, c.f. [EK, Thm.2-3.12] or [An5, Ex.2.11]. Here $a$ and $b$ are given by

$$a = (d-1)/(d+d^{-1}-1), b = (d^{-1}-1)(d+d^{-1}-1), \tag{7.49}$$

where $d \in (0,1)$ is a free parameter, so that there is a 1-parameter family of such metrics. Note that $a < 0$, $b > 0$ with $|a|, |b| \in (0,1)$ and $a + b \in (0,1)$. The potential function $u$ of the static vacuum solution is given, (up to a multiplicative constant), by $u = r^c, c = 1 - (a+b)$. We note that there are also Kasner metrics with $d \in (-1,0)$, so that $a + b > 1$. In this case however, the potential tends to 0 at infinity, which violates (7.9).

Observe also that the Kasner metrics (7.48) interpolate between the two (extreme) flat metrics, in the sense that $d = 0$ gives $a = 0$, $b = 1$ and hence the flat product $\mathbb{R}^2 \times S^1$, while $d = 1$ gives $a = 0$, $b = 0$ and so the flat product $\mathbb{R}^+ \times S^1 \times S^1$.

On the other hand, the metric $g$, and hence the metrics $g_i$, are $\mathcal{Z}_s^2$ solutions on the domains $A_i$ converging to $T_\infty$. Thus, renormalizing by $\alpha$, we have

$$\nabla \mathcal{Z}^2 = L^*(\frac{u}{\alpha}), \quad \frac{1}{4}|r|^2 = \Delta(\frac{u}{\alpha}); \tag{7.50}$$

all quantities in (7.50) are w.r.t. $g_i$, and $\alpha = \alpha_i \to 0$. The left sides of (7.50) are bounded as $i \to \infty$, and hence the right sides of (7.50) are bounded. Of course $u/\alpha \to \infty$ everywhere in $A_i$. We may decompose $u/\alpha$ as a sum

$$\frac{u}{\alpha} = \nu + \eta,$$

where $\nu = \nu_i$ remains uniformly bounded as $i \to \infty$, and $\eta = \eta_i$ becomes unbounded, but approaches Ker $L^*$. In other words, $\nabla \mathcal{Z}^2$ is uniformly bounded on $(A_i, g_i)$ and in Im $L^*$, and hence converges, (in a subsequence), to a limit $\nabla \mathcal{Z}^2$ on $(T_\infty, g_\infty)$ which is still in Im $L^*$. It follows that the limit $(T_\infty, g_\infty, \nu)$ is a solution to the $\mathcal{Z}_s^2$ equations

$$\nabla \mathcal{Z}^2 = L^* \nu, \quad \Delta \nu = \frac{1}{4}|r|^2. \tag{7.51}$$

This shows that the limit Kasner metric (7.48) must also be a $\mathcal{Z}_s^2$ solution, (with $\alpha = 1$), as the Schwarzschild metric is, (c.f. (6.17)). The following result shows this to be impossible.

**Lemma 7.12.** *The Kasner metric (7.48) is not a $\mathcal{Z}_s^2$ solution, for any value of $d \in (0, 1)$.*

*Proof.* The proof is a direct computation, since the Kasner metric has the explicit form (7.48). We defer the proof to Appendix B, since the methods involved are unrelated to the main developments here.

∎

Lemma 7.12 gives a contradiction to the construction of $T_\infty$ based on $\{p_i\}$ satisfying (7.47). It follows that there is no divergent sequence $\{p_i\} \in (E, g)$ satisfying (7.47), (near $x_i$) and hence the tangent cone $T_\infty$ is flat, as claimed.

∎

As observed in the proof above, Proposition 7.11 implies that

$$|r|(x) << t^{-2}(x), \tag{7.52}$$

as $t(x) \to \infty$, on sequences $x = x_i$ satisfying (7.46).

Further, the flat tangent cone $T_\infty$ is also a static vacuum solution. Thus, the potential function, (renormalized at the base sequence $\{x_i\}$ as above), is asymptotic either to a constant, or to a non-constant affine function. This means that either

$$t|\nabla log u|(x) << 1, \quad \text{or} \quad t|\nabla log u| - 1|(x) << 1, \tag{7.53}$$

as $t = t_i \to \infty$.

Having identified the structure of the tangent cones at infinity in this situation, we are now in position to rule out the existence of ends having a rank 2 collapse along some divergent sequence in $E$.

**Proposition 7.13.** *No end $(E, g)$ in $(N, g)$ has a tangent cone at infinity of the form $\mathbb{R}^+ \times S^1 \times S^1$.*

*Proof.* Suppose $(T_\infty, g_\infty, x_\infty)$ is such a tangent cone at infinity, based on a divergent sequence $\{x_i\} \in E$.

To understand this situation, we linearize the $\mathcal{Z}_s^2$ equations at infinity. Let

$$\delta = \delta_i = (\int_{B_{x_i}(\frac{1}{2})} |r_{g_i}|^2 dV_{g_i})^{1/2}, \tag{7.54}$$

where all metric quantities in (7.54) are taken w.r.t. $g_i$ converging to $g_\infty$. (Hence we have passed to covers to unwrap collapse). Thus $\delta_i \to 0$ as $i \to \infty$. Further, we may assume that $\delta_i > 0$ for each $i$, since $(N, g)$ is not flat on any open set.

Consider the $\mathcal{Z}_s^2$ equations in the $g_i$ scale, i.e., (since the metric is scalar-flat),

$$\alpha \nabla \mathcal{Z}^2 - L^*(u) = 0, \tag{7.55}$$

$$\Delta u = \frac{1}{4} \alpha |r|^2.$$

Here again all metric quantities are w.r.t. $g_i$, and $\alpha = \alpha_i = t_i^{-2} \cdot \alpha \to 0$, where the factor $\alpha = \alpha_g$ is from the $\mathcal{Z}_c^2$ solution $(N, g)$. Consider the linearization of these equations w.r.t. $\delta$ at $\delta = 0$. From the regularity given by Theorem 4.2/Remark 4.3(i) and the definition (7.54), we have

$$lim_{\delta \to 0} \frac{|r|^2}{\delta} = 0,$$

on $T_\infty$. Here and in the following, $\delta \to 0$ means that $\delta_i \to 0$ as $i \to \infty$. We will also pass to suitable sequences of $\{i\}$ as needed to obtain the limits, without being explicit about this. Hence,

$$lim_{\delta \to 0} \frac{\alpha |r|^2}{\delta} = 0, \text{ and } lim_{\delta \to 0} \frac{\alpha (\nabla \mathcal{Z}^2)}{\delta} = 0, \tag{7.56}$$

where the second equality again follows from Theorem 4.2/Remark 4.3(i).

Next, consider the linearization of the metrics $g_i$ at the limit $g_\infty$. Thus, set

$$g_i = g_\infty + \delta_i h_i,$$

where this equation defines $h_i$. Since, by construction, c.f. (7.54), the linearization of the curvature is bounded, it follows that $\{h_i\}$ is also bounded, and converges to the limit $h$, as $i \to \infty$. More precisely, one may need to alter $g_i$ by suitable diffeomorphisms $\psi_i$ of $T_\infty$, with $\psi_i$ converging to the identity map, so that the forms $h_i$ are locally bounded and converge to a limit $h$.

Thus, at least formally, the equations (7.55) linearized at infinity, i.e. at $\delta = 0$, correspond to the linearization of the static vacuum equations $L^*u = 0$ at the flat metric. We need to examine the behavior of the potential function to make this rigorous.

To do this, let $\bar{u}_i = u/u(x_i)$, so that $\bar{u}_i(x_i) = 1$. As noted above, since $(T_\infty, g_\infty)$ is a flat static vacuum solution, the limit function $\bar{u}_\infty$ on $T_\infty$ must either be a constant, thus 1, or a non-constant affine function on $T_\infty$.

Now consider the linearization of $\bar{u}_i$ at $\delta = 0$. Thus, write

$$\bar{u}_i = 1 + \delta_i \nu_i, \tag{7.57}$$

so that formally $\nu = (d\bar{u}/d\delta)|_{\delta=0}$. Observe from the above that $\nu_i(x_i) = 0$.

Suppose for the moment that $\nu_i$ is locally uniformly bounded on $\mathbb{R}^+ \times S^1 \times S^1$. It follows that $\nu_i$ (sub)-converges to a limit function $\nu$ and, by (7.54)-(7.57), one obtains a solution to the linearized static vacuum equations at the flat pair $(g_\infty, \bar{u}_\infty)$, i.e.

$$\bar{u}_\infty \cdot r' = D^2\nu, \quad \Delta\nu = 0. \tag{7.58}$$

where $r'$ is defined as $lim_{\delta \to 0} r/\delta$.

The solution is non-trivial since, by construction, $r'$ is not identically 0. The linearization $\nu$ is a harmonic function on $\mathbb{R}^+ \times S^1 \times S^1$, invariant under the $S^1 \times S^1$-action, since $\nu_i$ has been lifted to the covers to unwrap collapse. Hence $\nu$ is harmonic on $\mathbb{R}^+$, and so $\nu$ is an affine function. This means that $D^2\nu = r' = 0$, a contradiction.

Now essentially the same argument as above holds even if $\{\nu_i\}$ is unbounded. Thus from (7.55)-(7.56), we have

$$\bar{u}_i r = D^2 \bar{u}_i + o(\delta_i), \quad \Delta(\bar{u}_i) = o(\delta_i), \tag{7.59}$$

on $(A(\kappa^{-1}, \kappa), g_i) \to (A(\kappa^{-1}, \kappa), g_\infty) \subset (T_\infty, g_\infty)$. Divide this equation by $\delta = \delta_i$, to obtain

$$\bar{u}_i(\frac{r}{\delta}) = D^2(\frac{\bar{u}_i}{\delta}) + o(1), \ \ \Delta(\frac{\bar{u}_i}{\delta}) = o(1). \tag{7.60}$$

We have $\bar{u}_i \to \bar{u}_\infty$ and $r/\delta \to r'$. Hence $D^2(\bar{u}_i/\delta)$ is bounded on $(A(\kappa^{-1}, \kappa), g_i)$. For the same reasons as discussed in (7.50)-(7.51), we may write

$$\frac{\bar{u}_i}{\delta} = \bar{\nu}_i + \eta_i,$$

where $\bar{\nu}_i$ is uniformly bounded (locally) and $\eta_i$ diverges to infinity but approaches $\mathrm{Ker} D^2$. Thus, there is a limit function $\bar{\nu}$ satisfying

$$\bar{u}_\infty \cdot r' = D^2\bar{\nu}, \ \ \Delta\bar{\nu} = 0. \tag{7.61}$$

This gives the same contradiction as before.

Of course the fact that $\mathrm{Ker} D^2 \neq \{0\}$ is related to the possibility that the metrics $g_i$ may have to be altered by (affine) diffeomorphisms $\psi_i$ as above to obtain the linearization $h$ of $g_\infty$. ∎

Propositions 7.11 and 7.13 imply that no end $E$ has a divergent sequence on which there is a rank 2 collapse, and so, by (7.44)-(7.45), there is a constant $d_o > 0$, such that, for all $s$ large and $S_c(s) \subset E$,

$$diam S_c(s) \geq d_o \cdot s. \tag{7.62}$$

All tangent cones at infinity of $E$ are thus formed by unwrapping a rank 1 collapse. In analogy to Proposition 7.11, we next prove the following result in this situation.

**Proposition 7.14.** *Any tangent cone at infinity $T_\infty$ of $(E, g)$ satisfying (7.62) and the collapse assumption (7.40) is a flat product of the form $C \times S^1$, where $C$ is a flat Euclidean cone, possibly with an isolated singularity at $\{0\}$.*

*Proof.* As in the proof of Proposition 7.11, we claim that it suffices to prove that $T_\infty$ is flat. For in this case, the diameter assumption (7.62) and the collapse assumption (7.40), (c.f. also the remarks following (7.40)), imply that $T_\infty$ is a flat product $C \times S^1$, where $C$ is necessarily a flat surface of linear diameter growth. From the construction of the tangent cones preceding Corollary 7.6, $\mathrm{diam} S_C(s) \to 0$ as $s \to 0$, where $S_C(s) = t_\infty^{-1}(s) \cap C$. Thus the metric boundary $\partial C$ of $(C, g_\infty)$ is given by one point $\{0\}$, so that the metric completion $\bar{C}$ is given by $\bar{C} = C \cup \{0\}$. This implies that $C$ is a complete Euclidean cone, possibly with isolated singularity at $\{0\}$. Thus we need to prove $T_\infty$ is flat.

Now for the same reasons as in Proposition 7.11, (namely from Lemma 7.5), any tangent cone $(T_\infty, g_\infty)$ is a static vacuum solution. Since $T_\infty$ is formed by unwrapping a collapse, it follows that $(T_\infty, g_\infty, \bar{u}_\infty)$ has a free isometric $S^1$ action. By [An5,Prop.2.2], $(T_\infty, g_\infty, \bar{u}_\infty)$ is a *Weyl solution*, i.e.

$$g_\infty = g_V + f^2 d\theta^2,$$

where $g_V$ is a metric on a surface $V$ and $f$ denotes the length of the $S^1$ fiber. The discussion above on the structure of $T_\infty$ implies here that $(V, g_V)$, (or more precisely its completion), is a complete Riemannian surface, possibly with an isolated singularity at $\{0\}$, of linear diameter growth. Observe from the construction that both the fiber length $f$ and the potential function $u = \bar{u}_\infty$ are positive on $V$.

From the theory of Weyl solutions, c.f. [Wd, (7.1.21)], (or [An5,(2.10)]), it is well-known that the product $f \cdot u$ is a harmonic function on $V$. However, since the metric completion is given by $\bar{V} = V \cup \{0\}$, it is standard that $V$ admits no globally defined positive harmonic functions except constants. When $f \cdot u$ is constant, it is again standard that the Weyl solution is flat, c.f. [Wd, §7.1] or [An5, (2.9)ff] for example.

Proposition 7.14 shows that tangent cones at infinity $T_\infty$ of $E$, i.e. those based on *any* divergent sequence $\{p_i\}$, are always of the form $C \times S^1$. In particular, as in (7.52), we have

$$|r|(x) << t^{-2}(x), \tag{7.63}$$

uniformly as $t(x) \to \infty$ in $E$. As noted following (7.24), it follows that $E$ is diffeomorphic to $(\mathbb{R}^2 \setminus B) \times S^1$, $B$ a ball in $\mathbb{R}^2$, and the geodesic spheres $S_E(s) \subset E$ are connected, for $s$ large.

Referring to the discussion preceding (7.57), we note that in this case, the limit renormalized potential $\bar{u}_\infty$ always satisfies

$$\bar{u}_\infty \equiv 1, \tag{7.64}$$

regardless of the choice of base point sequence. For otherwise $\bar{u}_\infty$ must be a non-constant affine function, necessarily $S^1$ invariant, on $T_\infty$, (for some base point sequence). Since $C$ is a Euclidean cone, it would follow that $\bar{u}_\infty$ achieves both positive and negative values on $T_\infty \setminus \{0 \times S^1\}$. This is of course impossible, (since $u > 0$ everywhere in $E$).

In particular, the first alternative in (7.53) holds and thus $u$ has sublinear growth at infinity in $E$. (Note that (7.64) does *not* imply that $u$ is bounded on $E$). It also follows from (7.64), together with the global diameter bound (7.17) and the maximum principle for the trace equation (4.7), that $u$ increases to its supremum, possibly infinite, uniformly at infinity in $E$.

Next, as in Proposition 7.13, consider the linearization of the $\mathcal{Z}_s^2$ equations at the flat metric $(T_\infty, g_\infty)$, with $\bar{u}_\infty \equiv 1$, again for arbitrary base point sequences. The arguments in (7.57)-(7.61) give a non-trivial solution of the linearized static vacuum equations, i.e.

$$r' = D^2\nu, \quad \Delta\nu = 0.$$

Since, by construction, $r'$ is non-zero, (i.e. not identically 0), it follows that $\nu$ is a globally defined non-constant harmonic function on $T_\infty = (C \setminus \{0\}) \times S^1$, invariant under the $S^1$ action on the second factor, obtained by unwrapping the collapse. The metric on $T_\infty$ is a product, so that $\nu$ is harmonic on $C \setminus \{0\}$. Further, the fact that $u$ is everywhere increasing to sup$u$ at infinity in $E$ implies that $\nu$ extends to a globally defined subharmonic function on $\bar{C} = C \cup \{0\}$, (because the approximations $\nu_i$ in (7.57) tend uniformly to $-\infty$ near $\{0\} \subset \bar{C}$, as $i \to \infty$).

Now the only such functions, up to addition of affine functions, are

$$\nu = a \cdot log t_C, \tag{7.65}$$

for some constant $a$, where $t_C = dist_C(\cdot, 0)$. Although the exact value of $a$ is not determined, (i.e. it may depend on the base point sequence), note that it is necessarily bounded away from 0 and $\infty$, since $r'$ has size on the order of 1, by construction.

We claim that in this situation the approximations $\nu_i$ defined by (7.57) are locally bounded as $i \to \infty$, and in fact $\nu = \lim\nu_i$ for $\nu$ as in (7.65). To see this, it follows from (7.65) and the linearized static vacuum equations that on *any* tangent cone $T_\infty$, $|r'| \sim 1/t^2 log t$, where $t = t_\infty$ on $T_\infty$. From the construction of $r'$, i.e. from (7.54), it follows that

$$|r| \sim \frac{1}{t^2 log t}$$

on $(E, g)$, i.e. we have a more accurate estimate for the decay of $|r|$ than (7.63). This in turn implies from the $\mathcal{Z}_s^2$ equations that $|D^2u| \sim t^{-2}(log t)^{-1}$, and so, by integration, $u$ has logarithmic growth at infinity in $E$. The definition of the approximations $\nu_i$ in (7.57) then proves the claim.

In sum, it follows from the discussion above that, to leading order,

$$u \sim log t, \quad |\nabla u| \sim \frac{1}{t}, \quad |r| \sim \frac{1}{t^2 log t}, \tag{7.66}$$

as $t \to \infty$. In particular, we see that the potential function $u$ is unbounded on $E$.

We now need, finally, to divide this situation into two further cases, according to the behavior of the area growth of the geodesic spheres.

**Case II(B). (Rank 1 Collapse, sublinear area growth).** Suppose, given a divergent sequence $\{x_i\} \in E$, we have

$$diam S_E(s) \geq d_o \cdot s, \text{ and } area S_E(t_i) << t_i, \tag{7.67}$$

for all $s$ and $t_i = t(x_i)$.

In this case, we have the following analogue of Proposition 7.13.

**Proposition 7.15.** *No end $E \subset (N, g)$ has a divergent sequence $\{x_i\}$ satisfying (7.67).*

*Proof.* Consider first the quantity

$$\mu(t) = \int_{S_E(t)} < \nabla u, \nabla t > . \tag{7.68}$$

This corresponds to the mass of the positive measure $(\Delta u) dV$ on $B_E(t) \subset (E, g, u)$, for $t$ large, (compare with (7.34)). We claim that there is a constant $\mu(\infty)$ such that as $t \to \infty$,

$$\mu(t) \to \mu(\infty) > 0. \tag{7.69}$$

To see this, from the decay estimate (7.13) and the fact that the area growth of geodesic spheres in $E$ at most quadratic in $t$, (c.f. (7.18)), it follows from the trace equation (7.55) that

$$|\mu(t) - \mu(\infty)| \leq c \int_{A_E(t, \infty)} |r|^2 \to 0, \text{ as } t \to \infty. \tag{7.70}$$

This proves the first statement in (7.69). If $\mu(\infty) = 0$, then (7.70) implies $\mu(t) \sim t \int_{S(t)} |r|^2$, so that again by (7.13), $(vol S_E(t))^{-1} | \int_{S_E(t)} < \nabla u, \nabla t > | \leq c \cdot t^{-3}$. Since to leading order $u$ is of the form (7.66), this is impossible. Hence $\mu(\infty) > 0$.

However, (7.68)-(7.69) and the decay estimate (7.66) for $|\nabla u|$ imply that, as $i \to \infty$,

$$area(S_E(t_i)) \geq c \mu(\infty) t_i.$$

This contradicts the collapse assumption (7.67).

∎

The end $E$ satisfies (7.62) and hence Proposition 7.15 implies that the only remaining possibility for the structure of a collapsing end $E \subset N$ is the following.

**Case II(C).(Rank 1 Collapse, linear area growth).** Suppose there exist constants $d_o > 0$ and $\nu_o > 0$ such that

$$diam S_E(s) \geq d_o \cdot s, \text{ and } area S_E(s) \geq \nu_o \cdot s, \tag{7.71}$$

as $s \to \infty$.

Now Lemma 7.10, i.e. the estimate (7.42), together with the area estimate (7.71) and the mean-value inequality (7.21) from Lemma 7.7 imply, (since $u \geq u_o > 0$ in $E$),

$$\int^{\infty} t |\nabla log u|^2 (t) dt < \infty, \tag{7.72}$$

where $|\nabla log u|^2(t)$ is the maximum of $|\nabla log u|^2$ on $S_E(t)$. Hence, as in (7.38), (7.37) holds here.

It follows from Lemma 7.9 and the diameter/area growth estimate (7.71) again that the the conformally equivalent metric $(E, \widetilde{g})$, $\widetilde{g} = u^2 \cdot g$, is quasi-isometric to a flat product $\mathbb{R}^2 \times S^1$ near infinity. In particular, outside a large compact set, $E$ is topologically $(\mathbb{R}^+ \times S^1) \times S^1$. Let $F : E \to (\mathbb{R}^+ \times S^1) \times S^1$ be such a quasi-isometric diffeomorphism and let $S_2^1$, (resp. $S_1^1$), denote the circles $S^1$ in $E$ corresponding to the second $S^1$ factor, (resp. the first $S^1$ factor), in $(\mathbb{R}^+ \times S^1) \times S^1$ via the map $F$. Hence $S_2^1$ has uniformly bounded length w.r.t. $(E, \widetilde{g})$, i.e.

$$L_{\widetilde{g}}(S^1) = u \cdot f \leq K, \tag{7.73}$$

where $f$ is the length of the $S_2^1$ factor w.r.t. $(E, g)$.

The length of the first $S^1$ factor $S_1^1$ in $(E, g)$ grows at most linearly in $t$, since the diameter of geodesic spheres in $E$ grows at most linearly by (7.17). By the diameter and area assumption (7.71), this means that the length of the first $S^1$ factor $S_1^1$ in $E$ does grow linearly in $t$ and further $f$, the length of the $S_2^1$ factor, must also be bounded away from 0. Hence (7.73) implies that $u$ is bounded above on $(E, g)$ and $f$ is bounded below away from 0.

However, (7.66) implies that $u$ is unbounded. This contradiction implies that no end $E$ of $(N, g)$ satisfies the assumptions of Case II(C).

Since all other cases have been ruled out, it follows that the limit $\mathcal{Z}_c^2$ solution $N$ has only ends of Case I(A), i.e. asymptotically flat ends. Together with the previous remarks in Case I(A), this completes the proof of Theorem C.

∎

For the purposes of completeness, we note that if $(E, g)$ is any asymptotically flat end of a $\mathcal{Z}_c^2$ or static vacuum solution, then its linearization at infinity, as in Case II(A) above, (c.f. (7.57)-(7.58)), is determined by a potential $\nu$ of the form

$$\nu = -\frac{m}{t},$$

where $m$ is the mass of the end $E$. In other words, $\nu$ is a multiple of the fundamental solution of the Laplacian on $\mathbb{R}^3$.

**Remark 7.16.** The conclusions of Theorem C agree with the models for the degeneration of $\{g_\varepsilon\}$ in §6, compare in particular the difference between (6.27) and (6.28).

We conclude with a brief discussion regarding the validity of Theorem C when the assumption that the levels of $\omega$ are compact is removed. Observe that this assumption was used only in Lemmas 7.4 and 7.8. In place of Lemma 7.4, one may use the following:

**Proposition 7.17.** Let $(N, g', \omega)$ a complete $\mathcal{Z}_c^2$ or $\mathcal{Z}_s^2$ solution. Suppose there is a compact set $K \subset N$ such that, on $N \setminus K$,

$$\omega \leq \omega_o < 0. \tag{7.74}$$

Then there exists a constant $\delta_o > 0$, depending only on $\alpha / |\omega_o|$ such that if

$$|r| \leq \delta_o, \tag{7.75}$$

on $N \setminus K'$, for some compact set $K' \subset N$, then (7.10) holds, i.e.

$$limsup_{t \to \infty} |r| = 0. \tag{7.76}$$

*Proof.* We first claim that if $(\bar{N}, \bar{g}, \bar{\omega})$ is any complete $\mathcal{Z}_s^2$ solution with $|r| \leq \Lambda, \bar{\omega} \leq \bar{\omega}_o < 0$ everywhere and $|r|(x) = 1$, at some base point $x \in \bar{N}$, then there is a constant $\alpha_o > 0$, (independent of the $\mathcal{Z}_s^2$ solution), such that $\alpha / |\bar{\omega}_o| \geq \alpha_o$. The proof of this is essentially the same as that of Lemma 7.2 or 7.5.

Thus, suppose this were not the case. Then there exists a sequence of complete $\mathcal{Z}_s^2$ solutions $(N_i, g_i, \omega_i)$, with normalized coefficients $\alpha_i = \alpha_i / \bar{\omega}_{o,i} \to 0$, with uniformly bounded curvature and $|r|(x_i) = 1$ at some base point $x_i \in N_i$. We may pass to a limit $(N_\infty, g_\infty, x_\infty, \omega_\infty)$ of a subsequence, (unwrapping to the universal cover in case of collapse) of the pointed sequence $(N_i, g_i, \omega_i, x_i)$. As in the proof of Lemma 7.2, renormalize $\omega_i$ w.r.t. the base point $x_i$ so that the limit potential satisfies $\omega_\infty(x_\infty) = -1$.

The limit is then necessarily a complete solution of the static vacuum equations, with limit potential $\omega_\infty \leq 0$ everywhere. Since the convergence to the limit is smooth, (by Theorem 4.2, or [An4, Thm.3.6]), $|r|(x_\infty) = 1$, so the limit is non-flat. Also since $\omega_\infty(x_\infty) = -1$, the limit is not a super-trivial static vacuum solution. This situation contradicts Theorem 5.1(III), which thus proves the claim.

The assumption that $|r| \leq \Lambda$ above is not necessary, since if $|r|$ becomes unbounded somewhere, one may apply the same argument as in the proof of Lemma 4.1, (or Lemma 7.5), to obtain new sequences $(N_i, g_i, x_i, \omega_i)$ satisfying the same properties as before.

Now suppose that (7.76) were not true under the bounds (7.74)-(7.75). Then there exists a constant $\delta_1 > 0$ and a divergent sequence $x_i$ in $N$ with $|r_i|(x_i) \geq \delta_1$. As above, the pointed manifolds $(N, g, \omega, x_i)$ then give rise to a limit complete $\mathcal{Z}_s^2$ solution, $(N_\infty, g_\infty, x_\infty, \omega_\infty)$, with $|r|(x_\infty) \geq \delta_1 > 0$, with $|r| \leq \delta_o$ everywhere, (by continuity of (7.75), and with $\alpha$ fixed, i.e. $\alpha$ for the $\mathcal{Z}_s^2$ solution $N_\infty$ is the same as the $\alpha$ for the $\mathcal{Z}_s^2$ solution $N$. As above, renormalize $\omega$ and $\alpha$ by dividing by $|\omega_o|$, so that $\omega_\infty$ is bounded away from 0. Now rescale $(N_\infty, g_\infty)$ so $|r|(x_\infty) = 1$. This changes the coefficient $\alpha/|\omega_o|$ into $(\alpha/|\omega_o|) \cdot \delta_1^2 \leq (\alpha/|\omega_o|) \cdot \delta_o^2$. Hence if $\delta_o$ is sufficiently small, this contradicts the result above.

∎

Similarly, the assumption that the levels of $\omega$ are compact can be removed from Lemma 7.8. As before, the proof is based on the estimate (7.20), together with (7.9); we refer to [An5, Lem. 3.6, Rmk. 3.9] for details.

In particular, it follows that Theorem C is valid with the hypotheses (7.74)-(7.75) in place of the assumption that the levels of $\omega$ are compact.

## Appendix A.

It may seem that the analysis of the structure and behavior of the metrics $g_\varepsilon$ bears little relation with the work in the previous paper [AnI], where the discussion was concerned only with the behavior of Yamabe metrics. At a deeper level however this is not really the case. For instance, a fair amount of the work in this paper can be done on the space of Yamabe metrics. We make several comments on this below.

Let $\mathcal{C}$ denote the space of Yamabe metrics on $M$ and consider the functional

$$J_\varepsilon = \left( \varepsilon v^{1/3} \int |z|^2 - v^{2/3} \cdot s \right)|_\mathcal{C} : \mathcal{C} \to \mathbb{R}. \tag{A.1}$$

This functional of course closely resembles the functional $I_\varepsilon$,

$$I_\varepsilon = \varepsilon v^{1/3} \int |z|^2 + \left( v^{1/3} \int s^2 \right)^{1/2} \tag{A.2}$$

from (3.6). In fact both functionals were treated together in [An3,§8], (with the insignificant difference that the full curvature $R$ was used in place of the trace-free Ricci curvature $z$ here and $\mathcal{S}^2$ had no square root). Minimizers $\kappa_\varepsilon$ of $I_\varepsilon$ and $\gamma_\varepsilon$ of $J_\varepsilon$ were shown to have practically identical properties. In particular, Theorem 3.9 holds exactly with $(I_\varepsilon, k_\varepsilon)$ and $(J_\varepsilon, \gamma_\varepsilon)$ in place of $(I_\varepsilon^-, g_\varepsilon)$; in fact in this case the metrics are $C^\infty$ smooth and there is no junction $\Sigma$.

The associated Euler-Lagrange equations are also formally almost identical to those of $I_\varepsilon^-$. Thus, as noted in §3.1(III), the Euler-Lagrange equations for $I_\varepsilon^-$ are identical to those of $I_\varepsilon$, with $\tau$ replaced by $s/\sigma$, with $\sigma$ given by the $L^2$ norm of the scalar curvature. Similarly, either directly or from a slight modification of [An3,(8.7)], (c.f. also [An3,(8.52)]), the Euler-Lagrange equations for $J_\varepsilon$ at the unit volume minimizer $\gamma = \gamma_\varepsilon$ on $\Omega_\varepsilon = \Omega_\varepsilon(\gamma_\varepsilon)$ may be written as

$$\varepsilon \nabla \mathcal{Z}^2 + z = L^* h - \frac{\varepsilon}{6} \mathcal{Z}^2(\gamma) \cdot \gamma, \tag{A.3}$$

$$2\Delta h + sh = \frac{\varepsilon}{2}|z|^2 - d_\varepsilon, \tag{A.4}$$

where $h$ is a function of 0 mean value and $d_\varepsilon = (\varepsilon/2)\mathcal{Z}^2(\gamma)$. Referring to [AnI, §2], $L^*h \in N\mathcal{C}$, the normal space to $\mathcal{C}$ in $T\mathbb{M}$ and these equations reflect that at a critical point of $J_\varepsilon, \nabla J_\varepsilon \in N\mathcal{C}$. Recall that $r = L^*(-1)$ so that $z = L^*(-1) - \frac{s}{3} \cdot \gamma$. Hence, (A.3) may be rewritten as

$$\varepsilon\nabla\mathcal{Z}^2 = L^*(1 + h) + \psi \cdot g, \tag{A.5}$$

where $\psi = (\frac{s}{3} - \frac{\varepsilon}{6}\mathcal{Z}^2)\gamma$.

On the other hand, referring to (3.9) and (3.14), the Euler-Lagrange equation for $I_\varepsilon{}^-$ may be written as

$$\varepsilon\nabla\mathcal{Z}^2 = L^*(-\tau) - \chi \cdot g, \tag{A.6}$$

$\chi = (\frac{1}{4}s\tau + \frac{1}{12}\mathcal{S}^2 + \frac{\varepsilon}{6}\mathcal{Z}^2)$. Thus, the Euler-Lagrange equations (almost) agree formally under the correspondence

$$-\tau \leftrightarrow (1 + h). \tag{A.7}$$

As $\varepsilon \to 0$, the family $\{\gamma_\varepsilon\}$ gives a preferred maximizing sequence of Yamabe metrics on a maximal domain $\Omega_\varepsilon(\gamma_\varepsilon)$. Most of the results in §3 hold for $(J_\varepsilon, \gamma_\varepsilon)$ and $(I_\varepsilon, \kappa_\varepsilon)$, where however all *integral* estimates, as in (3.38) for example, must be taken over the set $\{\tau \le 0\}$, (under the correspondence (A.7)), in place of all $\Omega_\varepsilon$. Further, certainly in one respect the functional $J_\varepsilon$ is preferable to $I_\varepsilon{}^-$ or $I_\varepsilon$, in that it can be used to study manifolds $M$ with $\sigma(M) > 0$.

However, while certain analogues of the work in §4 and §5 can be shown to hold for $(J_\varepsilon, \gamma_\varepsilon)$ and $(I_\varepsilon, \kappa_\varepsilon)$, (by using the method of [AnI,Thm. 3.10]), one cannot obtain the full strength of these results for these functionals. This is closely related to the horizon problem discussed in §1(ii). One concrete example of the difference is the change of sign of the potential $\tau$ or $1 + h$ for $I_\varepsilon$ or $J_\varepsilon$, and the lack of strong convergence when the potential has the wrong sign, as discussed in Remark 4.3(ii).

Even so, beyond the close correspondence of these functionals, there are more general relations between the metrics $g_\varepsilon$ and Yamabe metrics. Thus, the main structural identity used to study Yamabe metrics in [AnI] is the $z$-splitting

$$z = L^*f + \xi, \tag{A.8}$$

with $\xi \in \mathrm{Ker}L$, or the associated $g$-splitting

$$-\xi = L^*u + \frac{s}{3} \cdot g, \tag{A.9}$$

where $u = 1 + f$.

Now it is clear that (A.9) also bears a strong formal relation with (A.5) or (A.6) under the correspondence

$$u \leftrightarrow (1 + h) \leftrightarrow -\tau, \tag{A.10}$$

$$\xi \leftrightarrow -\varepsilon\nabla\mathcal{Z}^2.$$

Of course while (A.9) gives an $L^2$ bound on $\xi$, which leads to the basic results of [AnI], there is no such apriori bound on the $L^2$ norm of $\varepsilon\nabla\mathcal{Z}^2$ as $\varepsilon \to 0$; ($\varepsilon\nabla\mathcal{Z}^2$ is not in Ker $L$). However, the behavior of the potential function $u$ plays a central role in the analysis of the degeneration of Yamabe metrics in [AnI], just as the behavior of $\tau$ does here in the degeneration of $(\Omega_\varepsilon, g_\varepsilon)$.

These similarities exist because, although we are treating different functionals, they are all (perturbations of) scalar curvature functionals and hence their derivatives or gradients are associated with the linearized operators $L$ and $L^*$.

## Appendix B.

In this appendix, we prove Lemma 7.12, i.e. that the Kasner metric

$$g = dt^2 + t^{2a}d\theta_1^2 + t^{2b}d\theta_2^2, \tag{B.1}$$

representing a non-flat static vacuum solution, is not also a $\mathcal{Z}_s^2$ solution.

Recall that $a = (d-1)/(d+d^{-1}-1), b = (d^{-1}-1)(d+d^{-1}-1), d \in (0,1)$, so that $a < 0$, $b > 0$ and $|a| = db$. Setting w.l.o.g. $\alpha = 1$, the $\mathcal{Z}_s^2$ equations have the form

$$\nabla \mathcal{Z}^2 = L^* u, \quad \Delta u = \frac{1}{4}|z|^2. \tag{B.2}$$

We observe first that we may assume that $u$ is a function of $t$ only. For since the metric is $S^1 \times S^1$ invariant, if $u$ is any solution of (B.2), then $\bar{u}$, the integral of $u$ over $S^1 \times S^1$, will also satisfy (B.2) with $\bar{u} = \bar{u}(t)$. The trace equation in (B.2) of course implies that $\bar{u}$ is not identically 0.

Now it is a straightforward but lengthy computation that (B.2) admits no solutions with $u = u(t)$. This is carried out in the following.

First, let $e_1 = \nabla t$ and let $e_2, e_3$ be unit vectors tangent to the first and second $S^1$ factors respectively. Then $\{e_i\}$ is an orthonormal frame for $(N, g)$, diagonalizing the Ricci curvature $r$, with eigenvalues $r_i$. A simple computation gives

$$r_1 = -\lambda_1/t^2, \quad r_2 = -\lambda_2/t^2, \quad r_3 = -\lambda_3/t^2. \tag{B.3}$$

where

$$\lambda_1 = ab < 0, \quad \lambda_2 = a(a-1) + ab < 0, \quad \lambda_3 = b(b-1) + ab > 0. \tag{B.4}$$

Note that $\sum \lambda_i = 0$, since the metric is scalar flat.

Next we solve the trace equation in (B.2). We have $\Delta u = u'\Delta t + u''$ and $\Delta t = (a+b)/t$. Let $\lambda^2 = \sum \lambda_i^2 > 0$, so that $|r|^2 = |z|^2 = \lambda^2/t^4$. Then the differential equation arising from (B.2) may be easily integrated and gives

$$u = \frac{\lambda^2}{8(3-(a+b))}t^{-2} + c_1 t^{1-(a+b)} + c_2, \tag{B.5}$$

where $c_1, c_2$ are arbitrary constants. (Note that $t^{1-(a+b)}$ is harmonic, giving the potential of the static vacuum solution (B.1)).

We next take the inner product of (B.2) with $r$, (since this is the only natural thing to do). Using the formula (3.8) for $\nabla \mathcal{Z}^2$ from [B, Ch. 4H], a straightforward computation then gives

$$-\frac{\lambda^2}{t^6}(10 - 2(a+b)) + |Dr|^2 - 6\lambda_1\lambda_2\lambda_3 = -\frac{u\lambda^2}{t^4} + \frac{\lambda_1 u''}{t^2} + (\frac{a\lambda_2}{t^3} + \frac{b\lambda_3}{t^3})u'. \tag{B.6}$$

Since the left side of (B.6) decays at $t^{-6}$, so must the right side, which implies that $c_1 = c_2 = 0$ in (B.5).

Computing $|Dr|^2$ and the right side of (B.6) leads to the equation

$$-\lambda^2(6 - 2(a+b)) + 2(a^2(\lambda_2 - \lambda_1)^2 + b^2(\lambda_3 - \lambda_1)^2) - 6\lambda_1\lambda_2\lambda_3 = \tag{B.7}$$

$$= \frac{\lambda^2}{8(3-(a+b))}(-\lambda^2 + 6\lambda_1 - 2a\lambda_2 - 2b\lambda_3).$$

Now using the fact that $|a| = db$ and the relations (B.4), one may compute that

$$2(a^2(\lambda_2 - \lambda_1)^2 + b^2(\lambda_3 - \lambda_1)^2) = (2d^2 + \frac{3}{d^2-d+1})\lambda^2 b^2. \tag{B.8}$$

Observe that $-6\lambda_1\lambda_2\lambda_3 \leq 0$. It follows from this, (B.7) and (B.8) that

$$-6 + 2b(1-d) + b^2(2d^2 + \frac{3}{d^2 - d + 1}) \geq \frac{1}{8(3 - (a+b))}(-\lambda^2 + 6\lambda_1 - 2a\lambda_2 - 2b\lambda_3),$$
(B.9)

where we have divided by $\lambda^2$. A straightforward but laborious computation now shows that (B.9) is impossible; essentially, the right side of (B.9) is small, while the left side is negative, and not small, due to the dominant $-6$ term. We leave further verification to the reader.

## Appendix C.

We present here the proof of Theorem 5.3, i.e. the only complete $\mathcal{Z}_c^2$ solution $(N, g, \omega)$, with a free isometric $S^1$ action is the flat metric, with constant potential $\omega$. Following this, we complete the proof of Theorem 3.12, c.f. Proposition C.3.

Theorem 5.3 is obvious when $N$ is compact, by using the maximum principle applied to the trace equation (4.7). Thus we assume that $(N, g)$ is complete and non-compact.

Let $V$ denote the orbit space of the free $S^1$ action, $V = N/S^1$, so that $V$ is a smooth surface. If $V$ is not simply connected, then $N$ has an essential torus. Hence the result of Schoen-Yau [SY3, Thm. 4] or Gromov-Lawson [GL, Thm. 8.4], (used in Prop. 7.1), shows that $(N, g)$ is flat. Thus, we may assume that $V$ is topologically $\mathbb{R}^2$.

The metric $g$ induces a complete Riemannian metric $g_V$ on $V$ for which the projection $\pi : N \to V$ is a Riemannian submersion. Let $f : V \to \mathbb{R}$ denote the length of the orbits, i.e. $f(x)$ is the length of the $S^1$ fiber through $x$. Standard submersion formulas, c.f. [B, 9.37], imply that that the scalar curvature $s_V$ of $g_V$, equal to twice the Gauss curvature, is given by

$$s_V = s + |A|^2 + 2|H|^2 + 2\Delta log f,$$
(C.1)

where $A$ is the obstruction to integrability of the horizontal distribution and $H = \nabla log f$ is the geodesic curvature of the fibers $S^1$.

Fix any base point $x_o \in V$ and let $D(s)$ denote the geodesic $s$-ball about $x_o$ in $(V, g_V)$ and $C(s)$ the corresponding geodesic sphere, (i.e. circle). Let $v(s) = \text{area} D_{x_o}(s)$, so that $\lambda(s) \equiv v'(s)$ is the length of $C(s)$.

We quote the following general result from [An4, Prop. 4.1], which uses only the fact that $(N, g)$ has non-negative scalar curvature $s_g \geq 0$. (The proof follows from suitable applications of the Gauss-Bonnet formula on domains in $(V, g_V)$).

**Proposition C.1.** *In the notation above, there exists a constant $c < \infty$ such that*

$$v(s) \leq c \cdot s^2, \lambda(s) \leq c \cdot s,$$
(C.2)

*and*

$$\int_V |\nabla log f|^2 dA_V \leq c, \int_V |A|^2 dA_V \leq c.$$
(C.3)

Using this result, together with the methods used in §7, we will prove that $(N, g, \omega)$ satisfies the hypotheses of Theorem C. Let $t(x) = \text{dist}(x, y)$, for some fixed base point $y \in N$.

**Proposition C.2.** *Let $(N, g, \omega)$ be a complete $\mathcal{Z}_c^2$ solution, with a free isometric $S^1$ action. Then there is a constant $\rho_o > 0$ such that on $(N, g)$*

$$\rho \geq \rho_o \cdot t.$$
(C.4)

*Further, $\omega$ is a proper exhaustion function on $N$ with*

$$\omega(x) \to inf \, \omega < 0,$$
(C.5)

*as $t(x) \to \infty$ in $N$.*

*Proof.* We begin with the proof of (C.4); the estimate (C.5) is then a straightforward consequence. The proof of (C.4) is formally analogous to the proof of Lemma 7.5 and is proved by contradiction.

Thus, let $\{x_i\}$ be any sequence in $(N, g)$ with $t_i = t(x_i) \to \infty$, such that $\rho(x_i) << t_i$. As in the proof of Lemma 7.5, we may assume w.l.o.g. that $x_i$ almost realizes the minimal value of the ratio $\rho(x)/\text{dist}(x, \partial A(\frac{1}{2}t_i, 2t_i))$, for $x \in A(\frac{1}{2}t_i, 2t_i), t_i \to \infty$. By Proposition 4.10, we may perturb $x_i$ slightly if necessary so that $\omega(x_i) \neq 0$.

We then rescale the metric $g$ based at $x_i$, by setting

$$g_i = \rho(x_i)^{-2} \cdot g,$$

so that $\rho_{g_i}(x_i) = 1$. In the $g_i$ metric, the $\mathcal{Z}_c^2$ equations (4.7) take the form

$$\alpha_i \nabla \mathcal{Z}^2 + L^* \tau = 0, \tag{C.6}$$

$$\Delta \omega = -\frac{\alpha_i}{4}|z|^2,$$

where $\alpha_i = \alpha \rho_i^{-2}, \rho_i = \rho(x_i, g)$; all other metric quantities in (C.6) are w.r.t. $g_i$.

It follows from the results in §4.2 that a subsequence of $(N, g_i, x_i)$ converges, in the strong $L^{2,2}$ topology to a non-flat limit $(N_\infty, g_\infty, x_\infty)$ which is either a complete $\mathcal{Z}^2$ solution, a complete $\mathcal{Z}_c^2$ solution, or a complete static vacuum solution, possibly with junction to a $\mathcal{Z}$ solution. In all cases, the limit $(N_\infty, g_\infty)$ has a free isometric $S^1$ action. The limit potential $\omega_\infty$ is obtained as the limit of the renormalized potentials $\omega_i = \omega/|\omega(x_i)|$ as in §4. As usual, if the sequence $\{(N, g_i, x_i)\}$ collapses at $x_i$, then we pass to sufficiently large covers to unwrap the collapse.

As in the proof of Theorem C, the proof here must be separated into non-collapse and collapse cases.

**Case (I).(Non-Collapse).** Suppose there is a constant $a_o > 0$ such that, for $s \in (\frac{1}{2}t(x_i), 2t(x_i))$,

$$areaD(s) \geq a_o s^2, \tag{C.7}$$

in $(V, g_V)$, i.e. the base space $(V, (g_i)_V, \pi(x_i))$ is not collapsing.

Now by Proposition C.1,

$$\int_{V \setminus D(s)} |\nabla log f|^2 dA_V \to 0, \quad \text{as} \ \ s \to \infty. \tag{C.8}$$

The integral in (C.8) is scale-invariant, and also invariant under (multiplicative) renormalizations of $f$. Thus, we normalize $f$ at each $x_i$ so that

$$f(x_i) = 1. \tag{C.9}$$

(This is basically equivalent to passing to suitable covering or quotient spaces of $N$).

It follows first that the sequence $(N, g_i, x_i)$ is not collapsing at $x_i$ and so one has convergence to the limit $(N_\infty, g_\infty, x_\infty)$. Second, from (C.8) and the normalization (C.9), we obtain

$$\int_{N_\infty} |\nabla log f|^2 dV = 0, \tag{C.10}$$

and so $f$, the length of the $S^1$ orbits in $(N_\infty, g_\infty)$, is a constant function. Exactly the same reasoning on the second estimate in (C.3) implies that $A \equiv 0$ on the limit. Since $A = 0$ and $\nabla f = 0$, it follows that the limit $(N_\infty, g_\infty)$ is a product Riemannian manifold $N_\infty = V_\infty \times S^1$, with $s = s_{V_\infty} \geq 0$. In any region $U$ where $V_\infty$ is flat, we then have $U \times S^1 \subset N$ is also flat. Hence, from the limit on $(N_\infty, g_\infty)$ of the trace equation (C.6), we have $\Delta \tau_\infty = 0$ on $U \times S^1$. But this means exactly that $\Delta s \leq 0$ everywhere on $N_\infty$ and hence $\Delta_{V_\infty} s \leq 0$ everywhere on $V_\infty$. Thus $s$ is a non-negative superharmonic function on the complete surface $(V_\infty, g_{V_\infty})$ of non-negative Gauss curvature. It is well-known, c.f. [CY] for example, that such surfaces are parabolic, i.e. admit no non-constant positive superharmonic functions. It follows that $s = $ const, and hence $s = 0$, so that $(N_\infty, g_\infty)$ is flat. This contradiction proves that (C.4) holds in the non-collapse case.

**Case (II). (Collapse).** If (C.7) does not hold, so that

$$areaD(s) << s^2, \tag{C.11}$$

for some $s = s_i \in A(\frac{1}{2}t(x_i), 2t(x_i))$, then one needs to argue differently, since in this case, the estimate (C.8) may arise from collapse of the area, and not the behavior of $\log f$.

Since the base space $V$ is collapsing asymptotically by (C.11), the construction of the limit $(N_\infty, g_\infty)$ requires that this collapse be unwrapped by passing to larger and larger covers. This means that in addition to the original isometric $S^1$ action on $N$, we obtain a second free isometric $S^1$ action by unwrapping the collapse on $V$. Thus, one now has a free isometric $S^1 \times S^1$ action on $(N_\infty, g_\infty, \omega_\infty)$, so that $N_\infty$ is a torus bundle over $\mathbb{R}$.

However, $(N_\infty, g_\infty)$ is complete, and of non-negative scalar curvature. Again [SY3, Thm. 4] implies that any such metric is flat; (this can also be derived easily from a submersion formula as in (C.1) in the $S^1 \times S^1$ case). This contradiction then implies (C.4) must hold in the collapse case also.

Having proved (C.4), we now turn to the behavior of the potential $\omega$. First, the maximum principle applied to the trace equation (4.7) implies that $inf_{S_y(s)}\omega$ is monotone decreasing in $s$ as $s \to \infty$. We may assume that $inf_N \omega < 0$, since otherwise $(N, g)$ is a complete $\mathcal{Z}^2$ solution and the result follows from Theorem 5.1(I). Hence there is a curve $\gamma(s)$ in $N$ diverging uniformly to infinity in $N$ such that $\tau(\gamma(s)) \leq \tau_o < 0$, for all $s$ large.

Consider the collection of all tangent cones $T_\infty$ of $(N, g)$, (as defined following (7.18)), based on sequences $x_i$, with $\tau(x_i) \leq \tau_o$, for some $\tau_o < 0$. We claim that all such tangent cones at infinity are flat.

To see this, it follows from (C.4), exactly as in the proof of Lemma 7.5, that any such tangent cone at infinity $(T_\infty, g_\infty)$ of $(N, g)$ is a static vacuum solution.

If $(T_\infty, g_\infty)$ is the limit of a non-collapsing sequence $(N, g_i, x_i), g_i = t(x_i)^{-2} \cdot g, t(x_i) \to \infty$, then the same argument as in Case (I), using (C.8), proves that $(T_\infty, g_\infty)$ is a flat static vacuum solution.

If instead $(T_\infty, g_\infty)$ is the limit of a collapsing sequence, then as in Case (II), $(T_\infty, g_\infty)$ is a static vacuum solution with a free isometric $S^1 \times S^1$ action. As stated in the proof of Proposition 7.11, the only such solution which is not flat is the Kasner metric (7.48). But the arguments in (7.50)-(7.51) and Lemma 7.12 have already ruled out this possibility.

Hence, all such tangent cones at infinity of $(N, g)$ are flat products, with $f \equiv 1$, (under the normalization (C.9)). Observe that from the smooth convergence to $(T_\infty, g_\infty)$, we have

$$|\nabla \log f|(x) << 1/t(x), \tag{C.12}$$

as $t(x) \to \infty$, so that $f$ grows slower than any positive power of $t$. It follows from this and from (C.2) that the annuli $A_y(s, 2s)$ satisfy

$$diam A_y(s, 2s) \leq c \cdot s, vol A_y(s, 2s) \leq c \cdot s^3, \tag{C.13}$$

for some fixed constant $c < \infty$ and for all $s$, (as in (7.17)-(7.18)).

We see from the above that for any sequence $x_i$ as above with $t(x_i) \to \infty$, the manifolds $(N, g_i, x_i)$, $g_i = t(x_i)^{-2} \cdot g_i$ converge smoothly to a flat static vacuum solution, (unwrapping a collapse if necessary). Hence the limit potential $\bar{\tau}_\infty = \lim \bar{\tau}_i, \bar{\tau}_i = \tau(x)/|\tau(x_i)|$, is either constant or a non-constant affine function. The latter case is impossible since $\bar{\tau}_\infty \leq 0$ everywhere.

It follows that on any divergent sequence $\{x_i\}$ on which $\tau(x_i)$ is bounded away from 0,

$$\tau \to const, \tag{C.14}$$

uniformly on compact subsets of $(T_\infty, g_\infty, x_\infty)$; apriori, the constant may depend on the choice of $T_\infty$.

However, since $V = \mathbb{R}^2$, the geodesic annuli $D_{\pi(y)}(s, 2s) \subset V$ are connected, for all $s$ large, and hence so are the annuli $A_y(s, 2s)$ in $N$. It then follows from (C.13) and (C.14) that

$$osc_{A_y(s,2s)}\omega \to 0, \quad \text{as} \quad s \to \infty. \tag{C.15}$$

This, together with the use of the maximum principle as above implies that $\tau$ decreases uniformly to its infimum at infinity in $E$. In particular, $\omega$ is a proper exhaustion function satisfying (C.5). ∎

The proof of Theorem 5.3 now follows easily from Theorem C. Thus, Proposition C.2 shows that the hypotheses of Theorem C are satisfied for $(N, g)$. Hence, all ends of $(N, g)$ must be asymptotically flat. But $(N, g)$ is of the form $\mathbb{R}^2 \times S^1$, which does not have an asymptotically flat end. ∎

We conclude with the following result, also of some independent interest, which completes the last part of the proof of Theorem 3.12.

**Proposition C.3** *Let $(\Omega, g)$ be a $\mathcal{Z}^2$ solution, i.e. a solution of the Euler-Lagrange equations*

$$\nabla \mathcal{Z}^2 = \frac{1}{6} c \cdot g, \tag{C.16}$$

$$\frac{1}{2}\Delta s + |z|^2 = c, \tag{C.17}$$

*where $c = \int_\Omega |z|^2 dV_g$. If $(\Omega, g)$ is complete, with bounded curvature and of finite volume, and $s \geq 0$ on $\Omega$, then $\Omega$ must be a closed 3-manifold.*

*Proof.* We will prove that if $\Omega$ is non-compact, then necessarily $c = 0$. This then forces $(\Omega, g)$ to be of constant, hence non-negative, sectional curvature which is impossible since $(\Omega, g)$ is non-compact and of finite volume.

If $(\Omega, g)$ is complete, of finite volume and bounded curvature, then it collapses everywhere at infinity. Hence, for any divergent sequence $\{x_i\} \in \Omega$ and any $R < \infty$, the pointed sequence $(B_{x_i}(R), g, x_i)$ in $(\Omega, g)$ collapses with bounded curvature along an injective F-structure. As in §4, we may then unwrap this collapse by passing to sufficiently large finite covers; choosing a sequence $R_j \to \infty$, and a suitable diagonal subsequence gives rise to a limit solution $(N, g_\infty, x_\infty)$ of (C.16)-(C.17) with a free isometric $S^1$ action. The convergence to such limits is smooth, by elliptic regularity applied to these equations, (as in §4). The constant $c$ of course remains unchanged in passing to such limits. Hence, we may evaluate $c$ on these simpler manifolds $(N, g_\infty)$.

Now the manifold $(N, g_\infty)$ is exactly of the type analysed above in the proof of Theorem 5.3. With the same notation as before, if the orbit space $V$ of the $S^1$ action is non-collapsing at infinity, i.e. (C.7) holds, then the same arguments as in Case (I) above imply that $(N, g_\infty)$ is flat. If instead $V$ collapses at infinity, so (C.11) holds, then again the same arguments as in Case (II) above again give $(N, g_\infty)$ is flat. However, when $(N, g_\infty)$ is flat, the equation (C.17) immediately implies $c = 0$, as required. ∎

This result also implies that the potential $\tau_\varepsilon$ cannot vanish identically on any component of $(\Omega_\varepsilon, g_\varepsilon)$, $\varepsilon > 0$.

December, 1998/December, 2000

Department of Mathematics

S.U.N.Y. at Stony Brook

Stony Brook, N.Y. 11794-3651

anderson@math.sunysb.edu